\documentclass{amsart}
\usepackage{amscd}
\usepackage{amsfonts}
\usepackage{amsmath}
\usepackage{amssymb}
\usepackage{latexsym}
\usepackage{amsthm,color}
\usepackage[all]{xy}
\usepackage{epsfig}
\usepackage{graphicx}
\usepackage{color}
\usepackage{tikz}
\usetikzlibrary{cd}
\usetikzlibrary{calc}
\newtheorem{lemma}{Lemma}[section]

\newtheorem{theorem}{Theorem}

\theoremstyle{definition}
\newtheorem{definition}{Definition}
\theoremstyle{definition}
\newtheorem{example}{Example}[section]
\theoremstyle{definition}
\newtheorem{remark}{Remark}[section]

\newtheorem{fact}{Fact}[section]

\begin{document}
\title[
Hyperplane families creating envelopes]
{Hyperplane families creating envelopes}
\author[T.~Nishimura]{Takashi Nishimura
}
\address{
Research Institute of Environment and Information Sciences,  
Yokohama National University, 
240-8501 Yokohama, JAPAN}
\email{nishimura-takashi-yx@ynu.ac.jp}
\begin{abstract}
{\color{black}A simple geometric mechanism: 
\lq\lq the locus of intersections of perpendicular bisectors 
and normal lines\rq\rq,\,   
often arises in many guises in Nonlinear Sciences.   
In this paper, a new application of this simple geometric mechanism is given.  
Namely, we show that this mechanism gives answers 
to all four basic problems on envelopes created by 
hyperplane families (existence problem, representation problem, 
equivalence problem of definitions, uniqueness problem) at once.}
\end{abstract}
\subjclass[2010]{57R45, 58C25} 
\keywords{Hyperplane family, Frontal, Envelope, 
Creative, Creator, Mirror-image mapping, Anti-orthotomic, Orthotomic, 
Cahn-Hoffman vector formula.}


\date{}

\maketitle

\section{Introduction\label{section1}}
Throughout this paper, let $n$ be a positive integer.    
Moreover, all manifolds, functions and mappings are of class $C^\infty$ 
unless otherwise stated.   
\smallskip 

{\color{black}{
A simple geometric mechanism: 
\lq\lq the locus of intersections of perpendicular bisectors 
and normal lines\rq\rq,  
often arises in many guises in Physical Sciences.    
For example, as Richard Feynman elegantly explained in \cite{feynman}, 
the orbit of a planet around the sun can be understood 
as a consequence of this mechanism 
under the assumption of the inverse-square law (see Figure \ref{feynman1} 
where the circle is the hodograph of the velocity vectors of a planet
{\color{black}\!\!, }
that is to say, the circle is a curve drawn by the end points 
of the vectors that are parallel to the velocity vectors and start 
at a fixed point {\color{black}$P$}.   
The orbit of the planet is similar to the locus of intersections $B_t$ 
of the perpendicular bisectors of velocity vectors $\overrightarrow{PA_t}$ 
and the normal lines to the circle at $A_t$).  
This is an example in Celestial Mechanics.       
In the same book \cite{feynman}, one can find that even the historical discovery 
of atomic nucleus due to  Ernest Rutherford  can be explained 
as a consequence of this simple geometric  mechanism 
(see Figure \ref{feynman2} 
where  the center of circle $O$ is an atomic nucleus.    
The orbit of an $\alpha$ particle is the locus of intersections $B_t$ 
of the perpendicular bisectors of the segment $\overline{PA_t}$ 
and the normal lines to the circle at $A_t$
).     
This is an example in Nuclear Physics.   
\begin{figure}[htbp]
  \begin{minipage}[b]{0.48\linewidth}
    \centering
    \includegraphics[width=8cm]
{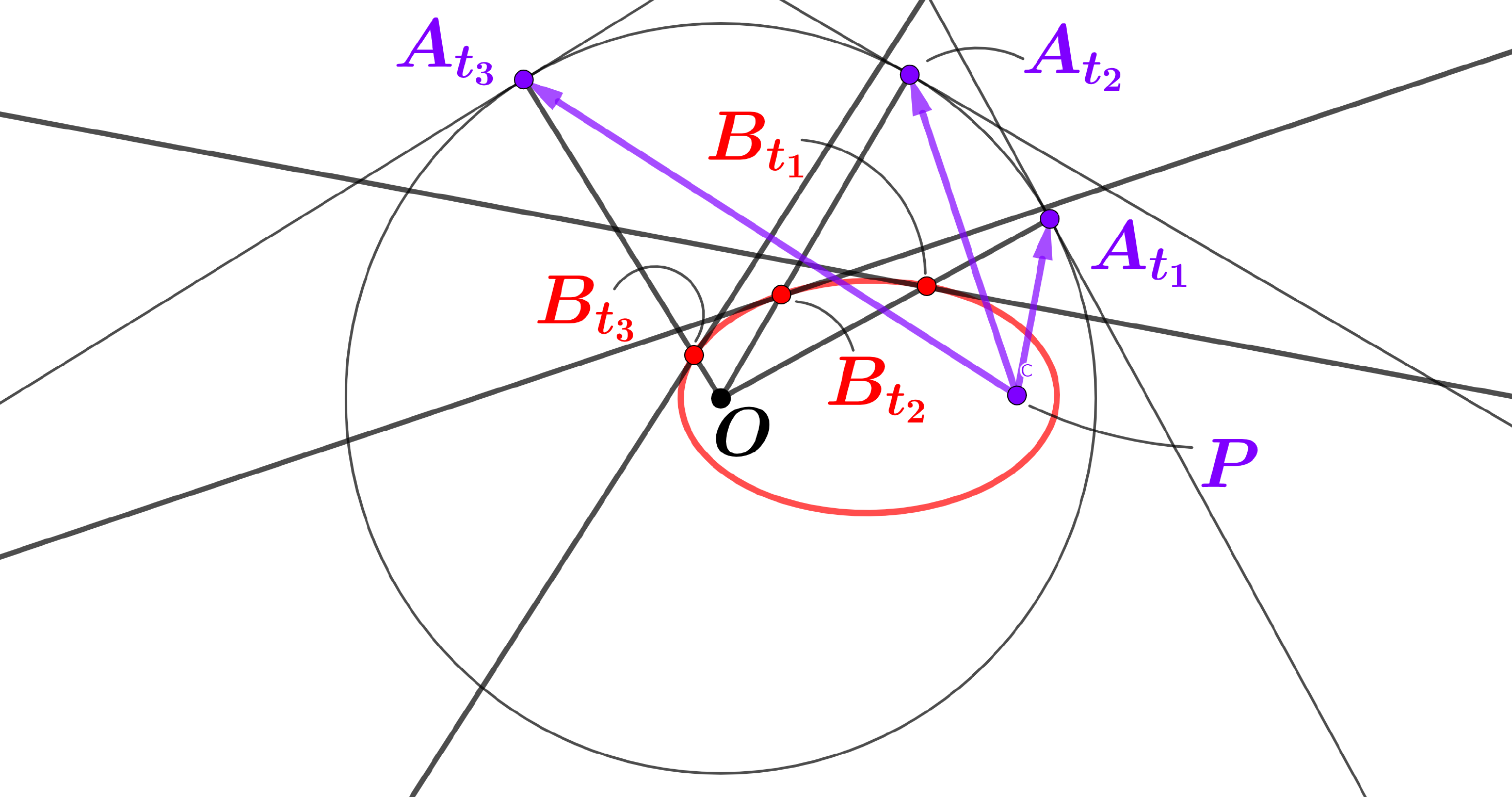}
    \caption{
Locus similar to the orbit of a planet.  
}\label{feynman1}
\end{minipage}
  \begin{minipage}[b]{0.48\linewidth}
    \centering
    \includegraphics[width=8cm]
{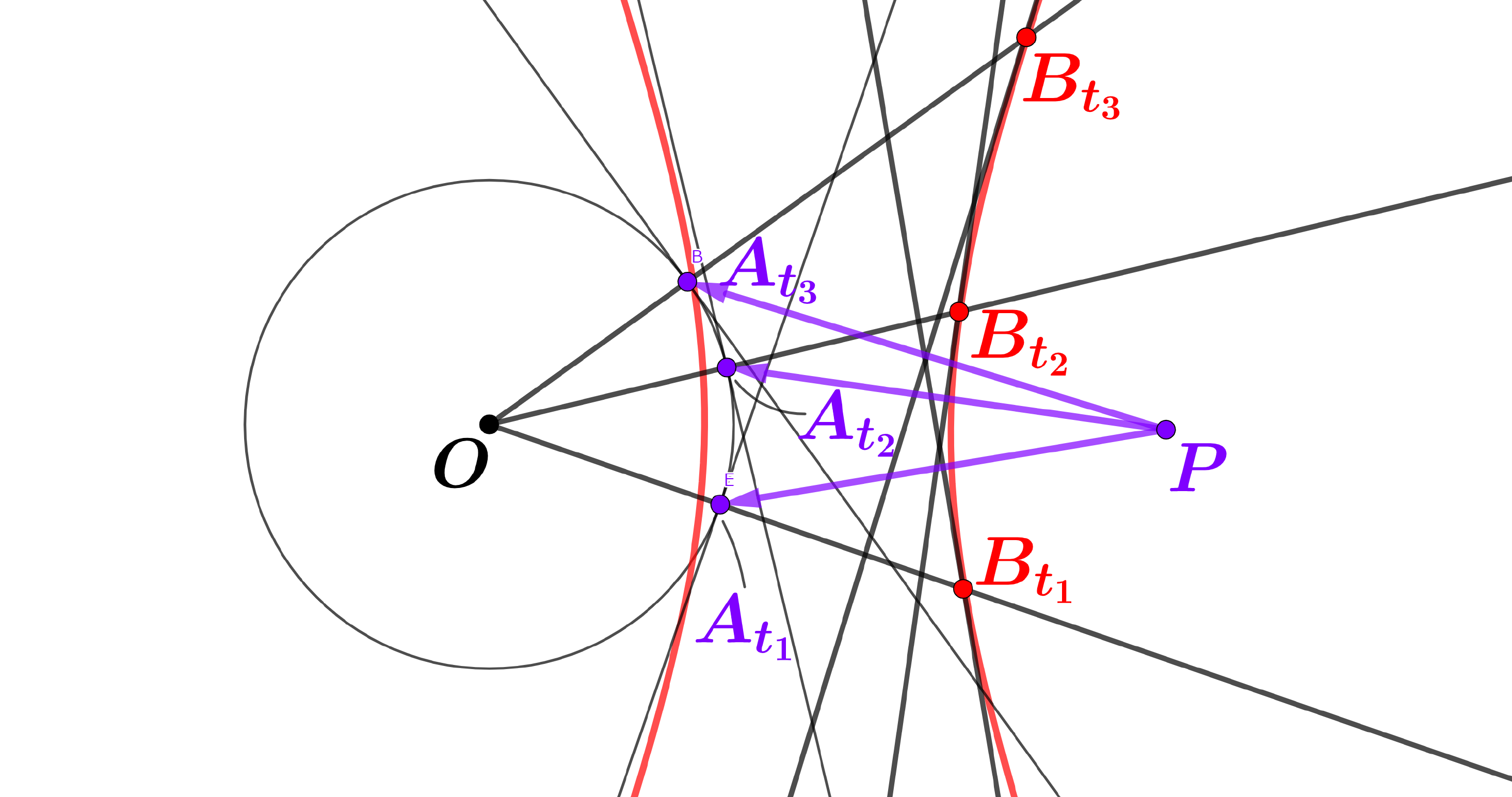}
    \caption{
Locus of an $\alpha$ particle.   
}\label{feynman2}
  \end{minipage}
\end{figure}
\par 
In Crystallography, one can find such the mechanism 
in the so-called \textit{Wulff construction} for the equilibrium shape of a crystal.      
A brief explanation of the Wulff construction is as follows.   
Given an equilibrium crystal, take an arbitrary point $P$ 
inside the crystal and fix it.   
Georg Wulff discovered in \cite{wulff} the so-called 
Gibbs-Wulff theorem which asserts that 
{\color{black}the length from the fixed point $P$ to the foot of the perpendicular 
to the tangent space to the face of the crystal} 
is proportional to its surface energy density of the face.   
Let $\gamma: S^2\to \mathbb{R}$ be the surface energy density function 
of the equilibrium crystal.   
The graph of $\gamma$ with respect to 
the polar coordinates about the point $P$
defines the mapping $g: S^2\to \mathbb{R}^3$.   
The mapping $g$ is often called the \textit{polar plot of $\gamma$} or 
the \textit{$\gamma$-plot} or the \textit{Wulff plot}.       
Set $f=2g$ and suppose that the image $f(S^2)$ has the well-defined 
normal vectors at any point $f(x)$.         
Then, by the Gibbs-Wulff theorem, 
the accurate shape of the crystal surface is proportional to the shape obtained 
by our simple geometric mechanism: \lq\lq the locus of the intersection 
of the perpendicular bisector of the vector $\overrightarrow{Pf(x)}$ 
and the normal line to $f(S^2)$ at $f(x)$\rq\rq.       
This is the Wulff construction and the constructed shape is called 
the \textit{Wulff shape}. 
Notice that in general $f$ is a continuous mapping 
and thus from the viewpoint of 
rigorous mathematics, 
the Wulff construction is not a well-defined construction method in general.      
Nevertheless, Hoffman and Cahn showed in \cite{hoffmancahn} that 
if $\gamma: S^2\to \mathbb{R}$ is {\color{black}differentiable}, then 
the image $f(S^2)$ has a well-defined normal vector at each point $f(x)$ and 
the set $\{\nabla\gamma(x)+\gamma(x)x\; |\; x\in S^2\}$ is exactly 
the shape obtained by our simple geometric mechanism for the point $P$ and 
the surface $f(S^2)$.  The Wulff construction 
and the Cahn-Hoffman formula in the plane is depicted in Figure \ref{wulff}.    
For details on the Wulff construction and Wulff shapes, see for instance 
\cite{giga, hannishimura}.      
\begin{figure}
\begin{center}
\includegraphics[width=10cm]
{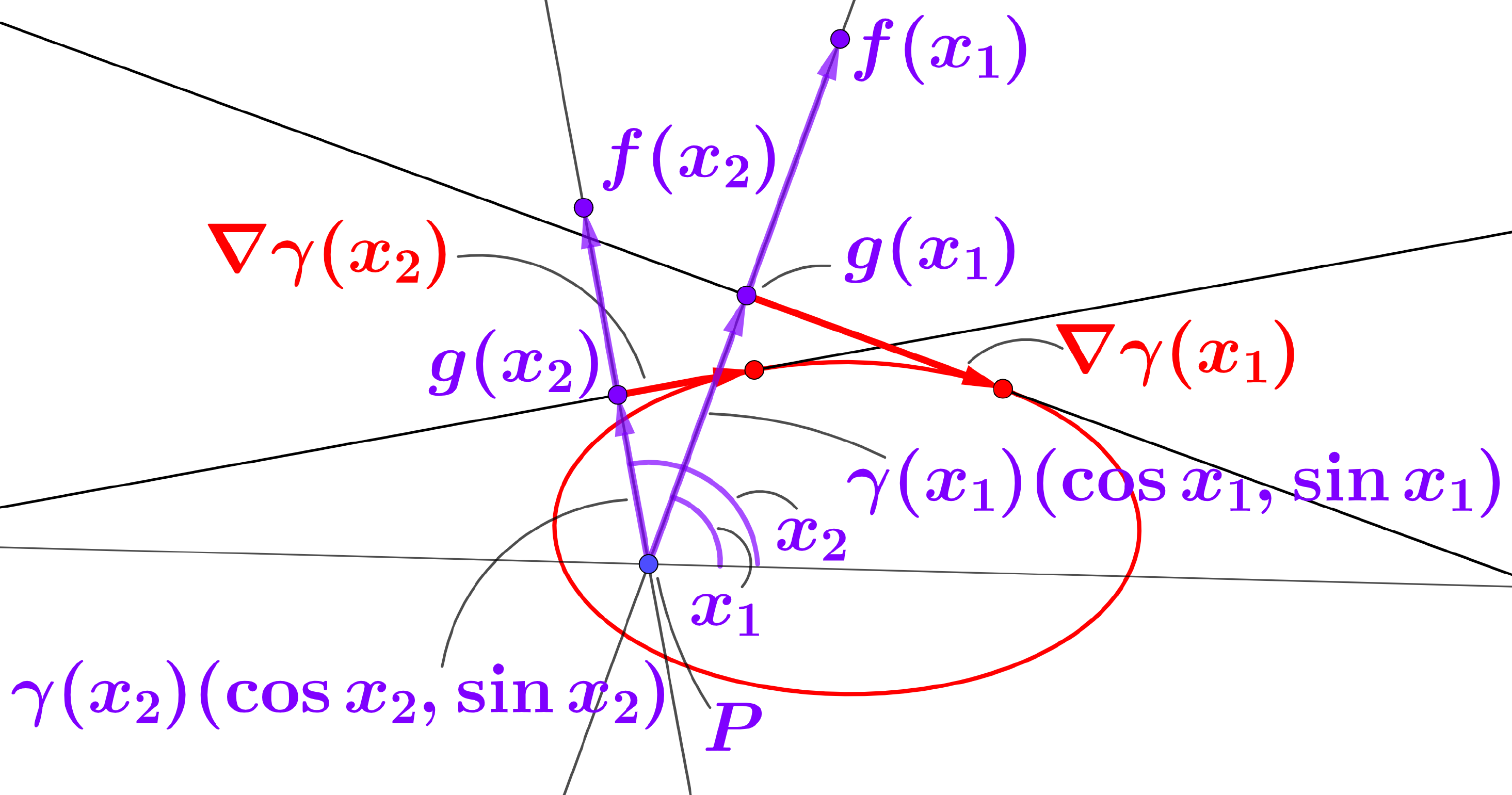}
\caption{The Wulff construction and the 
Cahn-Hoffman vector formula in the plane.    
} 
\label{wulff}
\end{center}
\end{figure}
\par 
Moreover, it is a surprising fact 
that 
our simple geometric mechanism: 
\lq\lq the locus of intersections of perpendicular bisectors 
and normal lines\rq\rq 
can be applied  even to Seismic Survey (see 7.14 (9) of \cite{brucegiblin}).   
\par 
\smallskip 
In Mathematics, our simple geometric mechanism:
\lq\lq the locus of intersections of perpendicular bisectors 
and normal lines\rq\rq 
is called the \textit{anti-orthotomic}  of a mapping $f$ having 
a well-defined normal vector to its image at each point  
(for details on anti-orthotomics, 
see 7.14 of \cite{brucegiblin}.   
See also \cite{janeczkonishimura} where anti-orthotomics are generalized to 
frontals and \cite{chaos} where more elementary explanation 
on anti-orthotomics 
can be found).       
In Mathematics as well, there are examples where anti-orthotomics 
are effectively applied (see \cite{brucegiblin}).    
\par 
\smallskip 
In order to understand better the powerfulness 
of the simple geometric mechanism, 
we would like to have more striking 
examples in Mathematics {\color{black}where} anti-orthotomics are 
effectively applied.   Namely, we want to seek 
mathematical problems 
which can be geometrically solved by our simple geometric mechanism  
though it seems difficult to solve them by other methods.    
This is the primitive motivation of this paper.   
In this paper, we show that the existence and uniqueness problem of envelopes 
for a given hyperplane family is one of such problems.    
Namely, we give a necessary and sufficient condition 
(see Definition \ref{definition2}) for a given hyperplane family  
to create an envelope.   And then, we give a necessary and sufficient 
condition for the uniqueness of created envelopes 
if the given hyperplane family creates an envelope.    
It seems difficult to prove that the condition 
given in Definition \ref{definition2} is actually a sufficient condition to create 
envelopes by other methods.    
In order to apply our simple geometric mechanism, 
we need some geometric objects 
to which the normal line can be reasonably well-defined at each point.   
Hyperplane families themselves are 
far from the reasonable geometric objects for our purpose.    
The reasonable geometric objects are frontals 
(the definition of frontal is given in the next paragraph).   
In order to obtain a frontal from a given hyperplane family, 
the mirror-image mapping will be locally introduced.    
Then, it turns out that if the given hyperplane family is creative 
(see Definition \ref{definition2} below), 
then the mirror-image mapping is actually a frontal such that the normal 
line at each point intersects the corresponding hyperplane.      
Thus, we can apply the anti-orthotomic method developed 
in \cite{janeczkonishimura} to obtain Theorem \ref{theorem1} 
and Theorem \ref{theorem2}.         
The existence and uniqueness 
problem of envelopes for a given hyperplane family 
can be easily interpreted as the existence 
and uniqueness problem of solutions for a 
certain type of system of first order differential equations 
with one constraint condition.    
In the author's opinion, one of the most attractive features 
of our simple geometric mechanism is 
that it can make all solutions and their precise 
expressions clear in one shot by geometry 
without the need to solve the corresponding system 
of differential equations with 
one constraint condition.  \\    
\par
\smallskip  
}}
Let $S^n$ be the $n$-dimensional unit sphere 
in the $(n+1)$-dimensional vector space 
$\mathbb{R}^{n+1}$.   Given a point $P$ of $\mathbb{R}^{n+1}$ 
and an $(n+1)$-dimensional 
unit vector ${\bf n}\in S^n\subset \mathbb{R}^{n+1}$, 
the \textit{hyperplane} $H_{\left(P, {\bf n}\right)}$ relative to $P$ and ${\bf n}$ 
is naturally defined as follows, 
where the dot in the center stands for the standard scalar product 
of two vectors $(X-P)$ and  ${\bf n}$ in the vector space $\mathbb{R}^{n+1}$.   
\[
H_{\left(P, {\bf n}\right)}=
\{X\in \mathbb{R}^{n+1}\; |\; (X-P)\cdot {\bf n} =0\}.   
\]
Let $N$ be an $n$-dimensional manifold without boundary.   
Given two mappings $\widetilde{\varphi}: N\to \mathbb{R}^{n+1}$ and 
$\widetilde{\nu}: N\to S^n$, 
the \emph{hyperplane family} 
$\mathcal{H}_{\left(\widetilde{\varphi}, \widetilde{\nu}\right)}$  
relative to $\widetilde{\varphi}$ and $\widetilde{\nu}$  
is naturally defined as follows.   
\[
\mathcal{H}_{\left(\widetilde{\varphi}, \widetilde{\nu}\right)}
=
{\color{black}\left\{H_{\left(\widetilde{\varphi}(x), \widetilde{\nu}(x)\right)}
\right\}_{x\in N}}.  
\] 
A mapping $\widetilde{f}: N\to \mathbb{R}^{n+1}$ is called a 
{\it frontal} if there exists a mapping 
$\widetilde{\nu}: N\to S^n$ such that 
$d\widetilde{f}_x({\bf v})\cdot \widetilde{\nu}(x)=0$ 
for any $x\in N$ and 
any ${\bf v}\in T_x N$,  
where two vector spaces 
$T_{\widetilde{f}(x)}\mathbb{R}^{n+1}$ and $\mathbb{R}^{n+1}$ are identified.   
By definition, 
it is natural to call 
$\widetilde{\nu}: N\to S^n$ 
a \textit{Gauss mapping} of the frontal $\widetilde{f}$.   
The notion of frontal has been {\color{black}recently} investigated (for instance, 
see \cite{ishikawa}).      
In this paper, as the definition of envelope created by a hyperplane family, 
the following is adopted.    
\begin{definition}\label{definition1}
Let $\mathcal{H}_{\left(\widetilde{\varphi}, \widetilde{\nu}\right)}$ 
be a hyperplane family.   
A mapping $\widetilde{f}: N\to \mathbb{R}^{n+1}$ is called 
an \textit{envelope created by} 
$\mathcal{H}_{\left(\widetilde{\varphi}, \widetilde{\nu}\right)}$ 
if the following two conditions are satisfied.   
\begin{enumerate}
\item[(a)] $\widetilde{f}(x)\in 
H_{\left(\widetilde{\varphi}(x), \widetilde{\nu}(x)\right)}$ 
for any $x\in N$.  
\item[(b)] $d\widetilde{f}_x({\bf v})\cdot \widetilde{\nu}(x)=0$ 
for any $x\in N$ and 
any ${\bf v}\in T_x N$.    
\end{enumerate}
\end{definition}
{\color{black}{
\noindent 
In other words, an {envelope created by} 
$\mathcal{H}_{\left(\widetilde{\varphi}, \widetilde{\nu}\right)}$ 
is a mapping $\widetilde{f}: N\to \mathbb{R}^{n+1}$ {\color{black}giving a solution} 
of the following system of first order differential equations with one constraint 
condition, where $\left(U, \left(x_1, \ldots, x_n\right)\right)$ 
is an arbitrary coordinate neighborhood of $N$ such that $x\in U$.      
\[
\left\{
\begin{array}{ccc}
\frac{\partial \widetilde{f}}{\partial x_1}(x)\cdot \widetilde{\nu}(x) & = & 0, \\ 
 \vdots & { } & { } \\
\frac{\partial \widetilde{f}}{\partial x_n}(x)\cdot \widetilde{\nu}(x) & = & 0, \\ 
\left(\widetilde{f}(x)-\widetilde{\varphi}(x)\right)\cdot \widetilde{\nu}(x) & = & 0.     
\end{array}
\right.
\]  
}}
\par 
By definition, any envelope $\widetilde{f}: N\to \mathbb{R}^{n+1}$ 
created by a hyperplane family 
$\mathcal{H}_{\left(\widetilde{\varphi}, \widetilde{\nu}\right)}$ 
must be a frontal with Gauss mapping $\widetilde{\nu}: N\to S^n$.     
For details on envelopes created by families of plane regular curves, 
refer to 
\cite{brucegiblin}.    
In Chapter 5 of \cite{brucegiblin}, 
several definitions for envelope are given.    
For a hyperplane family 
$\mathcal{H}_{\left(\widetilde{\varphi}, \widetilde{\nu}\right)}$, 
Definition \ref{definition1} is a generalization of their definition $E_2$ 
from a viewpoint of 
parametrization ({\color{black}$E_2$ envelope is a variety tangent to all 
lines of the given line family.   Thus, in the case of plane, 
an envelope defined by Definition \ref{definition1} is the same notion of 
$E_2$ envelope.   
For details on }the definition $E_2$, see 5.12 of \cite{brucegiblin}).      
The following definition, 
which may be regarded as a {\color{black}higher dimensional} 
generalization of $E_1$ 
from a viewpoint of parametrization 
({\color{black}$E_1$ envelope is the set of the limits of intersections 
with nearby members of the given line family. }  
For {\color{black}details on }the definition $E_1$, see 5.8 of \cite{brucegiblin} 
{\color{black}and for the relation  between Definition 2 in the plane case  
and $E_1$, see Subsection 2.3}), 
is the key notion for this paper.    
\begin{definition}\label{definition2}
{\rm 
Let $N$ be an $n$-dimensional manifold without boundary and let 
$\widetilde{\varphi}: N\to \mathbb{R}^{n+1}$, 
$\widetilde{\nu}: N\to S^n$ be mappings.   
Let $\widetilde{\gamma}: N\to \mathbb{R}$ be the function defined by 
$\widetilde{\gamma}(x)=\widetilde{\varphi}(x)\cdot \widetilde{\nu}(x)$.     
Let $T^*S^n$ be the cotangent bundle of $S^n$.   
A hyperplane family 
$\mathcal{H}_{\left(\widetilde{\varphi}, \widetilde{\nu}\right)}$ 
is said to be \textit{creative} if there exists a mapping 
$\widetilde{\Omega}: N\to T^*S^n$ with the form 
$\widetilde{\Omega}(x)
=\left(\widetilde{\nu}(x), \widetilde\omega(x)\right)$  such that 
for any $x_0\in N$ the equality     
$
d\widetilde{\gamma}=\widetilde\omega 
$ 
holds as germs of $1$-form at $x_0$.    
\[
\begin{tikzcd}
    & T^*S^n  \arrow[d]   \\
  N \arrow[ru, "{\widetilde{\Omega}}"] 
\arrow[r, "{\widetilde{\nu}}"]  & S^n 
\end{tikzcd}
\]
Namely,  $\mathcal{H}_{\left(\widetilde{\varphi}, \widetilde{\nu}\right)}$ 
is creative if there exists a $1$-form {\color{black}$\widetilde\Omega$} along 
$\widetilde{\nu}$ such that for any $x_0\in N$ 
by using of a coordinate neighborhood 
$\left(U, \left(x_1, \ldots, x_n\right)\right)$ of $N$ at $x_0$ and 
a \textit{normal} coordinate neighborhood 
$\left(V, \left(\Theta_1, \ldots, 
\Theta_n\right)\right)$ 
of  $S^n$ at  
$\widetilde{\nu}\left(x_0\right)$, 
the $1$-form germ 
$d\widetilde{\gamma}$ at $x_0$ 
is expressed 
as follows.   
\[
d\widetilde{\gamma} = 
\sum_{i=1}^n \left(\widetilde\omega(x)
\left({\color{black}\Pi}_{\left(\widetilde{\nu}(x), \widetilde{\nu}(x_0)\right)} 
\left(\frac{\partial}{\partial \Theta_i}\right) \right) \right)
d \left(\Theta_i\circ\widetilde{\nu}\right), 
\]
where {\color{black}
a normal coordinate neighborhood   
$\left(V, \left(\Theta_1, \ldots, 
\Theta_n\right)\right)$ is a 
local coordinate neighborhood at   
$\widetilde{\nu}\left(x_0\right)$ 
obtained by the inverse mappping of 
the exponential mapping at 
$\widetilde{\nu}\left(x_0\right)$,  
$S^n$ inherits its metric from the ambient space 
$\mathbb{R}^{n+1}$ and} 
${\color{black}\Pi}_{\left(\widetilde{\nu}(x), \widetilde{\nu}(x_0)\right)} : 
T_{\widetilde{\nu}\left(x_0\right)}S^n\to T_ {\widetilde{\nu}\left(x\right)}S^n$ 
is the Levi-Civita translation.      
{\color{black}Notice that our objective manifold 
is the unit sphere $S^n$ with metric inherited from 
$\mathbb{R}^{n+1}$.   Therefore, the Levi-Civita translation 
$\Pi_{\left(\widetilde{\nu}(x), \widetilde{\nu}(x_0)\right)}$ 
is the restriction of the rotation $R: \mathbb{R}^{n+1}\to \mathbb{R}^{n+1}$ 
satisfying $R(\widetilde{\nu}(x_0))=\widetilde{\nu}(x)$ 
to the tangent space $T_{\widetilde{\nu}\left(x_0\right)}S^n$.    
In particular, in the case $n=1$, a normal coordinate $\Theta$ at 
$\widetilde{\nu}\left(x\right)$ is nothing but the \textit{radian} 
(or, its negative) between two unit vectors $\widetilde{\nu}\left(x_0\right)$ and 
$\widetilde{\nu}\left(x\right)$ and the Levi-Civita translation 
$\Pi_{\left(\widetilde{\nu}(x), \widetilde{\nu}(x_0)\right)}$ is just 
the restriction of the plane rotation through $\Theta$ to the tangent space 
$T_{\widetilde{\nu}\left(x_0\right)}S^1$}.    
}
\end{definition} 
\begin{remark}\label{remark1}
\begin{enumerate}
{\color{black}{\item[(1)]}}\quad 
{\color{black}{
It is reasonable to say that $\widetilde{\gamma}$ is 
\textit{totally differentiable 
with respect to $\widetilde{\nu}$} if 
$\mathcal{H}_{\left(\widetilde{\varphi}, \widetilde{\nu}\right)}$ 
is creative.   
}}
%
\item[(2)]\quad 
For a creative hyperplane family 
$\mathcal{H}_{\left(\widetilde{\varphi}, \widetilde{\nu}\right)}$, 
the map-germ 
$\left(\widetilde{\nu}, \widetilde{\gamma}\right) 
: (N, x_0)\to S^n\times\mathbb{R}$ at any $x_0\in N$ is called 
an \textit{opening} of $\widetilde{\nu}: (N, x_0)\to S^n$ 
(for opening germs, see for example \cite{ishikawaaustralia}).     Thus, 
Definition \ref{definition2} may be regarded 
as a globalization of 
the notion of opening.     
\item[(3)]\quad 
Definition \ref{definition2} may be interpreted as follows.   
Let $\theta$ be a canonical contact $1$-form on 
$J^1(S^n, \mathbb{R})$, namely  
at any $(X_0, Y_0, P_0)\in J^1\left(S^n, \mathbb{R}\right)$ the $1$-form germ  
$\theta$ is expressed as $\theta=d Y - \sum_{i=1}^n C_i d \Theta_i$,  
where $\left(\Theta_1, \ldots, \Theta_n\right)$ 
is a \textit{normal} coordinate system at $X_0$ and 
$\left(\Theta_1, \ldots, \Theta_n, Y, C_1, \ldots, C_n\right)$ 
is a canonical coordinate system of $J^1\left(S^n, \mathbb{R}\right)$ at 
$(X_0, Y_0, P_0)$.    
Then, a hyperplane family 
$\mathcal{H}_{\left(\widetilde{\varphi}, \widetilde{\nu}\right)}$ 
is creative if there exists a mapping 
${\Omega}: N\to J^1\left(S^n, \mathbb{R}\right)$ with the form 
${\Omega}(x)=
\left(\widetilde{\nu}(x), \widetilde{\gamma}(x), \widetilde{c}_1(x), 
\ldots, \widetilde{c}_n(x)\right)$  
such that ${\Omega}^*\theta=0$, where 
$\widetilde{c}_1, \ldots, \widetilde{c}_n: N\to \mathbb{R}$ are some functions.    
\[
\begin{tikzcd}
    & J^1\left(S^n, \mathbb{R}\right) \arrow[d]   \\
  N \arrow[ru, "{\Omega}"] 
\arrow[r, "{\widetilde{\nu}}"]  & S^n 
\end{tikzcd}
\]
Notice that in Legendrian Singularity Theory, 
at any point $x_0\in N$, the map-germ 
${\Omega}: \left(N, x_0\right)\to J^1\left(S^n, \mathbb{R}\right)$ 
is assumed to be immersive and 
it is called a \textit{Legendrian immersion};  
and for Legendrian immersion ${\Omega}$, the mapping 
$N\ni x\mapsto \left(\widetilde{\nu}(x), \widetilde{\gamma}(x)\right)$ is called a 
\textit{wavefront} or \textit{front}  
(for details on Legendrian Singularity Theory and fronts, 
see for instance \cite{arnold, arnoldetall, sajiumeharayamada}).     
On the other hand, in Definition \ref{definition2}, 
${\Omega}$ is not assumed to be immersive in general and the mapping  
${\Omega}$ is called a \textit{Legendrian mapping}  
(for details on Legendrian mappings, 
see for instance \cite{ishikawaaustralia, ishikawa, takahashi}).       
Thus, in Definition \ref{definition2}, in general, 
the set-germ $\left({\Omega}(N), {\Omega}\left(x_0\right)\right)$ 
may be singular at some point 
$x_0\in N$     
(for example,  see Example \ref{example1}(4)).  
\item[(4)]\quad 
Notice that the $1$-form {\color{black}$\widetilde\Omega$} 
along $\widetilde{\nu}$ 
in Definition 
\ref{definition2} is not necessarily the pullback of a $1$-form over $S^n$ 
by $\widetilde{\nu}$ (for example,  see Example \ref{example1}(3), (4)) 
and 
{\color{black}the \lq\lq creativeness\rq\rq\, does not depend on the
particular choice of $\widetilde{\varphi}, \widetilde{\nu}$ 
and depends only on the hyperplane family 
$\mathcal{H}_{\left(\widetilde{\varphi}, \widetilde{\nu}\right)}$.}    
In the case that $N=S^n$ and $\widetilde{\nu}:S^n\to S^n$ 
is the identity mapping, 
for any $\widetilde{\varphi}: S^n\to \mathbb{R}^{n+1}$ the hyperplane family 
$\mathcal{H}_{\left(\widetilde{\varphi}, \widetilde{\nu}\right)}$ 
is always creative  by the following 
equality.   
\[
d\widetilde{\gamma}
= 
\sum_{i=1}^n 
\frac{\partial \widetilde{\gamma}
}{\partial \Theta_i}
d\Theta_i.   
\]
More generally, if $\widetilde{\gamma}: U\to \mathbb{R}$ 
may be expressed as the composition 
of $\widetilde{\nu}: U\to S^n$ and 
a certain function $\xi: S^n\to \mathbb{R}$ over an open set 
$U\subset N$, then  
%
the hyperplane family 
$\mathcal{H}_{\left(\widetilde{\varphi}|_U, \widetilde{\nu}|_U\right)}$ is creative.    
However, there are examples showing that {\color{black}there does not exist 
a function $\widetilde{\alpha}: S^n\to \mathbb{R}$ 
such that $\widetilde{\gamma}=\widetilde{\alpha}\circ \widetilde{\nu}$} 
although $\mathcal{H}_{\left(\widetilde{\varphi}, \widetilde{\nu}\right)}$ 
is creative {\color{black}(for example,  see Example \ref{example1}(3), (4))}.   
Moreover, there are many examples such that 
$\mathcal{H}_{\left(\widetilde{\varphi}|_U, \widetilde{\nu}|_U\right)}$ 
is not creative.  For instance, 
for any constant mapping $\widetilde{\nu}: \mathbb{R}\to S^1$, 
the line family $\mathcal{H}_{\left(\widetilde{\varphi}, \widetilde{\nu}\right)}$ 
is not creative  where $\widetilde{\varphi}: \mathbb{R}\to \mathbb{R}^2$ 
is defined by 
$\widetilde{\varphi}(t)=t^2\widetilde{\nu}(t)$.    And, it is clear in this case that 
 $\mathcal{H}_{\left(\widetilde{\varphi}, \widetilde{\nu}\right)}$ 
does not create an envelope 
in the sense of Definition \ref{definition1}.    However, it is easily seen that 
\begin{eqnarray*}
\qquad 
\qquad \quad \mathcal{D} & = & \left\{\left(X_1, X_2\right)\in \mathbb{R}^2\, 
|\, \exists t \mbox{ s.t. }
F\left(X_1, X_2, t\right)=
\frac{\partial F}{\partial t}\left(X_1, X_2, t\right)=0\right\}  \\ 
\qquad\quad { } & = & 
\left\{(X_1, X_2)\in 
\mathbb{R}^2\, |\, \left(X_1, X_2\right)\cdot \widetilde{\nu}(0)=0\right\}
\ne \emptyset,  
\end{eqnarray*} 
where $F\left(X_1, X_2, t\right)=
\left(\left(X_1, X_2\right)-\widetilde{\varphi}(t)\right)\cdot \widetilde{\nu}(t)$.      
Thus, for this example, the envelope defined by Definition \ref{definition1} 
is different from the envelope in the sense of classical definition 
(see 5.3 of \cite{brucegiblin}), 
For more examples on creative/non-creative hyperplane families 
and on comparison of 
Definition \ref{definition2} with the classical envelope $\mathcal{D}$, 
see Section 4.       
Therefore, it seems that  
the current situation on both the definitions of envelope 
and the relation of the creative condition (Defnition \ref{definition2}) 
with an envelope {\color{black}{seems 
to be wrapped in mystery}}.   
\end{enumerate} 
\end{remark}   
\par 
\medskip 
By definition, any frontal $\widetilde{f}: N\to \mathbb{R}^{n+1}$ with Gauss 
mapping $\widetilde{\nu}: N\to S^n$ is an envelope created by 
$\mathcal{H}_{(\widetilde{f}, \widetilde{\nu})}$.     
Therefore, the notion of envelope created by a hyperplane family 
is the same as the notion of frontal.   
Moreover, it is clear that for any mapping $\widetilde{\nu}: N\to S^n$, 
a constant mapping $\widetilde{f}: N\to \mathbb{R}^{n+1}$ is 
an envelope created by 
$\mathcal{H}_{\left(\widetilde{f}, \widetilde{\nu}\right)}$.    
On the other hand, for a constant mapping $\widetilde{\nu}: \mathbb{R}\to S^1$, 
if the line family $\mathcal{H}_{\left(\widetilde{\varphi}, \widetilde{\nu}\right)}$ 
does not create an envelope then 
$\widetilde{\varphi}: \mathbb{R}\to \mathbb{R}^2$ must be 
not constant.    
From these elementary observations, it is natural to ask 
to obtain 
a necessary and sufficient condition for a given hyperplane family 
$\mathcal{H}_{\left(\widetilde{\varphi}, \widetilde{\nu}\right)}$ 
to create an envelope 
$\widetilde{f}: N\to \mathbb{R}^{n+1}$ 
in terms of {\color{black}$\widetilde{\gamma}: N\to \mathbb{R}$} 
and $\widetilde{\nu}: N\to S^n$.  
{\color{black}Moreover, it is also desirable to solve the following two incidentally.     
\lq\lq Suppose that a given hyperplane family 
$\mathcal{H}_{\left(\widetilde{\varphi}, \widetilde{\nu}\right)}$ 
creates an envelope $\widetilde{f}: N\to \mathbb{R}^{n+1}$.   Then, 
obtain a representation formula of $\widetilde{f}$.\rq\rq \,    
\lq\lq Suppose that $n=1$.    Then, find the precise relation between 
$E_1$ envelope and $E_2$ envelope.\rq\rq         
}    
In this paper, {\color{black}as an application of our simple geometric mechanism,} 
{\color{black}all of these problems are} solved as follows.  
\begin{theorem}\label{theorem1}
Let $N$ be an $n$-dimensional manifold without boundary and let 
$\widetilde{\varphi}: N\to \mathbb{R}^{n+1}$, 
$\widetilde{\nu}: N\to S^n$ be mappings.   
Then, {\color{black}the following three hold.   
\begin{enumerate}
\item[(1)]\quad The hyperplane family 
$\mathcal{H}_{\left(\widetilde{\varphi}, \widetilde{\nu}\right)}$ creates 
an envelope if and only if it is creative. 
\item[(2)]\quad 
Suppose that the hyperplane family 
$\mathcal{H}_{\left(\widetilde{\varphi}, \widetilde{\nu}\right)}$ creates 
an envelope $\widetilde{f}: N\to \mathbb{R}^{n+1}$.    Then, for 
any $x\in N$, 
under the canonical identifications 
$ 
T^*_{\widetilde{\nu}(x)}S^n\cong T_{\widetilde{\nu}(x)}S^n 
\subset T_{\widetilde{\nu}(x)}\mathbb{R}^{n+1} \cong \mathbb{R}^{n+1}$, 
the $(n+1)$-dimensional vector $\widetilde{f}(x)$ is represented as follows.   
\[
\widetilde{f}(x)=\widetilde{\omega}(x)+\widetilde{\gamma}(x)\widetilde{\nu}(x),   
\]  
where the $(n+1)$-dimensional vector $\widetilde{\omega}(x)$ is identified 
with the corresponding $n$-dimensional 
cotangent vector $\widetilde{\omega}(x)$ under these identifications.   
\item[(3)]\quad 
Suppose that $n=1$.   Then, the line family 
$\mathcal{H}_{\left(\widetilde{\varphi}, \widetilde{\nu}\right)}$ creates 
an envelope ($E_2$-envelope) if and only if it creates an $E_1$ envelope.   
Moreover, these two envelopes are exactly the same.    
\end{enumerate}
}
\end{theorem}
\noindent 
By Theorem \ref{theorem1}, it is natural to call 
$\widetilde\omega$  
the \textit{creator} for an envelope $\widetilde{f}$ created by  
$\mathcal{H}_{\left(\widetilde{\varphi}, \widetilde{\nu}\right)}$.  
{\color{black}Recall that $E_1$ envelope (resp., $E_2$ envelope) 
is the set of the limit of intersections with nearby lines  
(resp., a parametrization tangent to all members of the given family).     
Thus, even in the case of plane, 
$E_2$ envelope is exactly the same as the envelope in Definition 
\ref{definition1}.  
}   
\par 
\medskip 
The key idea for the proof of Theorem \ref{theorem1}   
is to regard the given hyperplane family 
as a moving mirror parametrized by $x\in N$.     
Then, for any parameter $x_0\in N$, by taking a point   
$P\in \mathbb{R}^{n+1}$ outside 
the mirror 
$H_{\left(\widetilde{\varphi}(x_0), \widetilde{\nu}\left(x_0\right)\right)}$, 
the mirror-image  
\[
f_{{}_P}(x)=2\left(\left(\widetilde{\varphi}(x)-P\right)
\cdot \widetilde{\nu}(x)\right)\widetilde{\nu}(x)+P 
\]  
of $P$ by the mirror 
$H_{\left(\widetilde{\varphi}(x), \widetilde{\nu}(x)\right)}$ must have the same  
information as the mirror {\color{black}since the mirror is reconstructed}  
as the perpendicular {\color{black}bisectors} 
of the segment $\overline{Pf_{{}_P}(x)}$, 
where $x$ is a point in a sufficiently small neighborhood $U_{{}_P}$ of $x_0$.    
\begin{figure}
\begin{center}
\includegraphics[width=15cm]
{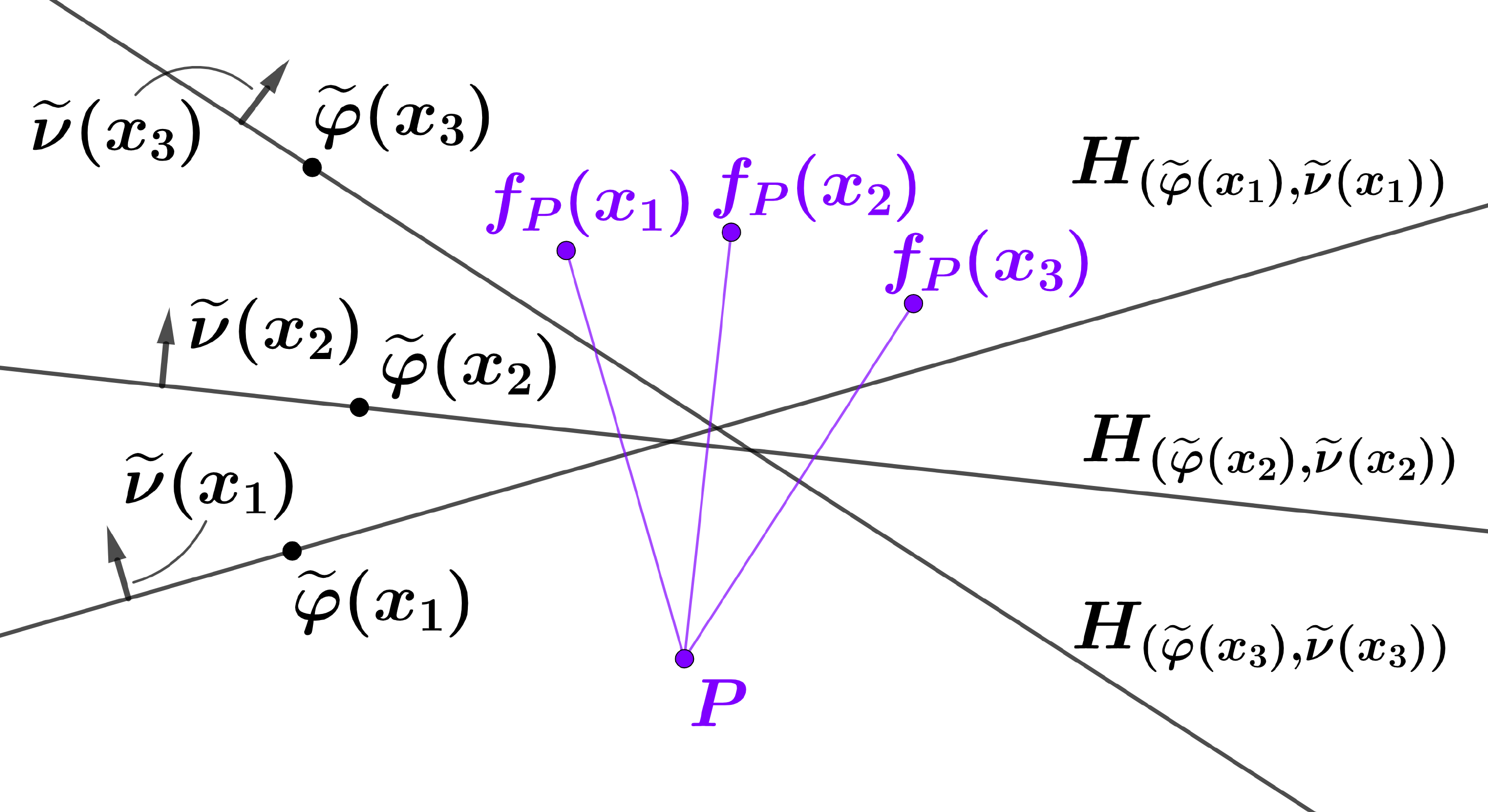}
\caption{
The mirror-image mapping $f_{{}_P}: U_{{}_P}\to \mathbb{R}^{n+1}$.    
} 
\label{figure1}
\end{center}
\end{figure}
Hence, investigation of the given hyperplane family 
$\mathcal{H}_{\left(\widetilde{\varphi}|_{{}U_{{}P}}, \widetilde{\nu}|_{{}U_{{}P}}\right)}$
may be replaced with   
analyzing the {\color{black}associated} \textit{mirror-image mapping}  
$f_{{}_P}: U_{{}_P}\to \mathbb{R}^{n+1}$ (see Figure \ref{figure1}).      
This suggests applicability of results  
in \cite{janeczkonishimura} to the problem {\color{black}{of this paper}}.      
\par
\begin{figure}
\begin{center}
\includegraphics[width=15cm]
{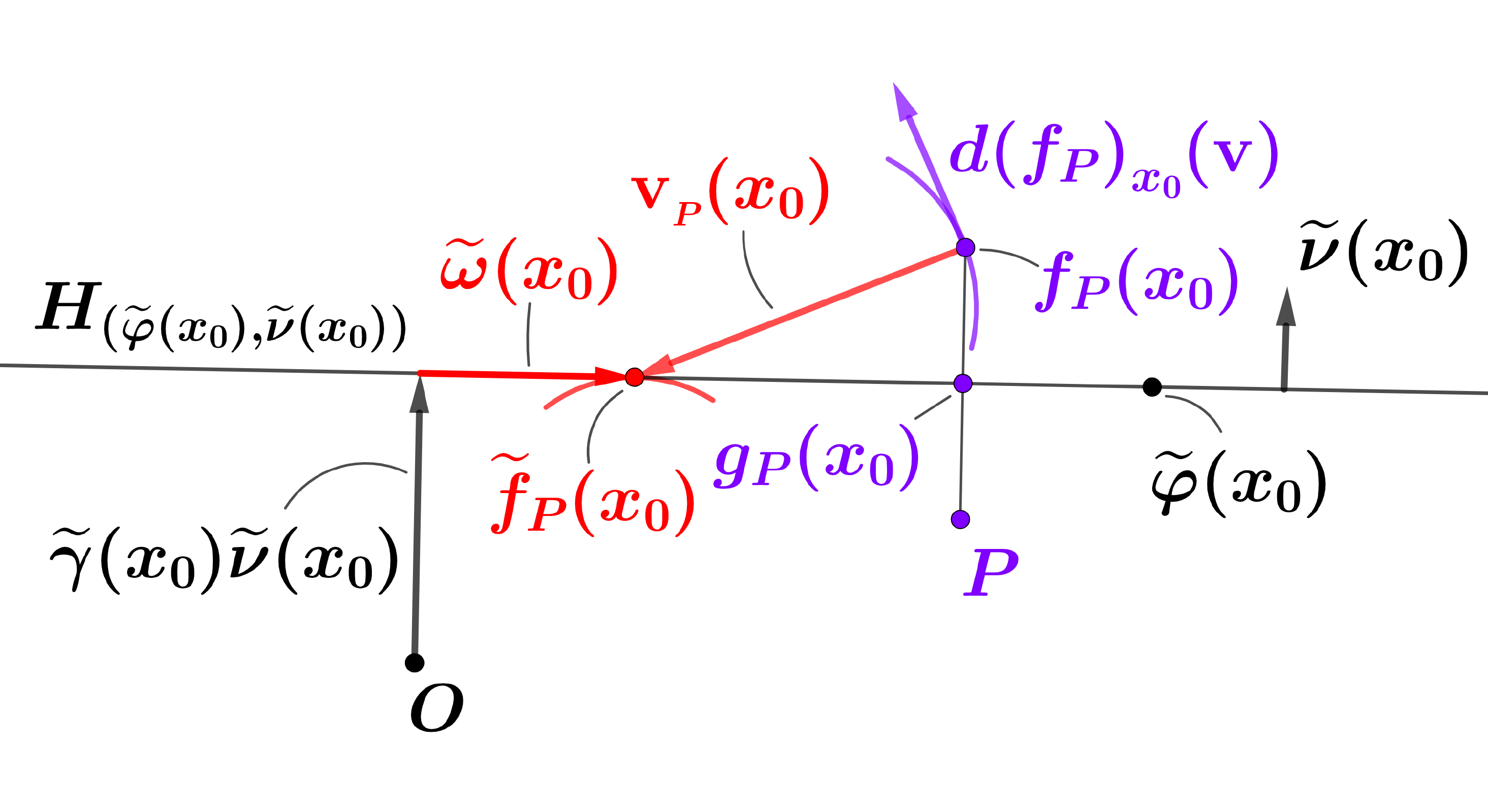}
\caption{The location 
$\widetilde{f}_{{}_P}(x_0)$ does not depend on 
the particular choice of $P$.    
} 
\label{figure3}
\end{center}
\end{figure}

 
A sketch of the proof of Theorem \ref{theorem1} {\color{black}(1)} 
may be given as follows.    
Suppose that the hyperplane family  
$\mathcal{H}_{\left(\widetilde{\varphi}, \widetilde{\nu}\right)}$ is creative.  
Then, by definition, there exists a mapping  
$\widetilde{\Omega} : N\to T^*S^n$ having the form 
$\widetilde{\Omega}(x)=\left(\widetilde{\nu}(x), \widetilde\omega(x)\right)$ 
such that 
the equality 
$
d\widetilde{\gamma}
= 
\widetilde\omega
$ 
holds as germs of $1$-form 
at $x_0$.     
Then, by investigating the Jacobian matrix of the mirror-image mapping 
$f_{{}_P}: U_{{}_P}\to \mathbb{R}^{n+1}$ at $x\in U_{{}_P}$ directly,  
it turns out that for any $x\in U_{{}_P}$ the non-zero vector   
\[
{\bf v}_{{}_P}(x)=
\sum_{i=1}^n \left(\left(\widetilde\omega\left(x\right)-P\right)
\left(\frac{\partial}
{\partial \Theta_{\left(i, \widetilde{\nu}(x)\right)}}\right) \right)
\frac{\partial}{\partial \Theta_{\left(i, \widetilde{\nu}(x)\right)}} 
- 
\left(\left(\widetilde{\varphi}(x)-P\right)\cdot 
\widetilde{\nu}\left(x\right)\right)\widetilde{\nu}\left(x\right)
\]
is perpendicular to the vector 
$ d \left(f_{{}_P}\right)_{x}\left({\bf v}\right)$ for any 
${\bf v}\in T_{x} N$,  
where $\mathbb{R}^{n+1}$, 
$T_{\widetilde{\nu}\left(x\right)}\mathbb{R}^{n+1}$ and 
$T^*_{\widetilde{\nu}\left(x\right)}\mathbb{R}^{n+1}$ are identified 
and $\frac{\partial}{\partial \Theta_{\left(i, \widetilde{\nu}(x)\right)}}=
P_{\left(\widetilde{\nu}(x), \widetilde{\nu}\left(x_0\right)\right)}
\left(\frac{\partial}{\partial \Theta_i}\right)$.   
Thus, $f_{{}_P}: U_{{}_P}\to \mathbb{R}^{n+1}$ is a frontal.    
From the construction,  
the mapping $\widetilde{f}_{{}_P}={\bf v}_{{}_P}+f_{{}_P}:U_{{}_P}\to \mathbb{R}^{n+1}$ 
must be exactly the same as the mapping $\widetilde{f}_{{}_P}$ 
given in Theorem 1 of \cite{janeczkonishimura}.    
Therefore, by Theorem 1 of \cite{janeczkonishimura} 
{\color{black}asserting that $\widetilde{f}_{{}_P}$ satisfies both conditions (a), (b) of 
Definition \ref{definition1}}, 
$\widetilde{f}_{{}_P}$ is an envelope created by the hyperplane family 
$\mathcal{H}_{\left(\widetilde{\varphi}|_{U_{{}_P}}, \widetilde{\nu}|_{U_{{}_P}}\right)}$.      
The mapping $\widetilde{f}_{{}_P}: U_{{}_P}\to \mathbb{R}^{n+1}$ is called 
the \textit{anti-orthotomic} of $f_{{}_P}: U_{{}_P}\to \mathbb{R}^{n+1}$  
relative to $P$.   
Calculation shows   
\[
\widetilde{f}_{{}_P}(x_0)  
 =  
\widetilde\omega\left(x_0\right)
+
\widetilde{\gamma}\left(x_0\right)
\widetilde{\nu}\left(x_0\right).  
\leqno{(*)}
\] 
Thus, unlike $f_{{}_P}(x_0)$,  
the location $\widetilde{f}_{{}_P}(x_0)$ does not depend on 
the particular choice of $P$.        
In other words, in order to discover the formula $(*)$, 
the role of $P$ is merely an auxiliary point 
just like an auxiliary line in elementary geometry (see Figure \ref{figure3}).     
Since $x_0$ is an arbitrary point of $N$, the hyperplane family  
$\mathcal{H}_{\left(\widetilde{\varphi}, \widetilde{\nu}\right)}$ 
creates an envelope 
$\widetilde{f}: N\to \mathbb{R}^{n+1}$.   
\par 
\noindent 
Conversely, suppose that the given hyperplane family 
$\mathcal{H}_{\left(\widetilde{\varphi}, \widetilde{\nu}\right)}$ 
creates an envelope 
$\widetilde{f}: N\to \mathbb{R}^{n+1}$.   Then, 
the mirror-image mapping $f_{{}_P}: U_{{}_P}\to \mathbb{R}^{n+1}$  
(resp., the mapping $g_{{}_P}: U_{{}_P}\to \mathbb{R}^{n+1}$ defined by 
$g_{{}_P}(x)=\left(\widetilde{f}(x)-P\right)\cdot \widetilde{\nu}(x)+P$)  
is called the \textit{orthotomic} (resp., \textit{pedal}) of $\widetilde{f}|_{U_{{}_P}}$ 
relative to the point $P$.    
It is known that both the orthotomic $f_{{}_P}$ and the pedal $g_{{}_P}$  
are frontals (see Proposition 1 and Corollary 1 of \cite{janeczkonishimura}).   
We prefer to investigate the orthotomic  $f_{{}_P}$ rather than the pedal 
$g_{{}_P}$ because its Gauss mapping 
$\nu_{{}_P}: U_{{}_P}\to S^n$ has characteristic properties:   
$\nu_{{}_P}(x)=\frac{\widetilde{f}(x)-f_{{}_P}(x)}{||\widetilde{f}(x)-f_{{}_P}(x)||}$ 
and $\widetilde{\nu}(x)\cdot\nu_{{}_P}(x)\ne 0$ for any $x\in U_{{}_P}$, and thus 
we can take a bird's eye view of $\widetilde{f}(x)$.    
Set $\widetilde\omega(x)=\widetilde{f}(x)-
\widetilde{\gamma}(x)\widetilde{\nu}(x)$ 
and $\widetilde{\Omega}(x)=\left(\widetilde{\nu}(x),\widetilde \omega(x)\right)$ 
for any $x\in U_{{}_P}$.    
Then, under the identification of $\mathbb{R}^{n+1}$ and 
$T^*_{\widetilde{\nu}(x)}\mathbb{R}^{n+1}$, 
$\widetilde{\Omega}$ having the form 
$\widetilde{\Omega}(x)=\left(\widetilde{\nu}(x), \widetilde{\omega}(x)\right)$ 
is a well-defined mapping 
$U_{{}_P}\to T^*S^n$.     
By investigating the Jacobian matrix of the mirror image mapping 
$f_{{}_P}$ at $x\in U_{{}_P}$ directly again, it turns out that 
$\widetilde\omega$ is actually 
the creator for the envelope $\widetilde{f}|_{U_{{}_P}}$.   
Since the vector $\widetilde\omega(x_0)$ does not depend 
on the particular choice of $P$ and the point $x_0$ is an arbitrary point of $N$, 
$\mathcal{H}_{\left(\widetilde{\varphi}, \widetilde{\nu}\right)}$ is creative.   
\par 
{\color{black}Theorem \ref{theorem1} (2) is 
a direct by-product of the proof of Theorem \ref{theorem1} (1)  
(see Figure \ref{figure3}).   
Theorem \ref{theorem1} (3) seems to be not a direct by-product 
of the proof of Theorem \ref{theorem1} (1) 
although 
it can be proved relatively easily 
by using the above argument (see Subsection \ref{proof_1(3)}).   }      
\par 
\medskip 
When $N=S^n$ and $\widetilde{\nu}: S^n\to S^n$ is the identity mapping, 
{\color{black}{
it is easily seen 
$\widetilde{\omega}(x)=\nabla\widetilde{\gamma}(x)$.   
Therefore, in the case that 
$N=S^n$ and $\widetilde{\nu}: S^n\to S^n$ is the identity mapping, }}
{\color{black}Theorem \ref{theorem1} (2)}  
has been known 
as the Cahn-Hoffman vector formula (\cite{hoffmancahn}).   
{\color{black}Theorem \ref{theorem1} (2)} is 
a {\color{black}comprehensive} generalization of their formula.   
{\color{black}{
Any Wulff shape is clearly a convex body and conversely it is known that 
any convex body can be constructed by the Wulff construction 
(for instance, see \cite{taylor}).    
There are many Wulff shapes such that 
the surface energy density functions $\gamma: S^n\to \mathbb{R}$ 
are not {\color{black}differentiable} 
(convex polytopes are typical examples).   
Thus, for studing Wulff shapes having non-smooth surface energy functions, 
it is very significant 
to answer the  following two problems: 
\lq\lq (a) Generalize Cahn-Hoffman vector formula to the corresponding formula 
for any $\widetilde{\nu}: N\to S^n$\rq\rq\. 
and \lq\lq (b) Resolution of singularities of the boundary 
of a convex body having 
non-smooth boundary 
by a frontal $\widetilde{f}: S^n\to \mathbb{R}^{n+1}$\rq\rq.    
By 
{\color{black}Theorem \ref{theorem1} (2), } 
the problem (a) is completely solved.     
As for the problem (b), to the best of author's knowledge, 
only the boundary of a square 
has been realized as a frontal $\widetilde{f}: S^1\to \mathbb{R}^2$ so far   
(see \cite{janeczkonishimura}).    
{\color{black}Although there are apparently no published proofs at present, 
it is a comparatively straightforward generalization of this result to show that 
the boundary of a convex polygon is realized 
as a frontal $\widetilde{f}: S^1\to \mathbb{R}^2$.    
However, even in the plane case, the problem (b) for 
the boundary of a convex body in general seems to be wrapped in mystery.}
\par  
Moreover, 
{\color{black}Theorem \ref{theorem1} (2)} 
might be useful 
even for the study of force problems in higher dimensional vector spaces.         
In \cite{blaschke}, 
Petr Blaschke discovered that pedal coordinates are more suitable 
settings to study force problems in $\mathbb{R}^2$.      
Readers who want to confirm their usefulness are recommeded 
to refer to \cite{blaschke} 
(see also 7.24 (6) of \cite{brucegiblin} though this is not 
a force problem but a very suitable problem for understanding 
how useful pedal coordinates are).    
{\color{black}Theorem \ref{theorem1} (2)}  
may be regarded 
as a higher dimensional generalization of 
pedal coordinates. 
Hence, it is expected that 
{\color{black}Theorem \ref{theorem1} (2)} 
is a very suitable expression 
to study force problems etc. in all finite-dimensional vector spaces over 
$\mathbb{R}$.   
Example \ref{example0} {\color{black}(2)} might be regarded as examples in which 
higher dimensional version of pedal coordinates are effectively used.   
}}
\par 
\smallskip 
As an application of Theorem \ref{theorem1},  
a characterization for a hyperplane family 
to create a unique envelope is given as follows.   
\begin{theorem}\label{theorem2}
Let $\widetilde{\varphi}: N\to \mathbb{R}^{n+1}$, 
$\widetilde{\nu}: N\to S^n$ be mappings.   
Then, the hyperplane family 
$\mathcal{H}_{\left(\widetilde{\varphi}, \widetilde{\nu}\right)}$ 
creates a unique envelope 
if and only if it is creative  
and the set consisting of regular points of $\widetilde{\nu}$ is dense in $N$.      
\end{theorem}  
{\color{black}{Under the assumption that $\Omega$ 
in Remark 1.1 (2) is immersive and  
some conditions are satisfied, 
a unique existence result of envelopes 
for hyperplane families 
has been obtained in \cite{fisher}.     
Since their assumptions clearly imply that   
the creative condition defined 
in Definition \ref{definition2} is satisfied and the set consisting of regular points 
of $\widetilde{\nu}$ is dense, 
their result follows from Theorem \ref{theorem1} and Theorem \ref{theorem2}.  
\par   
Notice that non-unique existence cases, too, are intriguing cases since 
Theorem \ref{theorem1} 
may be effectively applied 
even in such cases 
(see Example \ref{example0} {\color{black}(1), (2)}).    
}} 
\par 
\bigskip 
This paper is organized as follows. 
Theorem \ref{theorem1} and Theorem \ref{theorem2}   
are proved in Section \ref{section3} and Section \ref{section4} respectively.    
In Section \ref{section2}, 
examples 
are given.   
{\color{black}Section \ref{section5} is an appendix where an alternative proof  
of Theorem \ref{theorem1} except for Theorem \ref{theorem1} (3) is  
given.    
The alternative proof 
is a proof by a gauge theoretic approach.   In order to avoid 
unnecessary complication, 
the alternative proof is given only in the case 
$n=1$.   The author has no idea on how to prove Theorem 
\ref{theorem1} (3) by using the alternative proof.}      
\section{Proof of Theorem \ref{theorem1}}\label{section3} 
\subsection{Proof of Theorem \ref{theorem1} (1)}\label{proof_1(1)}
\subsubsection{Proof of \lq\lq if\rq\rq\, part}\label{subsection3.1}
Let $x_0$ be an arbitrary point of $N$.   
Take one point $P$ of 
$\mathbb{R}^{n+1}-
H_{\left(\widetilde{\varphi}\left(x_0\right), \widetilde{\nu}\left(x_0\right)\right)}$ 
and fix it.  
It follows 
$\left(\widetilde{\varphi}\left(x_0\right)-P\right)\cdot 
\widetilde{\nu}\left(x_0\right)\ne 0$.    
Let $\widetilde{U}_{{}_P}$ be the set of points $x\in N$ satisfying 
\[
\left(\widetilde{\varphi}(x)-P\right)\cdot \widetilde{\nu}(x)\ne 0.   
\leqno{(2.1)}
\label{(3.1)}
\] 
Then, it is clear that $\widetilde{U}_{{}_P}$ is an open neighborhood of $x_0$ and 
the mirror image of the fixed point $P$ 
by the mirror $H_{\left(\widetilde{\varphi}(x), \widetilde{\nu}(x)\right)}$ is given by 
\[
2\left(\left(\widetilde{\varphi}(x)-P\right)
\cdot\widetilde{\nu}(x)\right)\widetilde{\nu}(x)+P
\]     
for any $x\in \widetilde{U}_{{}_P}$.     
\begin{figure}
\begin{center}
\includegraphics[width=15cm]
{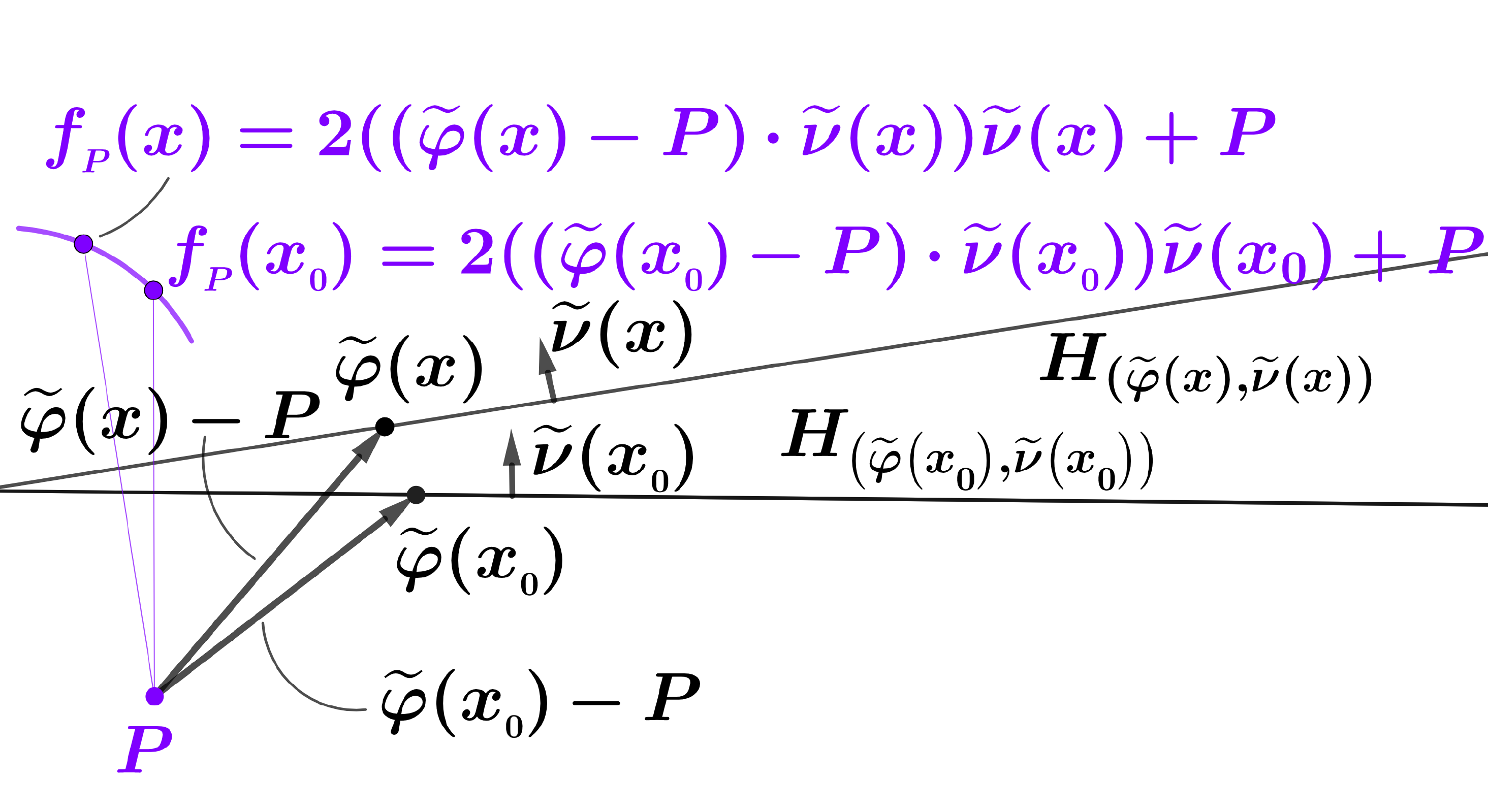}
\caption{
The mirror-image mapping $f_{{}_P}: U_P\to \mathbb{R}^{n+1}$.    
} 
\label{proof1}
\end{center}
\end{figure}
\par 
Since the hyperplane family 
$\mathcal{H}_{\left(\widetilde{\varphi}, \widetilde{\nu}\right)}$ is assumed to be 
creative, 
there exists a mapping 
$\widetilde{\Omega}: N\to T^*S^n$ with the form 
$\widetilde{\Omega}(x)=
\left(\widetilde{\nu}(x), \widetilde\omega(x)\right)$  such that 
for any $x\in N$ the following equality holds as $1$-form germs at $x$.     
\[
d\widetilde{\gamma}=\widetilde\omega.  
\]
Let $\left(V, \left(\Theta_1, \ldots, \Theta_n\right)\right)$ 
be a normal coordinate 
neighborhood of $S^n$ at $\widetilde{\nu}\left(x_0\right)$.    
Set $U_{{}_P}=\widetilde{U}_{{}_P}\cap \widetilde{\nu}^{-1}(V)$.     
Consider the mirror-image mapping 
$f_{{}_P}: U_{{}_P}\to \mathbb{R}^{n+1}$ defined by 
\[
f_{{}_P}(x)=
2\left(\left(\widetilde{\varphi}(x)-P\right)
\cdot\widetilde{\nu}(x)\right)\widetilde{\nu}(x)+P 
\]
for any $x\in U_{{}_P}$.   
In order to show that $f_{{}_P}$ is a frontal, it is sufficient 
to construct a Gauss mapping with respect to $f_{{}_P}$.    
By using the mapping $\widetilde{\Omega}|_{U_{{}_P}}$, 
a Gauss mapping for $f_{{}_P}$ is constructed as follows.    
For any $x\in U_{{}_P}$ set $X=\widetilde{\nu}(x)$.     
Let ${\color{black}\Pi}_{\left(X, X_0\right)}: T_{X_0}S^n\to T_X S^n$ 
be the Levi-Civita translation.  
For any $i$ $(1\le i\le n)$, 
set $\frac{\partial}{\partial \Theta_{\left(i, X\right)}}=
{\color{black}\Pi}_{\left(X, X_0\right)}\left(\frac{\partial}{\partial \Theta_i}\right)$.     
Then notice that for any $x\in U_{{}_P}$, 
under the identification of $\mathbb{R}^{n+1}$ and 
$T_{f_{{}_P}(x)}\mathbb{R}^{n+1}$, 
\[
\left\langle \frac{\partial}{\partial \Theta_{\left(1, X\right)}}, \ldots, 
\frac{\partial}{\partial \Theta_{\left(n, X\right)}}, 
\widetilde{\nu}(x) \right\rangle
\] 
is an orthonormal basis 
of the tangent vector space $T_{f_{{}_P}(x)}\mathbb{R}^{n+1}$.  
\begin{figure}
\begin{center}
\includegraphics[width=15cm]
{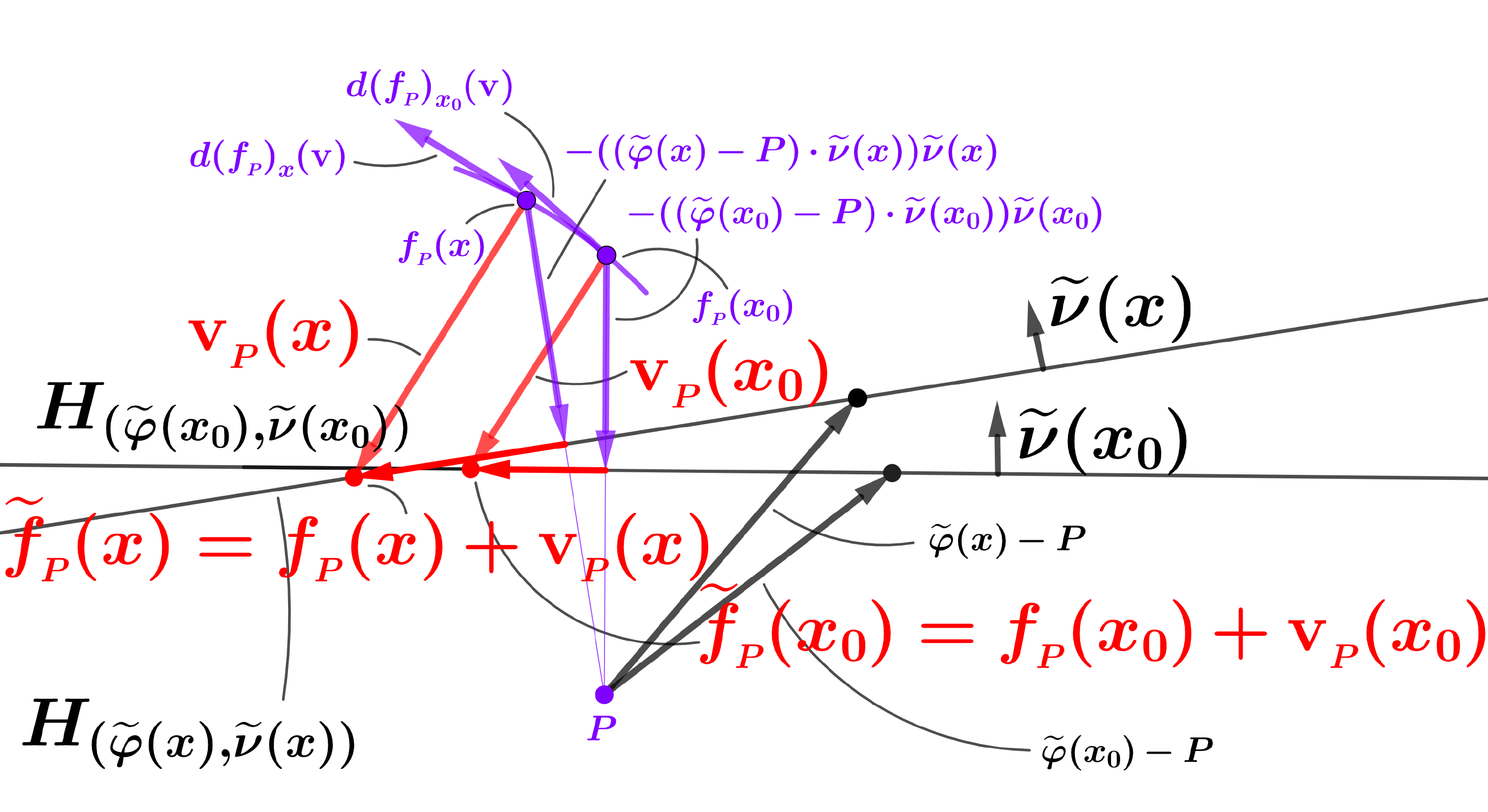}
\caption{
Figure for Proof of \lq\lq if\rq\rq\, part. 
} 
\label{proof2}
\end{center}
\end{figure}
\begin{lemma}\label{lemma1} 
For any $x\in U_{{}_P}$, the following equality holds.   
\[
d\left(P\cdot \widetilde{\nu}\right) =  
\sum_{i=1}^n\left(P\cdot \frac{\partial}{\partial \Theta_{\left(i, X\right)}}\right)
d\left(\Theta_i\circ \widetilde{\nu}\right).   
\] 
\end{lemma}
\noindent 
\underline{\it Proof of Lemma \ref{lemma1}}.\quad 
\begin{eqnarray*}
d\left(P\cdot \widetilde{\nu}\right) & = & 
\sum_{j=1}^n 
\frac{\partial \left(P\cdot \widetilde{\nu}\right)}{\partial x_j}(x) d x_j \\ 
{ } & = & 
\sum_{j=1}^n \left(P\cdot \left(\sum_{i=1}^n
\frac{\partial \left(\Theta_i\circ \widetilde{\nu}\right)}{\partial x_j}(x)
\frac{\partial}{\partial \Theta_{(i,X)}}\right)\right)d x_j \\ 
{ } & = & 
\sum_{i=1}^n\left(P\cdot \frac{\partial}{\partial \Theta_{(i,X)}}\right)
\left(\sum_{j=1}^n
\frac{\partial \left(\Theta_i\circ \widetilde{\nu}\right)}
{\partial x_j}(x)d x_j\right) \\ 
{ } & = & 
\sum_{i=1}^n\left(P\cdot \frac{\partial}{\partial \Theta_{\left(i, X\right)}}\right)
d\left(\Theta_i\circ \widetilde{\nu}\right).   
\end{eqnarray*}
\hfill $\Box$
\par 
\noindent 
\smallskip 
By Lemma \ref{lemma1}, under the identification of 
$T_{\widetilde{\nu}(x)}S^n$ and $T^*_{\widetilde{\nu}(x)}S^n$, 
it follows 
\begin{eqnarray*}
d\left(\left(\widetilde{\varphi}-P\right)\cdot \widetilde{\nu}\right) 
& = & 
d\left(\widetilde{\varphi}\cdot \widetilde{\nu}\right)  -  
d\left(P\cdot \widetilde{\nu}\right) \\  
{ } & = & d\widetilde{\gamma}  -  
d\left(P\cdot \widetilde{\nu}\right) \\  
{ } & = & \widetilde\omega - 
d\left(P\cdot \widetilde{\nu}\right) \\  
{ } & = & 
\sum_{i=1}^n\left(\widetilde\omega(x)\cdot 
\frac{\partial}{\partial \Theta_{\left(i,X\right)}}\right)
 d \left(\Theta_i\circ \widetilde{\nu}\right) 
- \sum_{i=1}^n\left(P\cdot \frac{\partial}{\partial \Theta_{\left(i,X\right)}}\right)
d \left(\Theta_i\circ \widetilde{\nu}\right) \\
{ } & = & 
\sum_{i=1}^n\left(
\left(\widetilde\omega(x) - P\right)\cdot 
\frac{\partial}{\partial \Theta_{\left(i,X\right)}}
\right)
d \left(\Theta_i\circ \widetilde{\nu}\right)  
\end{eqnarray*} 
for any $x\in U_{{}_P}$.   
Set 
\[
{\bf v}_{{}_P}(x)=\sum_{i=1}^n\left(
\left(\widetilde\omega(x) - P\right)\cdot 
\frac{\partial}{\partial \Theta_{\left(i,X\right)}}
\right)
\frac{\partial}{\partial\Theta_{(i,X)}} 
- 
\left(\left(\widetilde{\varphi}(x)-P\right)\cdot 
\widetilde{\nu}(x)\right)\widetilde{\nu}(x)
\]   
for any $x\in U_{{}_P}$ 
where $\mathbb{R}^{n+1}$ and $T_{f_{{}_P}\left(x\right)}\mathbb{R}^{n+1}$ 
are identified and $T_{f_{{}_P}\left(x\right)}S^n$ 
and $T_{f_{{}_P}\left(x\right)}^*S^n$ are identified.    
By (\ref{(3.1)}), ${\bf v}_{{}_P}(x)$ is not the zero vector.           
Moreover, the following holds.  
\begin{lemma}\label{perpendicular}
For any ${\bf v}\in T_{x_0} N$, ${\bf v}_{{}_P}(x_0)$ 
is perpendicular to $d \left(f_{{}_P}\right)_{{}x_0}({\bf v})$.   
\end{lemma}
\noindent 
\underline{\it Proof of Lemma \ref{perpendicular}}.\quad 
Calculation of the product of the vector ${\bf v}_{{}_P}\left(x_0\right)$ 
and the Jacobian matrix of 
$f_{{}_P}$ at $x_0$  (denoted by $J\left(f_{{}_P}\right)_{x_0}$) 
is carried out as follows, 
where $\mathbb{R}^{n+1}$ and $T_{f_{{}_P}\left(x_0\right)}\mathbb{R}^{n+1}$ 
are identified and $T_{f_{{}_P}\left(x_0\right)}S^n$ 
and $T_{f_{{}_P}\left(x_0\right)}^*S^n$ are identified.        
\begin{eqnarray*}
{ } & { } & 
{\bf v}_{{}_P}\left(x_0\right)J\left(f_{{}_P}\right)_{x_0} \\ 
{ } & = & 
2\sum_{i=1}^n\left(\left(\widetilde\omega\left(x_0\right)-P\right)\cdot 
\frac{\partial}{\partial \Theta_i}\right)
\left(\left(\widetilde{\varphi}\left(x_0\right)-P\right)
\cdot \widetilde{\nu}\left(x_0\right)\right)
d\left(\Theta_i\circ \widetilde{\nu}\right)  \\ 
{ } & { } & \qquad 
- 
2\left(\left(\widetilde{\varphi}\left(x_0\right)-P\right)
\cdot \widetilde{\nu}\left(x_0\right)\right)
d\left(\left(\widetilde{\varphi}-P\right)\cdot 
\widetilde{\nu}\right)_{\mbox{at }x_0} \\ 
{ } & = & 
2\left(\left(\widetilde{\varphi}\left(x_0\right)-P\right)
\cdot \widetilde{\nu}\left(x_0\right)\right)
\sum_{i=1}^n\left(\left(\widetilde\omega\left(x_0\right)-P\right)\cdot 
\frac{\partial}{\partial \Theta_i}\right)
d\left(\Theta_i\circ \widetilde{\nu}\right)  \\ 
{ } & { } & \qquad 
- 
2\left(\left(\widetilde{\varphi}\left(x_0\right)-P\right)
\cdot \widetilde{\nu}\left(x_0\right)\right)
\sum_{i=1}^n
\left(\left(\widetilde\omega\left(x_0\right)-P\right)\cdot 
\frac{\partial}{\partial \Theta_i}\right)
d\left(\Theta_i\circ \widetilde{\nu}\right)  \\ 
{ } & = & 0.
\end{eqnarray*}
\hfill $\Box$
\par 
\noindent 
\smallskip 
We may consider that the point $x_0$ is an arbitrary point of $U_{{}_P}$.   
Thus we have the following.   
\begin{lemma}\label{lemma2}
The mapping $f_{{}_P}: U_{{}_P}\to \mathbb{R}^{n+1}$ 
is a frontal with Gauss mapping 
$\nu_{{}_P}: U_{{}_P}\to S^n$ such that 
$\nu_{{}_P}(x)\cdot \widetilde{\nu}(x)\ne 0$, 
where $\nu_{{}_P}(x)=\frac{{\bf v}_{{}_P}(x)}{\Vert {\bf v}_{{}_P}(x)\Vert}$.     
\end{lemma}
\noindent 
By Lemma \ref{lemma2}, 
the hyperplane $H_{\left(\widetilde{\varphi}(x), \widetilde{\nu}(x)\right)}$ and 
the line $\ell_x=\left\{f_{{}_P}(x)+t\nu_{{}_P}(x) \left|\right. t\in \mathbb{R}\right\}$ 
must intersect only at one point for any $x\in U_{{}_P}$.          
Define the mapping 
$\widetilde{f}_{{}_P}: U_{{}_P}\to \mathbb{R}^{n+1}$ 
by 
\[
\left\{\widetilde{f}_{{}_P}(x)\right\}= 
H_{\left(\widetilde{\varphi}(x), \widetilde{\nu}(x)\right)}\cap \ell_x.
\]
Then, from the construction, 
$\widetilde{f}_{{}_P}$ must have the following form 
(see p.7 of \cite{janeczkonishimura}).    
\[
\widetilde{f}_{{}_P}(x)=f_{{}_P}(x)-\frac{||f_{{}_P}(x)-P||^2}{2\left(f_{{}_P}(x)-P\right)
\cdot \nu_{{}_P}(x)}\nu_{{}_P}(x).   
\]
By Theorem 1 of \cite{janeczkonishimura} 
{\color{black}{
(more precisely, by 3.1 in p.9 of \cite{janeczkonishimura})}} 
and Lemma \ref{lemma2}, 
we have the following.   
\begin{lemma}\label{lemma3}
The mapping $\widetilde{f}_{{}_P}$ is a frontal with Gauss mapping 
$\widetilde{\nu}|_{U_{{}_P}}: U_{{}_P}\to S^n$.    
In other words, $\widetilde{f}_{{}_P}: U_{{}_P}\to \mathbb{R}^{n+1}$ is an envelope 
created by the hyperplane family 
$\mathcal{H}_{\left(\widetilde{\varphi}|_{U_{{}_P}}, \widetilde{\nu}|_{U_{{}_P}}\right)}$.     
\end{lemma}
\par 
\smallskip 
On the other hand, it is easily seen that  
$ 
\left(
f_{{}_P}\left(x_0\right)+{\bf v}_{{}_P}\left(x_0\right)
-\widetilde{\varphi}\left(x_0\right)\right)\cdot 
\widetilde{\nu}\left(x_0\right) =0$ (see Figure \ref{proof2}).   
Thus, the vector $f_{{}_P}\left(x_0\right)+{\bf v}_{{}_P}\left(x_0\right)$ 
must belong to 
$H_{\left(\widetilde{\varphi}\left(x_0\right), 
\widetilde{\nu}\left(x_0\right)\right)}$. 
From the construction and by using the equality 
$P=\sum_{i=1}^n 
\left(P\cdot \frac{\partial}{\partial \Theta_i}\right)
\frac{\partial}{\partial \Theta_i}
+ 
\left(P\cdot \widetilde{\nu}\left(x_0\right)\right)\widetilde{\nu}\left(x_0\right),
$ 
we have the following. 
\begin{eqnarray*}
 \widetilde{f}_{{}_P}\left(x_0\right) 
 & = & f_{{}_P}(x)+{\bf v}_{{}_P}\left(x_0\right) \\ 
{ } & = & 
2\left(\left(\widetilde{\varphi}\left(x_0\right)-P\right)
\cdot\widetilde{\nu}\left(x_0\right)\right)\widetilde{\nu}\left(x_0\right)+P \\ 
{ } & { } & \quad\quad 
+ 
\sum_{i=1}^n\left(
\left(\widetilde\omega\left(x_0\right) - P\right)\cdot 
\frac{\partial}{\partial \Theta_i}
\right)
\frac{\partial}{\partial\Theta_i} 
- 
\left(\left(\widetilde{\varphi}\left(x_0\right)-P\right)\cdot 
\widetilde{\nu}\left(x_0\right)\right)\widetilde{\nu}\left(x_0\right) \\ 
{ } & = &  \left(\left(\widetilde{\varphi}\left(x_0\right)-P\right)
\cdot\widetilde{\nu}\left(x_0\right)\right)\widetilde{\nu}\left(x_0\right)+P+
\sum_{i=1}^n\left(
\left(\widetilde\omega\left(x_0\right) - P\right)\cdot 
\frac{\partial}{\partial \Theta_i}
\right)
\frac{\partial}{\partial\Theta_i}  \\ 
{ } & = & 
\left(\widetilde{\varphi}\left(x_0\right)
\cdot\widetilde{\nu}\left(x_0\right)\right)\widetilde{\nu}\left(x_0\right)+
\sum_{i=1}^n\left(\widetilde\omega\left(x_0\right) \cdot \frac{\partial}
{\partial \Theta_i}
\right)
\frac{\partial}{\partial\Theta_i} \\ 
{ } & = & 
\widetilde{\gamma}\left(x_0\right)
\widetilde{\nu}\left(x_0\right)+
\widetilde\omega\left(x_0\right).   
\end{eqnarray*}
This proves the following lemma.    
\begin{lemma}\label{lemma4}
The following equality holds.   
\[
\widetilde{f}_{{}_P}\left(x_0\right) = 
\widetilde{\gamma}\left(x_0\right)\widetilde{\nu}\left(x_0\right)
+
\widetilde\omega\left(x_0\right). 
\]
\end{lemma}
\noindent 
Lemma \ref{lemma4} shows that 
$\widetilde{f}_{{}_P}\left(x_0\right)$ does not depend on the particular choice of 
$P\in 
\mathbb{R}^{n+1}-
H_{\left(\widetilde{\varphi}\left(x_0\right), \widetilde{\nu}\left(x_0\right)\right)}$.   
Define the mapping  $\widetilde{f}: N\to \mathbb{R}^{n+1}$ by 
$\widetilde{f}(x) = 
\widetilde{\gamma}(x)\widetilde{\nu}(x)+
\widetilde\omega(x)$.    
Since $x_0$ is an arbitrary point of $N$, 
by Lemma \ref{lemma3} and Lemma \ref{lemma4}, 
it follows that 
the mapping $\widetilde{f}: N\to \mathbb{R}^{n+1}$ is an envelope created by 
$\mathcal{H}_{\left(\widetilde{\varphi}, \widetilde{\nu}\right)}$.  
This completes the proof of \lq\lq if\rq\rq\, part.   
\hfill $\Box$
\par 
\medskip 
\noindent 
\subsubsection{Proof of \lq\lq only if\rq\rq\, part}
\begin{figure}
\begin{center}
\includegraphics[width=15cm]
{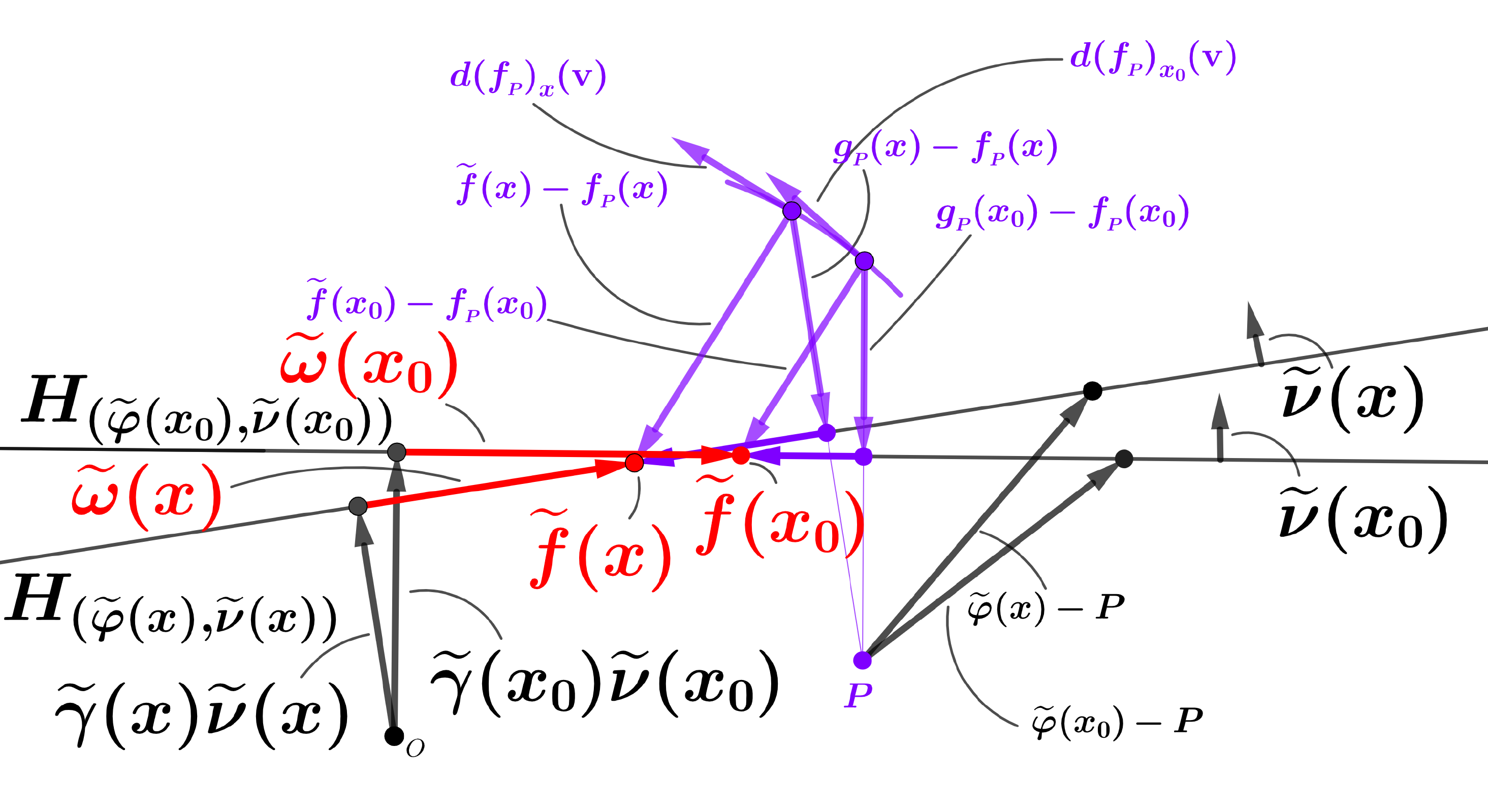}
\caption{
Figure for Proof of \lq\lq only if\rq\rq\, part.       
} 
\label{proof3}
\end{center}
\end{figure}

Suppose that the hyperplane family 
$\mathcal{H}_{\left(\widetilde{\varphi}, \widetilde{\nu}\right)}$ creates an envelope 
$\widetilde{f}: N\to \mathbb{R}^{n+1}$.    Then, by definition, 
$\widetilde{f}$ is a frontal such that the inclusion 
$\widetilde{f}(x)+d\widetilde{f}_x(T_x N)
\subset H_{\left(\widetilde{\varphi}(x), \widetilde{\nu}(x)\right)}$ 
holds for any $x\in N$.      
Let $\widetilde\omega: N\to \mathbb{R}^{n+1}$ be the mapping defined by 
$\widetilde\omega(x)=\widetilde{f}(x)
-\widetilde{\gamma}(x)\widetilde{\nu}(x)$ (see Figure \ref{proof3}).    
It is sufficient to show that under some identifications,  
$\widetilde\omega$ is actually a creator for the envelope $\widetilde{f}$.  
\par 
It is easily seen that $\widetilde\omega(x)\cdot \widetilde{\nu}(x)=0$ 
for any $x\in N$. 
Thus, under the identification of $\mathbb{R}^{n+1}$ and 
$T^*_{\widetilde{\nu}(x)}\mathbb{R}^{n+1}$, we have 
\begin{lemma}\label{lemma4.5}
For any $x\in N$, $\widetilde\omega(x)\in T^*_{\widetilde{\nu}(x)}S^n$ holds.   
\end{lemma} 
Let $\widetilde{\Omega} : N\to T^*S^n$ be the mapping defined by 
$\widetilde{\Omega}(x)=\left(\widetilde{\nu}(x), \widetilde\omega(x)\right)$.   
Let $x_0$ be an arbitrary point of $N$ and let $P$ be a point of 
$\mathbb{R}^{n+1}-
H_{\left(\widetilde{\varphi}(x_0), \widetilde{\nu}(x_0)\right)}$.  
Again, we consider the mirror-image mapping 
$f_{{}_P}: \widetilde{U}_{{}_P}\to \mathbb{R}^{n+1}$ defined by 
\[
f_{{}_P}(x)=2\left(\left(\widetilde{\varphi}(x)-P\right)
\cdot\widetilde{\nu}(x)\right)\widetilde{\nu}(x)+P, 
\]    
where $\widetilde{U}_{{}_P}=\left\{x\in N\,  \left|\, 
\left(\widetilde{\varphi}(x)-P\right)\cdot \widetilde{\nu}(x)\ne 0 \right\}\right.$.   
The mapping $f_{{}_P}$ is exactly 
the orthotomic of $\widetilde{f}|_{\widetilde{U}_{{}_P}}$ 
relative to the point $P$.    
Thus, by Proposition 1 of \cite{janeczkonishimura} {\color{black}{
(more precisely, by 2.1 in pp. 7--8 of \cite{janeczkonishimura}) 
}}, 
$f_{{}_P}$ is a frontal and the mapping  
$\nu_{{}_P}: \widetilde{U}_{{}_P}\to S^n$ define by 
\[
\nu_{{}_P}(x)=\frac{\widetilde{f}(x)-f_{{}_P}(x)}{\Vert \widetilde{f}(x)-f_{{}_P}(x)\Vert}
\]
is its Gauss mapping.    
In particular, we have the following.   
\begin{lemma}\label{lemma5}
For any $x\in \widetilde{U}_{{}_P}$ and any 
${\bf v}\in T_x N$, the following holds.   
\[
\left(\widetilde{f}(x)-f_{{}_P}(x)\right)\cdot 
d \left(f_{{}_P}\right)_x({\bf v})=0.   
\]
\end{lemma}
For any $x\in \widetilde{U}_{{}_P}$, 
set 
\[
g_{{}_P}(x)=\frac{1}{2}\left(f_{{}_P}(x)-P\right)+P=
\left(\left(\widetilde{\varphi}(x)-P\right)
\cdot\widetilde{\nu}(x)\right)\widetilde{\nu}(x)+P.   
\] 
Then, since $f_{{}_P}(x)$ is the mirror-image of $P$ with respect to the mirror 
$H_{\left(\widetilde{\varphi}(x), \widetilde{\nu}(x)\right)}$, 
the following clearly holds.  
\begin{lemma}\label{lemma6}
The vector $\widetilde{f}(x)-g_{{}_P}(x)$ is perpendicular to the vector 
$g_{{}_P}(x)-f_{{}_P}(x)=-\left(\left(\widetilde{\varphi}(x)-P\right)
\cdot \widetilde{\nu}(x)\right)\widetilde{\nu}(x)$ 
for any $x\in \widetilde{U}_{{}_P}$.    
\end{lemma}
\noindent 
Thus, 
\[
\widetilde{f}(x)-f_{{}_P}(x)=
\left(\widetilde{f}(x)-g_{{}_P}(x)\right)+\left(g_{{}_P}(x)-f_{{}_P}(x)\right)
\]
is an orthogonal decomposition of $\widetilde{f}(x)-f_{{}_P}(x)$ 
for any $x\in \widetilde{U}_{{}_P}$  (see Figure \ref{proof3}).    
\par 
In order to decompose the vector $\widetilde{f}(x)-g_{{}_P}(x)$ 
reasonably, the open neighborhood  
$\widetilde{U}_{{}_P}$ of $x_0$ is reduced as follows.   
Let $\left(V, \left(\Theta_1, \ldots, 
\Theta_n\right)\right)$ be a normal coordinate neighborhood of $S^n$ at 
$\widetilde{\nu}\left(x_0\right)$.    
Set again 
$U_{{}_P}=\widetilde{U}_{{}_P}\cap \widetilde{\nu}^{-1}(V)$.    
Notice that 
$ 
\left\langle 
d \Theta_1, \ldots, 
d \Theta_n
\right\rangle
$ 
is an orthonormal basis of the cotangent space 
$T^*_{\widetilde{\nu}\left(x_0\right)}S^n$.   
\begin{lemma}\label{lemma7}
The equality 
\[
\widetilde{f}\left(x_0\right)-g_{{}_P}\left(x_0\right) = 
\widetilde\omega\left(x_0\right) -
\sum_{i=1}^n 
\left(P\cdot 
\frac{\partial}{\partial \Theta_i}\right)
\frac{\partial}{\partial\Theta_i} 
\]
holds where three vector spaces $\mathbb{R}^{n+1}$, 
$T_{\widetilde{\nu}\left(x_0\right)}\mathbb{R}^{n+1}$ and 
$T^*_{\widetilde{\nu}\left(x_0\right)}\mathbb{R}^{n+1}$ 
are identified. 
\end{lemma}
\noindent 
\underline{Proof of Lemma \ref{lemma7}}.\quad     
\begin{eqnarray*}
\widetilde{f}\left(x_0\right)-g_{{}_P}\left(x_0\right) & = & 
\widetilde{f}\left(x_0\right) - 
\left(\left(\left(\widetilde{\varphi}\left(x_0\right)-P\right)
\cdot \widetilde{\nu}\left(x_0\right)\right)\widetilde{\nu}\left(x_0\right)
+P\right) \\ 
{ } & = & 
\left(\widetilde{f}\left(x_0\right) 
-
\left(\widetilde{\varphi}\left(x_0\right)
\cdot \widetilde{\nu}\left(x_0\right)\right)\widetilde{\nu}\left(x_0\right) \right)
+
\left(\left(P\cdot \widetilde{\nu}\left(x_0\right)\right)
\widetilde{\nu}\left(x_0\right)
-P\right)
\\ 
{ } & = & 
\left(\widetilde{f}\left(x_0\right) 
- 
\widetilde{\gamma}\left(x_0\right)
\widetilde{\nu}\left(x_0\right) \right)
+
\left(\left(P\cdot \widetilde{\nu}\left(x_0\right)\right)
\widetilde{\nu}\left(x_0\right)
-P\right)
\\ 
{ } & = & 
\widetilde\omega\left(x_0\right) -
\sum_{i=1}^n 
\left(P\cdot 
\frac{\partial}{\partial \Theta_i}\right)
\frac{\partial}{\partial\Theta_i}.    
\end{eqnarray*}
\hfill $\Box$
\par 
\noindent 
By Lemma \ref{lemma7}, the following holds.  
\begin{eqnarray*} 
\widetilde{f}\left(x_0\right)-f_{{}_P}\left(x_0\right) & = & 
\left(\widetilde{f}\left(x_0\right)-g_{{}_P}\left(x_0\right)\right)
+\left(g_{{}_P}\left(x_0\right)-f_{{}_P}\left(x_0\right)\right)\\ 
{ } & = & 
\widetilde\omega\left(x_0\right) -
\sum_{i=1}^n 
\left(P\cdot 
\frac{\partial}{\partial \Theta_i}\right)
\frac{\partial}{\partial\Theta_i}
- 
\left(\left(\widetilde{\varphi}\left(x_0\right)-P\right)\cdot 
\widetilde{\nu}\left(x_0\right)\right)\widetilde{\nu}\left(x_0\right).
\end{eqnarray*}
Hence, by Lemma \ref{lemma1} and Lemma \ref{lemma5}, the germ of $1$-form 
 $d\widetilde{\gamma}$ at $x_0$ is calculated 
as follows, where $X=\widetilde{\nu}(x)$, $
\frac{\partial}{\partial \Theta_{\left(i, X\right)}}=
P_{\left(X, X_0\right)}\left(\frac{\partial}{\partial \Theta_i}\right).    
$ and $P_{\left(X, X_0\right)}: T_{X_0}S^n\to T_X S^n$ 
is the Levi-Civita translation.  
\begin{eqnarray*}
d\widetilde{\gamma} 
& = & 
d\widetilde{\gamma}
- 
d\left(P\cdot \widetilde{\nu}\right) 
+ 
d\left(P\cdot \widetilde{\nu}\right) \\ 
{ } & = & 
d\left(\left(\widetilde{\varphi}-P\right)\cdot \widetilde{\nu}\right) 
+ 
d\left(P\cdot \widetilde{\nu}\right) \\ 
{ } & = & 
\sum_{i=1}^n 
\left(\left(\widetilde{\omega}-P\right)\cdot 
\frac{\partial}{\partial \Theta_{\left(i, X\right)}}\right)
d\left(\Theta_i\circ\widetilde{\nu}\right)
+ 
\sum_{i=1}^n 
\left(P\cdot 
\frac{\partial}{\partial \Theta_{\left(i, X\right)}}\right)
d\left(\Theta_i\circ\widetilde{\nu}\right)
\\ 
{ } & = & 
\left(
\widetilde\omega 
-
\sum_{i=1}^n 
\left(P\cdot 
\frac{\partial}{\partial \Theta_{\left(i, X\right)}}\right)
d\left(\Theta_i\circ\widetilde{\nu}\right)
\right)
+ 
\sum_{i=1}^n 
\left(P\cdot 
\frac{\partial}{\partial \Theta_{\left(i, X\right)}}\right)
d\left(\Theta_i\circ\widetilde{\nu}\right)
\\ 
{ } & = & 
\widetilde\omega.    
\end{eqnarray*}  
This calculation proves the following lemma.   
\begin{lemma}\label{lemma8}
The equality 
\[
d\widetilde{\gamma} = 
\widetilde\omega 
\] 
holds as germs of $1$-form at $x_0$. 
\end{lemma} 
Since $x_0$ is an arbitrary point of $N$, 
by Lemma \ref{lemma8}, $\widetilde\omega$ 
is actually the creator 
for the given envelope $\widetilde{f}: N\to \mathbb{R}^{n+1}$.      
This completes the proof of \lq\lq only if\rq\rq\, part.   
\hfill $\Box$
{\color{black}\subsection{Proof of Theorem \ref{theorem1} (2)}\label{proof_1(2)}
Theorem \ref{theorem1} (2) is a direct by-product of the proof of 
Theorem \ref{theorem1} (1).   
\qed 
\subsection{Proof of Theorem \ref{theorem1} (3)}\label{proof_1(3)}
Recall that the line family 
$\mathcal{H}_{\left(\widetilde{\varphi}, \widetilde{\nu}\right)}$ is said to 
{\it create 
an $E_1$ envelope} (denoted by $(E_1)$ in this subsection) 
if for any fixed $t_0\in N$ and any $t\in N$ near $t_0$ 
the limit 
$
\lim_{t\to t_0}
H_{\left(\widetilde{\varphi}(t), \widetilde{\nu}(t)\right)}\cap 
H_{\left(\widetilde{\varphi}(t_0), \widetilde{\nu}(t_0)\right)}
$  
exists.    
On the other hand, the line family 
$\mathcal{H}_{\left(\widetilde{\varphi}, \widetilde{\nu}\right)}$ is said to 
{\it create 
an $E_2$ envelope}  (denoted by $(E_2)$ in this subsection) 
if it creates an envelope in the sense of Definition 
\ref{definition1}.    
\par 
\noindent 
\medskip 
\underline{$(E_1)\Rightarrow (E_2)$}\quad 
Let $t_0$ be a point of $N$ and let $t_i\in N$ $(i=1, 2, \ldots)$ 
be a sequence conversing to $t_0$.    
Since $(E_1)$ is assumed, 
we can assume that a point $X_{t_i}$ can be taken from the intersection 
$
H_{\left(\widetilde{\varphi}(t), \widetilde{\nu}(t)\right)}\cap 
H_{\left(\widetilde{\varphi}(t_0), \widetilde{\nu}(t_0)\right)}
$
such that $\lim_{t_i\to t_0}X_{t_i}$ exists.  Denote the limit by $X_{t_0}$.      
Then, we have the following.   
\begin{eqnarray*}
\left(X_{t_i}-\widetilde{\varphi}({t_i})\right)\cdot \widetilde{\nu}(t_i) & = & 0, \\
\left(X_{t_i}-\widetilde{\varphi}(t_0)\right)\cdot \widetilde{\nu}(t_0) & = & 0. 
\end{eqnarray*}
This implies 
\[
X_{t_i}\cdot \left(\widetilde{\nu}(t_i)-\widetilde{\nu}(t_0)\right) = 
\widetilde{\gamma}(t_i)-\widetilde{\gamma}(t_0).     
\]
Thus we have 
\[
X_{t_0}\cdot \frac{\partial \widetilde{\nu}}{\partial t}(t_0)=
\frac{\partial \widetilde{\gamma}}{\partial t}(t_0).   
\]   
This implies that there exists a real number $\alpha(t_0)$ such that 
the following identity holds where  
$d \left(\Theta\circ \widetilde{\nu}\right)$ and 
$d\widetilde{\gamma}$ stand for the 
$1$-dimensional cotangent vectors in  
$T_{t_0}^*N$, namely the following identity is nothing but 
the identity of two real numbers.  
\[
\alpha(t_0)d \left(\Theta\circ \widetilde{\nu}\right)=d\widetilde{\gamma}.   
\]
It is not difficult to see that the function $\alpha: N\to \mathbb{R}$ is of 
class $C^\infty$.    
This means that the line family 
 $\mathcal{H}_{\left(\widetilde{\varphi}, \widetilde{\nu}\right)}$ 
is creative.   
Therefore, by Theorem \ref{theorem1} (1), the line family creates an 
$E_2$ envelope.   
\qed 
\par 
\noindent 
\medskip 
\underline{$(E_2)\Rightarrow (E_1)$}\quad 
For the proof of this implication, 
it is used the notions and notations 
introduced in the proof of Theorem \ref{theorem1} (1).    
The assumption $(E_2)$ implies that $\widetilde{\gamma}$ is totally 
differentiable with respect to $\widetilde{\nu}$.     
Take an arbitrary point $t_0\in N$ and fixed it.   
Since  $\widetilde{\gamma}$ is totally 
differentiable with respect to $\widetilde{\nu}$ at $t_0$, 
for any $t$ near $t_0$ if the length of the vector 
$\overrightarrow{f_P(t_0)f_P(t)}$ is positive, 
then the horizontal vector of $\overrightarrow{f_P(t_0)f_P(t)}$ 
must be non-zero, where $P$ is a point taken outside 
the line $H_{\left(\widetilde{\varphi}(t_0), \widetilde{\nu}(t_0)\right)}$ and 
$f_P$ is a mirror-image mapping introduced in the proof of Theorem 
\ref{theorem1} (1).   
Denote the intersection of the perpendicular bisector of 
$\overrightarrow{f_P(t_0)f_P(t)}$ and the line 
$H_{\left(\widetilde{\varphi}(t_0), \widetilde{\nu}(t_0)\right)}$ 
by $J_t$.    
Then, from the construction, it follows that the triangre 
$\triangle J_tf_P(t_0)f_P(t)$ is an isosceles triangle with 
legs $J_tf_P(t_0)$ and $J_tf_P(t)$.   
This implies the following (see Figure \ref{figure14}).   
\[
J_t\in 
H_{\left(\widetilde{\varphi}(t), \widetilde{\nu}(t)\right)}
\cap 
H_{\left(\widetilde{\varphi}(t_0), \widetilde{\nu}(t_0)\right)}.
\] 
\begin{figure}[htbp]
    \centering
    \includegraphics[width=15cm]
{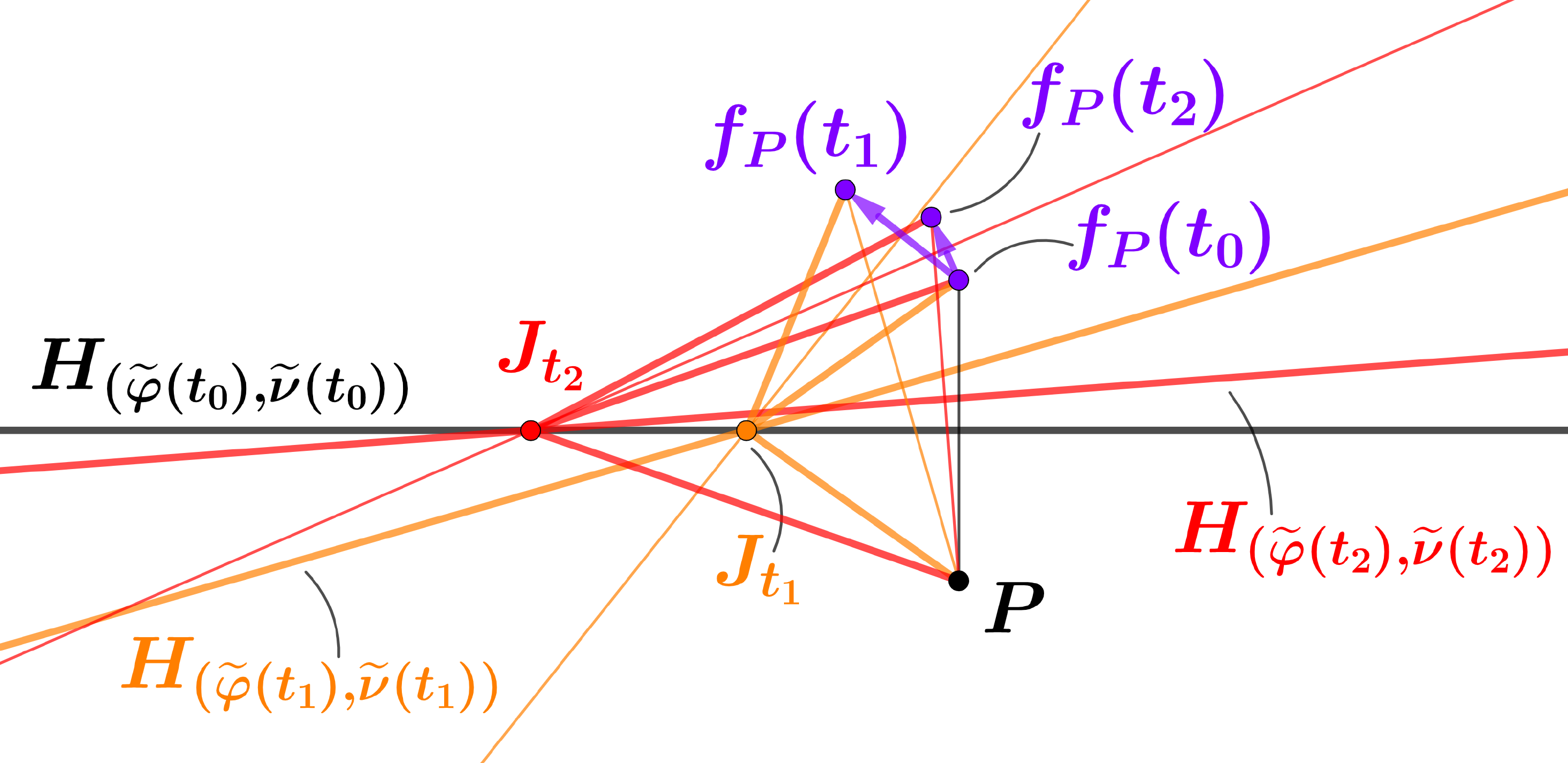}
    \caption{
Figure for $(E_2)\Rightarrow  (E_1)$.}
\label{figure14}
\end{figure}
Notice that $\lim_{t\to t_0}||J_tf_P(t_0)||$ is positive.   
Thus, we have 
\[
\lim_{t\to t_0}\angle J_tf_P(t_0)f_P(t)
= 
\lim_{t\to t_0}\angle J_tf_P(t)f_P(t_0)
=\frac{\pi}{2}.    
\]  
By Proposition 1 of \cite{janeczkonishimura} asserting that $f_P$ is a frontal 
with its Gauss mapping 
$\frac{f_P(t_0)-\widetilde{f}_P(t_0)}{||f_P(t_0)-\widetilde{f}_P(t_0)||}$, 
it follows 
\[
\lim_{t\to t_0}J_t=\widetilde{f}_P(t_0), 
\] 
where $\widetilde{f}_P$ is the anti-orthotomic of $f_P$ relative to the point $P$ 
introduced in the proof of Theorem \ref{theorem1} (1).   
Since $t_0$ is an arbitrary point of $N$, 
the given $E_2$ envelope must be an $E_1$ envelope 
by Theorem \ref{theorem1} (1). 
\qed  
} 
\section{Proof of Theorem \ref{theorem2}}\label{section4}
\noindent 
\underline{\it Proof of \lq\lq if\rq\rq\, part}. \quad     
Since the hyperplane 
$\mathcal{H}_{\left(\widetilde{\varphi}, \widetilde{\nu}\right)}$ is creative, 
by Theorem \ref{theorem1}, it creates an envelope.    
Let $\widetilde{f}_1, \widetilde{f}_2: N\to \mathbb{R}^{n+1}$ be 
envelopes created by 
$\mathcal{H}_{\left(\widetilde{\varphi}, \widetilde{\nu}\right)}$. 
\par  
Let $x_0\in N$ be a regular point of $\widetilde{\nu}$.      
Then, there exists an open coordinate neighborhood 
$\left(U, \left(x_1, \ldots, x_n\right)\right)$ such that $x_0\in U$ and 
$\widetilde{\nu}|_U: U\to \widetilde{\nu}(U)$ 
is a diffeomorphism. 
Then, the germ of $1$-form 
$d\left(\widetilde{\varphi}\cdot\widetilde{\nu}\right)$ 
at $x_0\in U$ is  
\begin{eqnarray*}
d\left(\widetilde{\varphi}\cdot\widetilde{\nu}\right) 
& = & 
\sum_{j=1}^n 
\frac{\partial\left(\widetilde{\varphi}\cdot \widetilde{\nu}\right)}{\partial x_j}(x)
d x_j \\ 
{ } & = & 
\sum_{j=1}^n 
\frac{\partial\left(\widetilde{\varphi}\cdot \widetilde{\nu}\right)}{\partial x_j}(x)
\left(\sum_{i=1}^n \frac{\partial \left(x_j\circ \widetilde{\nu}^{-1}\right)}
{\partial \Theta_{\left(i, \widetilde{\nu}(x)\right)}}
\left(\widetilde{\nu}(x)\right) d \Theta_i\right) \\ 
{ } & = & 
\sum_{i=1}^n\left(\sum_{j=1}^n
\frac{\partial\left(\widetilde{\varphi}\cdot \widetilde{\nu}\right)}{\partial x_j}(x)
\frac{\partial \left(x_j\circ \widetilde{\nu}^{-1}\right)}
{\partial \Theta_{\left(i, \widetilde{\nu}(x)\right)}}
\left(\widetilde{\nu}(x)\right) 
\right)d \Theta_i.    
\end{eqnarray*}
Let $\widetilde{\Omega} : N\to T^*S^n$ be the mapping with the form 
$\widetilde{\Omega}(x)=\left(\widetilde{\nu}(x), \widetilde\omega(x)\right)$ 
such that 
$\widetilde\omega$ is the creator for $\widetilde{f}$.     
Then, by the above calculation, 
$\widetilde\omega|_U$ must have the following form.   
\[
\widetilde\omega|_U(x)= 
\sum_{i=1}^n\left(\sum_{j=1}^n
\frac{\partial\left(\widetilde{\varphi}\cdot \widetilde{\nu}\right)}{\partial x_j}(x)
\frac{\partial \left(x_j\circ \widetilde{\nu}^{-1}\right)}
{\partial \Theta_{\left(i, \widetilde{\nu}(x)\right)}}\left(\widetilde{\nu}(x)\right) 
\right)
d \Theta_i.      
\]
Hence, by {\color{black}Theorem \ref{theorem1} (2)}, we have the following.      
\begin{lemma}\label{lemma9}
At a regular point $x_0\in N$ of $\widetilde{\nu}$, the equality 
$\widetilde{f}_1(x_0)=\widetilde{f}_2(x_0)$ holds.  
\end{lemma}
Let $x_0\in N$ be a singular point of $\widetilde{\nu}$.    
Then, since we have assumed that 
the set of regular points of $\widetilde{\nu}$ is dense, 
there exists a point-sequence $\left\{y_i\right\}_{i=1, 2, \ldots}\subset N$ 
such that 
$y_i$ is a regular point of $\widetilde{\nu}$ 
for any $i\in \mathbb{N}$ and $\lim_{i\to \infty}y_i=x_0$.     
Then, by Lemma \ref{lemma9}, we have  
\[
\widetilde{f}_1(x_0)=\widetilde{f}_1\left(\lim_{i\to \infty}y_i\right)
= \lim_{i\to \infty}\widetilde{f}_1(y_i)=\lim_{i\to \infty}\widetilde{f}_2(y_i)
= \widetilde{f}_2\left(\lim_{i\to \infty}y_i\right)=\widetilde{f}_2(x_0).   
\]
Thus, we have the following.   
\begin{lemma}\label{lemma10}
Even at a singular point $x_0\in N$ of $\widetilde{\nu}$, the equality 
$\widetilde{f}_1(x_0)=\widetilde{f}_2(x_0)$ holds.  
\end{lemma}
\hfill $\Box$
\par 
\smallskip 
\noindent 
\underline{\it Proof of \lq\lq only if\rq\rq\, part}. \quad  
Suppose that the hyperplane 
$\mathcal{H}_{\left(\widetilde{\varphi}, \widetilde{\nu}\right)}$ 
is creative and the set of regular points of $\widetilde{\nu}$ is not dense in $N$.   
Then, there exists an open set $U$ of $N$ such that any point $x\in U$ 
is a singular point of $\widetilde{\nu}$.   Then, there exist 
an integer $k$ $(0\le k <n)$ 
and an open set $U_k$ such that  $U_k\subset U$ 
and the rank of $\widetilde{\nu}$ at $x$ is $k$ 
for any  $x\in U_k$.   
Let $x_0$ be a point of $U_k$.   
We may assume that $U_k$ is sufficiently small open neighborhood of $x_0$.     
Then, by the rank theorem (for the rank theorem, 
see for example \cite{brocker}),    
we have the following.    
\begin{lemma}\label{lemma11}
There exist functions $\eta_1, \ldots, \eta_k: N\to \mathbb{R}$ such that 
the following three hold.   
\begin{enumerate}
\item[(1)]\quad For any $i$ $(1\le i\le n)$, $\eta_i(x)=0$ if $x\not\in U_k$.   
\item[(2)]\quad There exists an $i$ $(1\le i\le n)$ 
such that $\eta_i\left(x_0\right)\ne 0$.   
\item[(3)]\quad  
The following equality holds for any $x\in N$.   
\[
\sum_{i=1}^n
\eta_i(x) 
d\left(\Theta_i\circ \widetilde{\nu}\right)
=0.    
\] 
\end{enumerate}
\end{lemma}
Since we have assumed that 
$\mathcal{H}_{\left(\widetilde{\varphi}, \widetilde{\nu}\right)}$ is creative, 
there exists a mapping $\widetilde{\Omega}: N\to T^*S^n$ with he form 
$\widetilde{\Omega}(x)=\left(\widetilde{\nu}(x), \widetilde\omega(x)\right)$ 
such that 
$d\left(\widetilde{\varphi}\cdot \widetilde{\nu}\right)=\widetilde\omega$.        
By Lemma \ref{lemma11}, the following holds.   
\begin{lemma}\label{lemma12}
For any function $\alpha: N\to \mathbb{R}$ and any $x\in N$, 
the following equality holds as germs of $1$-form at $x$.  
\[ 
d\left(\widetilde{\varphi}\cdot \widetilde{\nu}\right)= \widetilde\omega(x) +
\alpha(x)
\sum_{i=1}^n
\eta_i(x) 
d\left(\Theta_i\circ \widetilde{\nu}\right).  
\]
\end{lemma}
Therefore, by 
{\color{black}Theorem \ref{theorem1} (2), }
uncountably many distinct 
envelopes $\widetilde{f}$ are    
created by the same hyperplane family 
$\mathcal{H}_{\left(\widetilde{\varphi}, \widetilde{\nu}\right)}$.     
\hfill $\Box$  
\section{Examples}\label{section2}
\begin{example}[Uniform spin of affine tangent lines]\label{example1}\quad 
\begin{enumerate}
\item[(1)] Let $\alpha: \mathbb{R}\to \mathbb{R}$ be a non-constant function.    
{\color{black}Notice that $\alpha$ is of class $C^\infty$ as stated 
at the top of Section \ref{section1}}.   
Let $\widetilde{\varphi}: \mathbb{R}\to \mathbb{R}^2$ be the 
mapping defined by 
$\widetilde{\varphi}(t)=(\alpha(t),0)$.    
Let $\widetilde{\nu}: \mathbb{R}\to S^1$ be the constant mapping 
$\widetilde{\nu}(t)=(0,1)$.    
For any fixed $\theta_0\in \mathbb{R}$, 
let $R_{\theta_0}:\mathbb{R}^2\to \mathbb{R}^2$ be 
the linear mapping representing the rotation through angle $\theta_0$.    
Set $\widetilde{\nu}_{\theta_0}(t)=R_\theta\circ \widetilde{\nu}(t)=
\left(-\sin \theta_0, \cos \theta_0\right)$ and 
$\widetilde{\gamma}_{\theta_0}(t)=\widetilde{\varphi}(t)\cdot 
\widetilde{\nu}_{\theta_0}(t)=
-\alpha(t)\sin \theta_0$.      
Figure is depicted in Figure \ref{example1(1)}.     
\begin{figure}
\begin{center}
\includegraphics[width=15cm]
{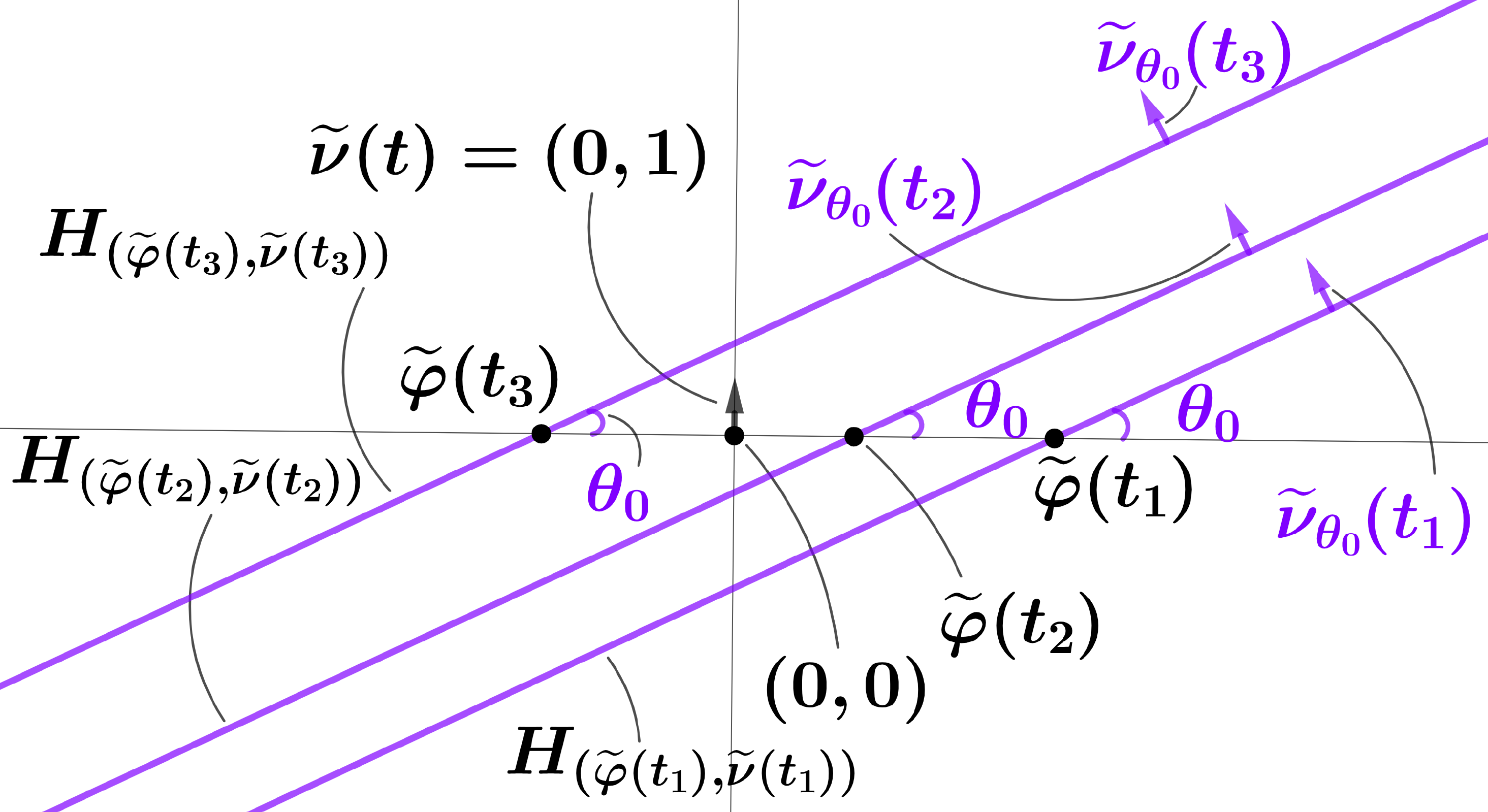}
\caption{
Figure for Example \ref{example1} (1).       
} 
\label{example1(1)}
\end{center}
\end{figure}
It follows $d\left(\Theta\circ \widetilde{\nu}_{\theta_0}\right)\equiv 0$ and 
$d\gamma_{\theta_0}=-\sin \theta_0 d\alpha$.    
Since $\alpha$ is non-constant, 
there exists a regular point of $\alpha${\color{black}, 
that is to say, there exists a $t\in \mathbb{R}$ such that $\alpha'(t)\ne 0$}.   
Therefore, by Theorem \ref{theorem1}, the line family 
$\mathcal{H}_{\left(\widetilde{\varphi}, \widetilde{\nu}_{\theta_0}\right)}$ creates 
an envelope if and only if $\theta_0\in \pi\mathbb{Z}$.    
Suppose that $\theta_0\in \pi\mathbb{Z}$.       
In this case, by Theorem \ref{theorem2}, uncountably many distinct envelope 
$\widetilde{f}: \mathbb{R}\to \mathbb{R}^2$ can be created by the given 
line family 
$\mathcal{H}_{\left(\widetilde{\varphi}, \widetilde{\nu}_{\theta_0}\right)}$.   
Let $\beta: \mathbb{R}\to \mathbb{R}$ be a function.     
Since  $d\left(\Theta\circ \widetilde{\nu}_{\theta_0}\right)\equiv 0$ and 
$d\gamma_{\theta_0}\equiv 0$ in this case, the $1$-form 
$t\mapsto \beta(t)d \left(\Theta\circ \widetilde{\nu}_{\theta_0}\right)$ 
along $\widetilde{\nu}_{\theta_0}$ may be 
a creator $\widetilde\omega$ for the line family.   
By {\color{black}Theorem \ref{theorem1} (2), }
the envelope $\widetilde{f}$ 
has the following form.  
\[
\widetilde{f}(t) 
=
\widetilde\omega(t)+
{\color{black}{\left(
{\color{black}{
\widetilde{\gamma}_{\theta_0}(t)\cdot \widetilde{\nu}_{\theta_0}(t) 
}}
\right)\widetilde{\nu}_{\theta_0}(t) 
}}
= 
(
{\color{black}{
{\color{black}\pm}
}}
\beta(t), 0)+(0,0)
=
(
{\color{black}{
{\color{black}\pm}
}}
\beta(t), 0),  
\]
where {\color{black}double sign should be read in the same order and} 
$\beta(t)d \left(\Theta\circ \widetilde{\nu}_{\theta_0}\right)${\color{black},} 
$\beta(t)R_{\frac{\pi}{2}}\circ \widetilde{\nu}_{\color{black}{\theta_0}}(t)$ 
are identified 
(both are denoted by the same symbol $\widetilde{\omega}(t)$).
\par 
Set $F_{\theta_0}\left(X_1, X_2, t\right)=
\left(X_1-\alpha(t), X_2\right)\cdot \widetilde{\nu}_{\theta_0}(t)$.    
Suppose that $\theta_0\not\in \pi\mathbb{Z}$.    
In this case, the classical common definition of envelope 
$\mathcal{D}$ relative to $F_{\theta_0}$ is as follows.   
\[
\mathcal{D}=\left\{\left(X_1, X_2\right)\, |\, \exists t \mbox{ s.t. }
\alpha'(t)=0, X_1=\cot \theta_0 X_2 +\alpha(t)\right\}.   
\] 
Therefore, in this case, $\mathcal{D}=E_1=E_2=\emptyset$ 
if and only if $\alpha$ is non-singular.    
Suppose that $\theta_0\in \pi\mathbb{Z}$.    Then, 
\[
\mathcal{D}=\left\{\left(X_1, X_2\right)\, |\, X_2=0\right\}.   
\] 
Therefore, in this case, $E_1=E_2=\mathcal{D}$ if and only if $\beta$ 
is surjective.    
\item[(2)] 
Let $\widetilde{\nu}: \mathbb{R}\to S^1$ be the mapping given by 
$\widetilde{\nu}(t)=(\cos t, \sin t)$.    
Set $\widetilde{\nu}_{\theta_0}=R_{\theta_0}\circ \widetilde{\nu}$,  
where $R_{\theta_0}$ is the rotation defined in the above example.  
Then, since $\frac{d\left(\Theta\circ \widetilde{\nu}_{\theta_0}\right)}{dt}(t)=1$,  
it follows $d \left(\Theta\circ \widetilde{\nu}_{\theta_0}\right)=dt$.   
Thus, by Theorem \ref{theorem1}{\color{black}(1)} and Theorem \ref{theorem2}, 
for any $\widetilde{\varphi}: \mathbb{R}\to \mathbb{R}^2$ the line family 
$\mathcal{H}_{\left(\widetilde{\varphi}, \widetilde{\nu}_{\theta_0}\right)}$ 
creates a unique envelope $\widetilde{f}_{\theta_0}$.   
For any $\widetilde{\varphi}: \mathbb{R}\to \mathbb{R}^2$, set 
$\widetilde{\gamma}_{\theta_0}(t)=
\widetilde{\varphi}(t)\cdot \widetilde{\nu}_{\theta_0}(t)$.   
Since $d \widetilde{\gamma}_{\theta_0}=
\frac{d \widetilde{\gamma}_{\theta_0}}{dt}(t)
d \left(\Theta\circ \widetilde{\nu}_{\theta_0}\right)$, by 
{\color{black}Theorem \ref{theorem1} (2), }  
it follows 
\begin{eqnarray*}
\widetilde{f}(t)
& = & 
\frac{d \widetilde{\gamma}_{\theta_0}}{dt}(t)
R_{\pi/2}\circ \widetilde{\nu}_{\theta_0}\left(t\right)
+\widetilde{\gamma}_{\theta_0}(t)\widetilde{\nu}_{\theta_0}(t) \\ 
{ } & = &  
\frac{d \widetilde{\gamma}_{\theta_0}}{dt}(t)
R_{\pi/2}\circ \widetilde{\nu}_{\theta_0}\left(t\right) + 
\widetilde{\gamma}_
{\theta_0}(t)\left(\cos\left(t+\theta_0\right), \sin\left(t+\theta_0\right)\right), 
\end{eqnarray*}
where the $1$-form $d\left(\Theta\circ \widetilde{\nu}\right)$ 
and the vector field  
$R_{\pi/2}\circ \widetilde{\nu}_{\theta_0}\left(t\right)$ are identified.   
Let $\alpha: \mathbb{R}\to \mathbb{R}$ be a function and set 
$\widetilde{\varphi}(t)=
\widetilde{\nu}(t)+\alpha(t)R_{\pi/2}\circ \widetilde{\nu}_{\theta_0}(t)$.    
\begin{figure}
\begin{center}
\includegraphics[width=15cm]
{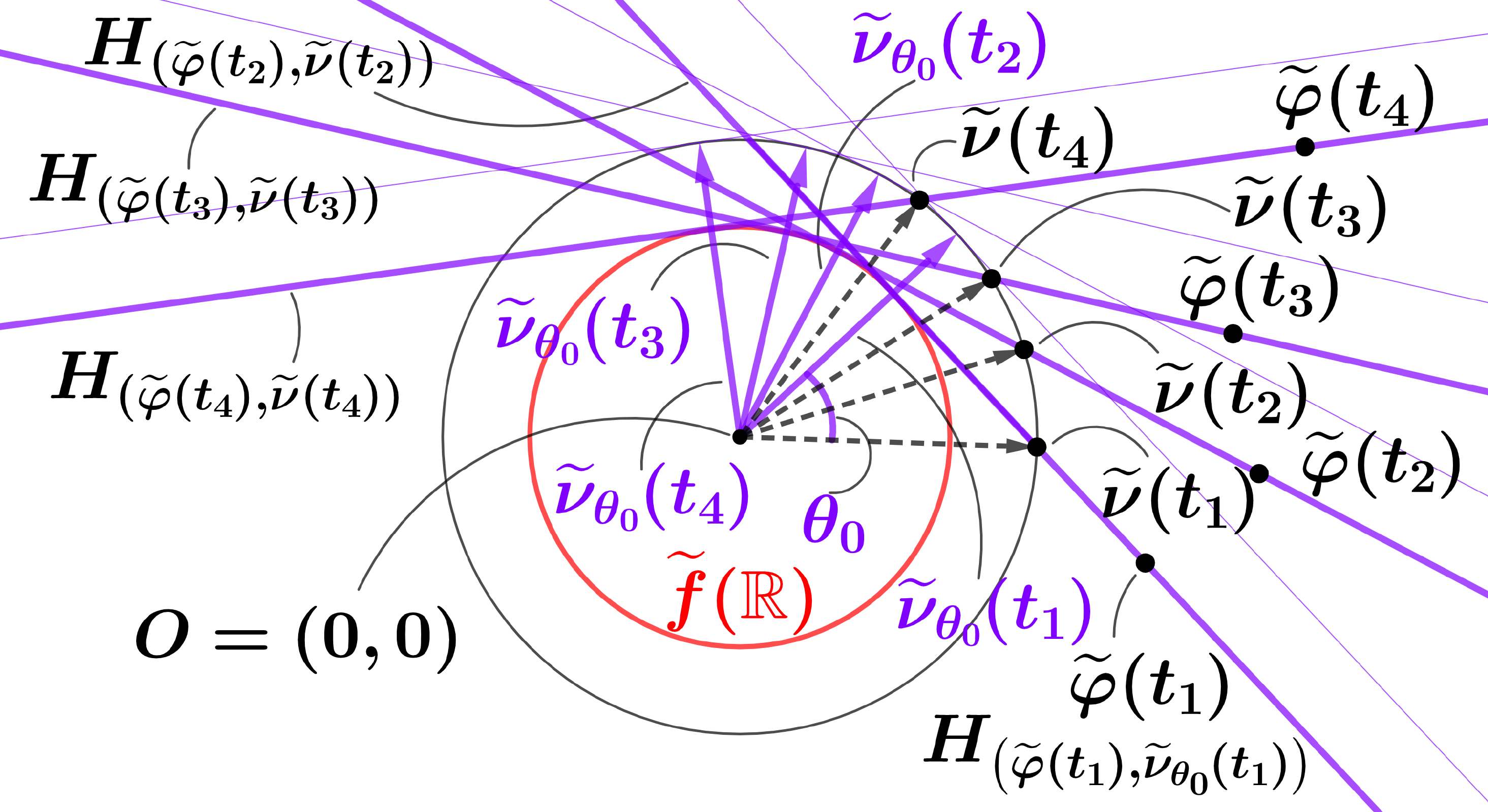}
\caption{
Figure for Example \ref{example1} (2).       
} 
\label{example1(2)}
\end{center}
\end{figure}
Then, it follows $\frac{d \widetilde{\gamma}_{\theta_0}}{dt}(t)\equiv 0$.   
Thus, as expected, the envelope created by 
the line family 
$\mathcal{H}_{\left(\widetilde{\varphi}, \widetilde{\nu}_{\theta_0}\right)}$ 
in this case 
is actually the circle 
with radius $|c|$ centered at the origin, where 
$c=\widetilde{\gamma}_
{\theta_0}(t)=\cos\theta_0$ (see Figure \ref{example1(2)}).   
\item[(3)] 
Let $\widetilde{\nu}: \mathbb{R}\to S^1$ be the mapping defined by  
$\widetilde{\nu}(t)=\frac{1}{\sqrt{1+9t^4}}\left(-3t^2, 1\right)$.   
Set $\widetilde{\nu}_{\theta_0}=R_{\theta_0}\circ \widetilde{\nu}$ 
where $R_{\theta_0}$ is as above.   
Let $\alpha: \mathbb{R}\to \mathbb{R}$ be a function and set 
$\widetilde{\varphi}_{\theta_0}(t)=
(t, t^3)+\alpha(t)R_{\pi/2}\circ \widetilde{\nu}_{\theta_0}(t)$.     
Set $\widetilde{\gamma}_{\theta_0}(t)=
\widetilde{\varphi}_{\theta_0}(t)\cdot \widetilde{\nu}_{\theta_0}(t)$. 
It is easily seen that $0$ is a singular point of 
$\widetilde{\gamma}_{\theta_0}$ if and only if 
$\theta_0\in \pi\mathbb{Z}$.  On the other hand, by calculation, we have     
$\frac{d\left(\Theta\circ \widetilde{\nu}_{\theta_0}\right)}{dt}(t)=
\frac{6t}{1+9t^4}$ and thus 
$0$ is a unique singular point of $\widetilde{\nu}_{\theta_0}$ for any $\theta_0$.    
Therefore, by Theorem \ref{theorem1}, the hyperplane family 
$\mathcal{H}_{\left(\widetilde{\varphi}, \widetilde{\nu}_{{}_\theta}\right)}$ 
does not create an envelope if $\theta_{\color{black}0}\not\in \pi\mathbb{Z}$.   
\par 
Next, suppose that $\theta_0\in \pi\mathbb{Z}$.    
Then, calculations show 
\[  
\quad\quad d\left(\widetilde{\gamma}_{\theta_0}\right)
=
\frac{\color{black}\mp(6t^2+18t^6)}{(1+9t^4)^{\frac{3}{2}}}dt 
=
 \frac{\color{black}\mp(t+3t^5)}{\sqrt{1+9t^4}}
\frac{d\left(\Theta\circ \widetilde{\nu}_{\theta_0}\right)}{dt}(t)dt 
= 
\frac{\color{black}\mp(t+3t^5)}{\sqrt{1+9t^4}}
d\left(\Theta\circ \widetilde{\nu}_{\theta_0}\right){\color{black},} 
\] 
{\color{black}where double sign should be read in the same order}.  
Set $\widetilde\omega(t)=\frac{\color{black}\mp(t+3t^5)}
{\sqrt{1+9t^4}}d\left(\Theta\circ \widetilde{\nu}_{\theta_0}\right)$.   
By Theorem \ref{theorem1} {\color{black}and} Theorem \ref{theorem2},  
the hyperplane family  
$\mathcal{H}_{(\widetilde{\varphi}, \widetilde{\nu}_{\theta_0})}$ 
creates a unique envelope with the desired form 
\begin{eqnarray*}
\widetilde{f}(t) & = & 
\widetilde\omega(t)+
\widetilde{\gamma}_{\theta_0}(t)\widetilde{\nu}_{\theta_0}(t) \\ 
{ } & = & \frac{\color{black}\mp(t+3t^5)}{1+9t^4}
\left({\color{black}\mp} 1, 
{\color{black}\mp} 3t^2\right)
{\color{black}\mp} 
\frac{2t^3}{1+9t^4}({\color{black}\mp} 3t^2, {\color{black}\pm} 1) \\ 
{ } & = & \frac{1}{1+9t^4}\left(t+3t^5+6t^5, 3t^3+9t^7-2t^3\right)\\ 
{ } & = & \left(t, t^3\right), 
\end{eqnarray*} 
where for each $t\in \mathbb{R}$ 
the cotangent vector $\frac{\color{black}\mp(t+3t^5)}{\sqrt{1+9t^4}}
d\left(\Theta\circ \widetilde{\nu}_{\theta_0}\right)$ 
and the vector 
$\frac{\color{black}\mp(t+3t^5)}{\sqrt{1+9t^4}}R_{\pi/2}
\circ \widetilde{\nu}_{\theta_0}(t)$ 
in the vector space $\mathbb{R}^2$
are identified (see Figure \ref{example1(3)}).     
\begin{figure}
\begin{center}
\includegraphics[width=15cm]
{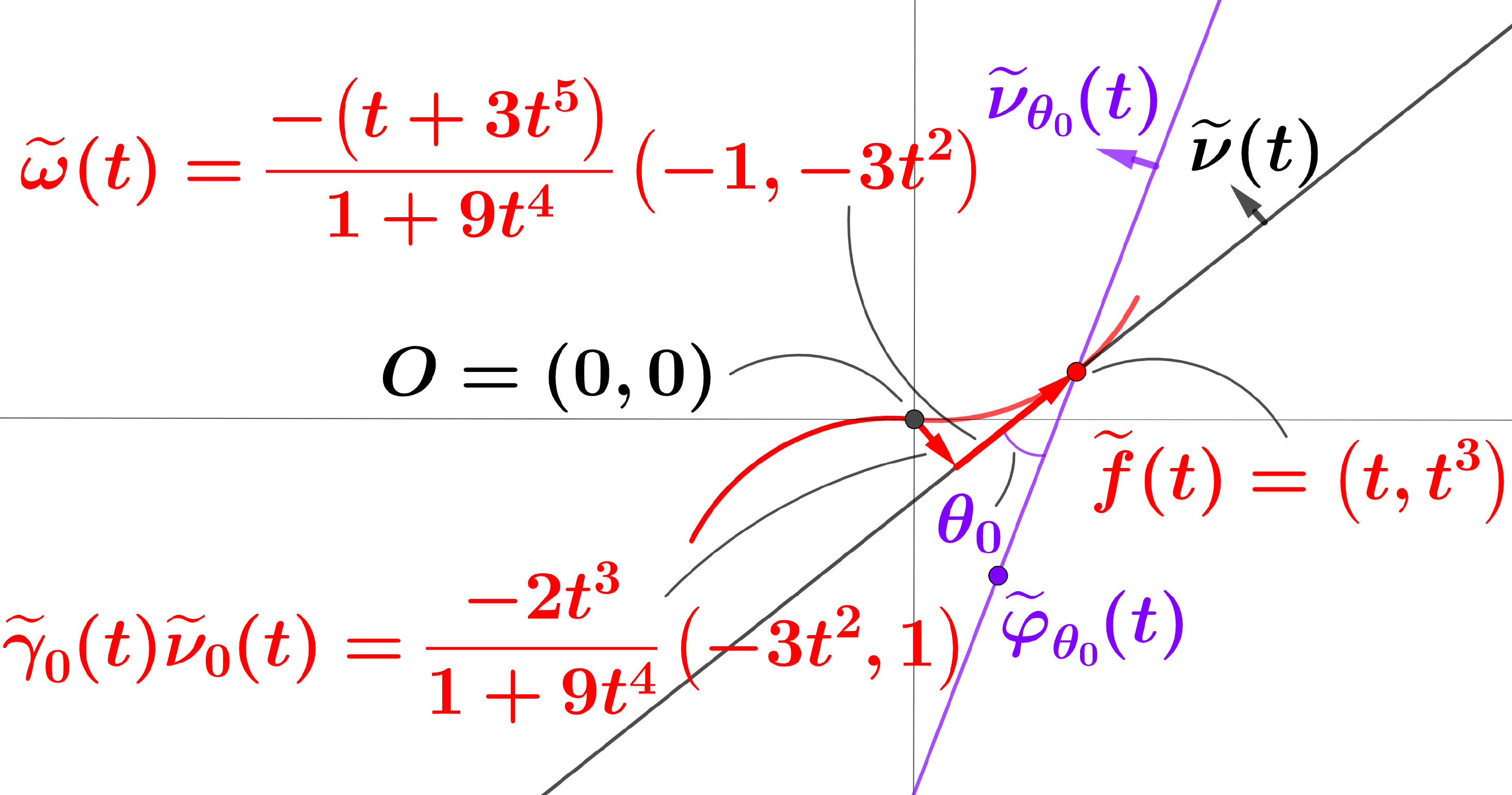}
\caption{
Figure for Example \ref{example1} (3) 
{\color{black}in the case $\theta_0\in 2\pi\mathbb{Z}$}.       
} 
\label{example1(3)}
\end{center}
\end{figure}
\par 
Set $U=\mathbb{R}-\{0\}$.    
It is easily seen that $\widetilde{\nu}_{\theta_0}|_U$ is non-singular 
even in the case $\theta_0\not\in \pi\mathbb{Z}$.       
Hence, by Theorem \ref{theorem1} and Theorem \ref{theorem2}, 
the hyperplane family  
$\mathcal{H}_{\left(\widetilde{\varphi}|_U, 
\widetilde{\nu}_{\theta_0}|_U\right)}$ creates 
a unique envelope $\widetilde{f}_{\theta_0}: U\to \mathbb{R}^2$ 
even when $\theta_0\not\in \pi\mathbb{Z}$ and 
$\lim_{t\to 0}\Vert\widetilde{f}_{\theta_0}(t)\Vert=\infty$ 
when $\theta_0\not\in \pi\mathbb{Z}$.   
\item[(4)] 
Let $\widetilde{\nu}: \mathbb{R}\to S^1$ be the mapping defined by  
$\widetilde{\nu}(t)=\frac{1}{\sqrt{4+25t^6}}\left(-5t^3, 2\right)$.   
Set $\widetilde{\nu}_{\theta_0}=R_{\theta_0}\circ \widetilde{\nu}$ 
where $R_{\theta_0}$ 
is as above.   
Let $\alpha: \mathbb{R}\to \mathbb{R}$ be a function and set 
$\widetilde{\varphi}_{\theta_0}(t)=
(t^2, t^5)+\alpha(t)R_{\pi/2}\circ \widetilde{\nu}_{\theta_0}(t)$.     
Set $\widetilde{\gamma}_{\theta_0}(t)=
\widetilde{\varphi}_{\theta_0}(t)\cdot \widetilde{\nu}_{\theta_0}(t)=
\frac{-3t^5\cos \theta_0-2t^2\sin\theta_0-5t^{\color{black}8}
\sin\theta_0}{\sqrt{4+25t^6}}$.    
By calculation, we have     
$\frac{d\left(\Theta\circ \widetilde{\nu}_{\theta_0}\right)}{dt}(t)=
\frac{30t^2}{4+25t^6}$.  
Therefore, 
the hyperplane family 
$\mathcal{H}_{\left(\widetilde{\varphi}, \widetilde{\nu}_{{}_\theta}\right)}$ 
is not creative if $\theta\not\in \pi\mathbb{Z}$ and 
it creates no envelope in this case by Theorem \ref{theorem1}.   
\par 
Next, suppose that $\theta_0\in \pi\mathbb{Z}$.    
Then, calculation shows  
\begin{eqnarray*}
d\left(\widetilde{\gamma}_{\theta_0}\right)
& = & 
\frac{\color{black}\mp 30t^2\left(2t^2+5t^8\right)}{(4+25t^6)\sqrt{4+25t^6}}dt  \\ 
& = & 
 \frac{\color{black}\mp(2t^2+5t^8)}{\sqrt{4+25t^6}}
\frac{d\left(\Theta\circ \widetilde{\nu}_{\theta_0}\right)}{dt}(t)dt  
 = 
\frac{\color{black}\mp(2t^2+5t^8)}{\sqrt{4+25t^6}}
d\left(\Theta\circ \widetilde{\nu}_{\theta_0}\right){\color{black},}  
\end{eqnarray*}
{\color{black}where double sign should be read in the same order}.  
Therefore, 
the hyperplane family 
$\mathcal{H}_{\left(\widetilde{\varphi}, \widetilde{\nu}_{{}_\theta}\right)}$ 
is creative.    
Set $\widetilde\omega(t)=
\frac{\color{black}\mp(2t^2+5t^8)}{\sqrt{4+25t^6}}.    
d\left(\Theta\circ \widetilde{\nu}_{\theta_0}\right)$. 
By Theorem \ref{theorem1} {\color{black}and} Theorem \ref{theorem2},  
$\mathcal{H}_{(\widetilde{\varphi}, \widetilde{\nu}_{\theta_0})}$ 
creates a unique envelope with the desired form 
\begin{eqnarray*}
\widetilde{f}(t) & = & 
\widetilde\omega(t)+
\widetilde{\gamma}_{\theta_0}(t)\widetilde{\nu}_{\theta_0}(t) \\ 
{ } & = & \frac{\color{black}\mp(2t^2+5t^8)}{4+25t^6}
\left({\color{black}\mp}2, {\color{black}\mp}5t^3\right)
+ \frac{{\color{black}\mp}3t^5}{4+25t^6}({\color{black}\mp}5t^3, {\color{black}\pm}2) \\ 
{ } & = & \frac{1}{4+25t^6}\left(4t^2+10t^8+15t^8, 10t^5+25t^{11}-6t^5\right)\\ 
{ } & = & \left(t^2, t^5\right),  
\end{eqnarray*}   
where for each $t\in \mathbb{R}$  
the cotangent vector $\frac{\color{black}\mp(2t^2+5t^8)}{\sqrt{4+25t^6}}
d\left(\Theta\circ \widetilde{\nu}_{\theta_0}\right)$ and the vector 
$\frac{\color{black}\mp(2t^2+5t^8)}{\sqrt{4+25t^6}}R_{\pi/2}\circ 
\widetilde{\nu}_{\theta_0}(t)$ in the 
vector space $\mathbb{R}^2$ are identified (see Figure \ref{example1(4)}).   
\begin{figure}
\begin{center}
\includegraphics[width=15cm]
{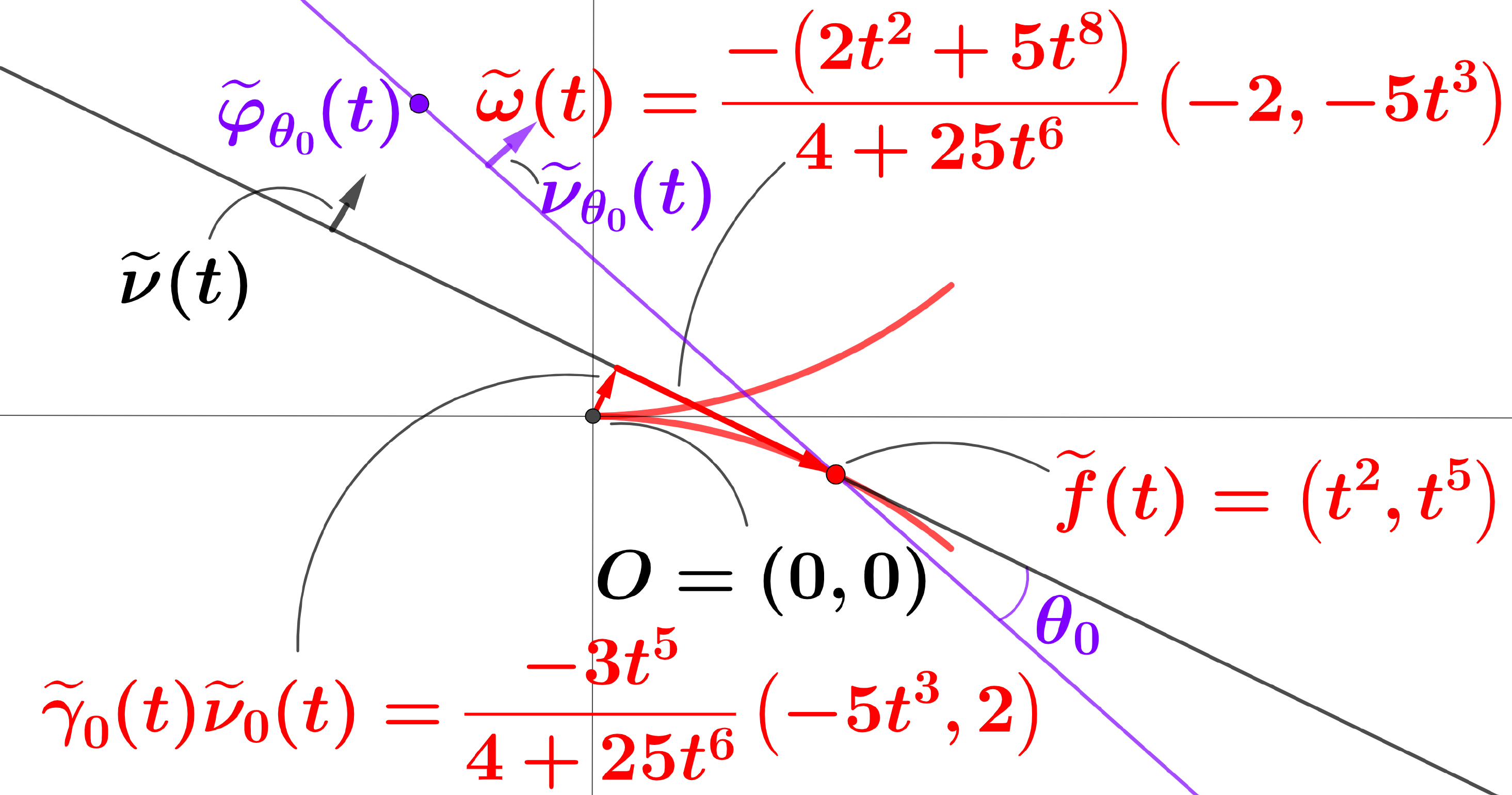}
\caption{
Figure for Example \ref{example1} (4) {\color{black}in the case 
$\theta_0\in 2\pi\mathbb{Z}$}.       
} 
\label{example1(4)}
\end{center}
\end{figure}
In the case $\theta_0=0$, consider the mapping 
$\widetilde{\Omega}: \mathbb{R}\to T^*S^1$ 
given in Definition \ref{definition2} and   
$\Omega: \mathbb{R}\to J^1\left(S^1, \mathbb{R}\right)$ given in Remark 
\ref{remark1}(1).   
Namely, consider the following two mappings.  
\begin{eqnarray*}
\widetilde{\Omega}(t) & = & 
\left(\frac{1}{\sqrt{4+25t^6}}\left({\color{black}\mp}5t^3, {\color{black}\pm}2\right), 
\frac{\color{black}\mp30t^2(2t^2+5t^8)}
{\color{black}\left(4+25t^6\right)^{\frac{3}{2}}}\right),  \\ 
\Omega(t) & = & \left(\frac{1}{\sqrt{4+25t^6}}
\left({\color{black}\mp}5t^3, {\color{black}\pm}2\right), 
\frac{{\color{black}\mp}3t^5}{\sqrt{4+25t^6}}, 
\frac{\color{black}\mp30t^2(2t^2+5t^8)}{\color{black}\left(4+25t^6\right)^{\frac{3}{2}}}
\right).  
\end{eqnarray*}   
Since $d \left(\widetilde{\gamma}_{\theta_0}\right)
=\frac{\color{black}\mp(2t^2+5t^8)}{\sqrt{4+25t^6}}
d\left(\Theta\circ \widetilde{\nu}_{\theta_0}\right)$, 
the map-germ of $\Omega$ at any $t$ is 
nothing but an opening of the map-germ 
$\widetilde{\Omega}: (\mathbb{R}, t)\to T^*S^1$.     
At $t=0$, the map-germ of each of them is not immersive 
and has singular images.    
\par 
Set $U=\mathbb{R}-\{0\}$.    
It is easily seen that $\widetilde{\nu}_{\theta_0}|_U$ is non-singular 
even in the case $\theta_0\not\in \pi\mathbb{Z}$.       
Hence, by Theorem \ref{theorem1} and Theorem \ref{theorem2}, 
the hyperplane family  
$\mathcal{H}_{\left(\widetilde{\varphi}|_U, 
\widetilde{\nu}_{\theta_0}|_U\right)}$ creates 
a unique envelope $\widetilde{f}_{\theta_0}: U\to \mathbb{R}^2$ 
even when $\theta_0\not\in \pi\mathbb{Z}$ and 
$\lim_{t\to 0}\Vert\widetilde{f}_{\theta_0}(t)\Vert=\infty$ 
when $\theta_0\not\in \pi\mathbb{Z}$.   
\end{enumerate}  
\end{example}
{\color{black}{
\begin{example}[Unit speed curves]\label{example0}\quad 
\begin{enumerate}
\item[(1)] Let ${\bf r}: \mathbb{R}\to \mathbb{R}^2$ be a unit speed curve.   
As usual, set ${\bf t}(s)={\bf r}'(s)$ and ${\bf n}(s)$ is defined 
from {\bf t}(s) by rotating anticlockwise through $\frac{\pi}{2}$.    
The Serret-Frenet formulas for the plane curve ${\bf r}$ is as follows.    
\[
\left\{
\begin{array}{cccc}   
{\bf t}'(s) & = & { } & \kappa(s){\bf n}(s) \\
{\bf n}'(s) & = & -\kappa(s){\bf t}(s). & { }  
\end{array} 
\right. 
\]
Set $\widetilde{\varphi}={\bf r}$ and $\widetilde{\nu}={\bf n}$.    
Then, the line family 
$\mathcal{H}_{\left(\widetilde{\varphi}, \widetilde{\nu}\right)}
= 
\mathcal{H}_{\left({\bf r}, {\bf n}\right)}
$ 
is the affine tangent line family of the curve ${\bf r}$.    
In this case, the correspondence ${\bf r}\mapsto 
\mathcal{H}_{\left({\bf r}, {\bf n}\right)}$ may be regarded as 
the Legendre transformation of the given curve ${\bf r}$.  
Set $\widetilde{\gamma}(s)=\widetilde{\varphi}(s)\cdot \widetilde{\nu}(s)$.   
Then,  
\[
\widetilde{\gamma}'(s)={\bf r}(s)\cdot 
\left(-\kappa(s){\bf t}(s)\right)
= 
{\color{black}-}
\left({\bf r}(s)\cdot {\bf t}(s)\right)
\left(\Theta_{{\bf t}}\circ \widetilde{\nu}\right)'(s){\color{black},} 
\]
{\color{black}where $\widetilde{\nu}(s)=
\left(\cos \Theta_{{\bf t}}\circ \widetilde{\nu}(s), 
\sin \Theta_{{\bf t}}\circ \widetilde{\nu}(s)\right)$}.   
Therefore, by Theorem \ref{theorem1}, 
the line family $\mathcal{H}_{\left(\widetilde{\varphi}, \widetilde{\nu}\right)}$ 
creates an envelope.   
\par 
Suppose that the set of regular points of $\widetilde{\nu}$ 
is dense, that is to say, 
the set $\{s\in \mathbb{R}\; |\; \kappa(s)\ne 0\}$ is dense.   
Then, by Theorem \ref{theorem2}, the created envelopes are unique.   
By {\color{black}Theorem \ref{theorem1}, } 
the unique envelope is as follows  
(see Figure \ref{example3(1)}). 
\begin{eqnarray*}
\widetilde{f}(s) & = & \widetilde{\omega}(s)
+\widetilde{\gamma}(s)\cdot \widetilde{\nu}(s) \\ 
{ } & = & \left({\bf r}(s)\cdot {\bf t}(s)\right){\bf t}(s)
+\left({\bf r}(s)\cdot {\bf n}(s)\right){\bf n}(s) \\ 
{ } & = & {\bf r}(s).   
\end{eqnarray*}   
{\color{black}Notice that if there is a point $s\in \mathbb{R}$ 
such that $\kappa(s)=0$, 
then the full discriminant of the line family is different 
from the unique desired envelope since the full discriminant includes 
the affine tangent line at $s$.     
This is one of advantages of our method.}               
The correspondence 
\[
\mathcal{H}_{\left({\bf r}, {\bf n}\right)}\mapsto {\bf r}
\]
may be regarded as the inverse Legendre transformation 
for plane curves.       
\par 
Next, suppose that the set of regular points of $\widetilde{\nu}$ is not dense.  
Then, there exists an open interval $(a, b)$ such that 
$\kappa(s)=0$ for any $s\in (a,b)$.     
Then, for any $s\in (a,b)$ and any function 
$\alpha: \mathbb{R}\to \mathbb{R}$ 
such that $\alpha(\mathbb{R}-(a, b))=\{0\}$, 
it follows 
\[
\widetilde{\gamma}'(s)=
\alpha(s)\left(\Theta_{{\bf t}}\circ \widetilde{\nu}\right)'(s).
\] 
By {\color{black}Theorem \ref{theorem1}, }
\begin{eqnarray*}
\widetilde{f}(s) & = & \widetilde{\omega}(s)
+\widetilde{\gamma}(s)\cdot \widetilde{\nu}(s) \\ 
{ } & = & 
\alpha(s){\bf t}(s)
+\left({\bf r}(s)\cdot {\bf n}(s)\right){\bf n}(s) \\ 
{ } & = & 
\left(\left(\alpha(s)-\left({\bf r}(s)\cdot{\bf t}(s)\right)\right)
+\left({\bf r}(s)\cdot{\bf t}(s)\right)\right){\bf t}(s)
+\left({\bf r}(s)\cdot {\bf n}(s)\right){\bf n}(s) \\ 
{ } & = & {\bf r}(s)+\beta(s){\bf t}(s),    
\end{eqnarray*}
where $\beta(s)=\alpha(s)-\left({\bf r}(s)\cdot{\bf t}(s)\right)$.   
Hence, in this case, 
the inverse Legendre transformation does not work well.    
\end{enumerate}
\begin{figure}[htbp]
    \centering
    \includegraphics[width=15cm]
{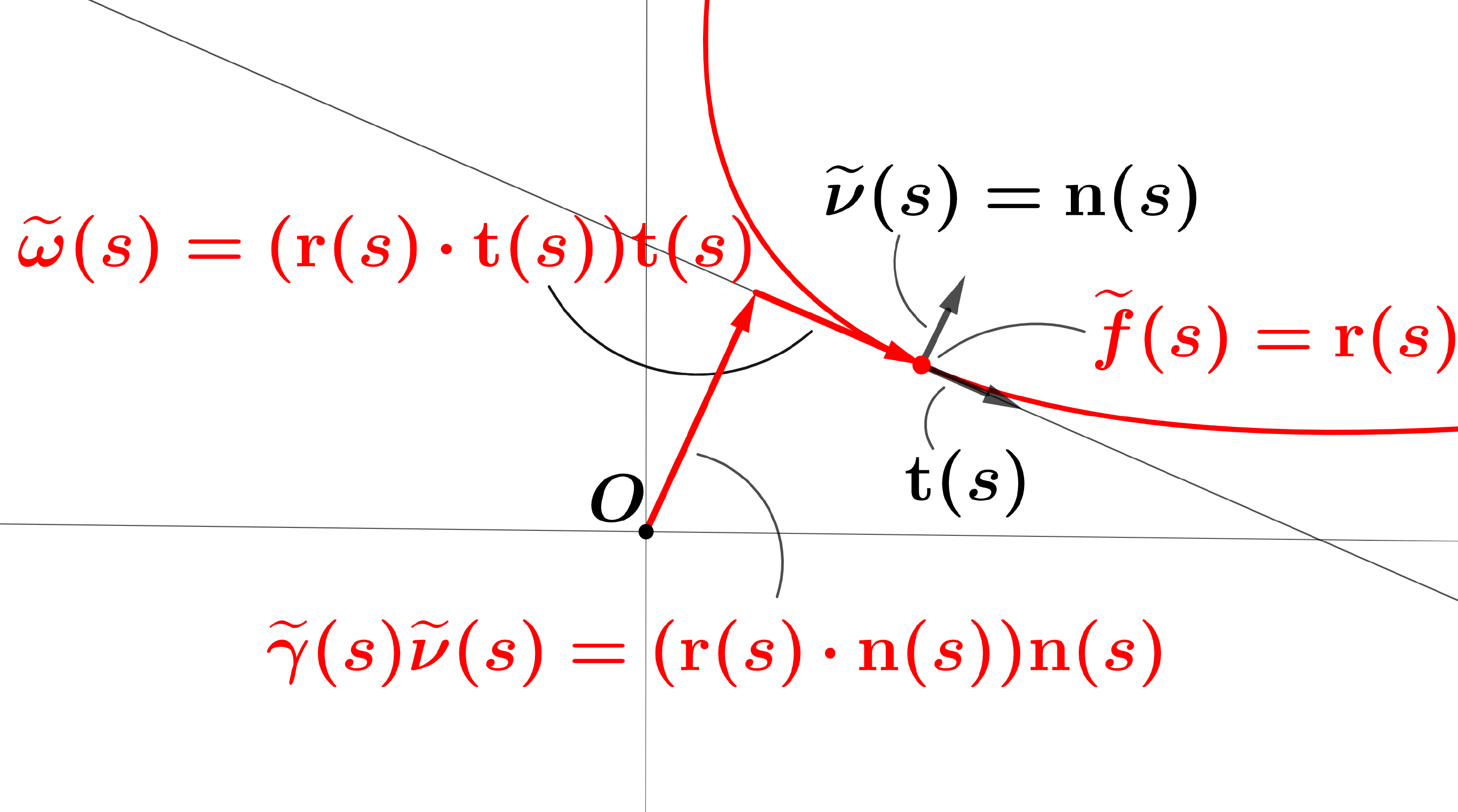}
    \caption{
Example \ref{example0} (1).  
}\label{example3(1)}
\end{figure}
\begin{enumerate}
\item[(2)] Let ${\bf r}: \mathbb{R}\to \mathbb{R}^3$ 
be a unit speed space curve.   
As usual, set ${\bf t}(s)={\bf r}'(s)$ and assume $||{\bf t}'(s)||> 0$ for any 
$s\in \mathbb{R}$ so that the principal normal vector ${\bf n}(s)$  
can be defined by ${\bf t}'(s)=||{\bf t}'(s)||{\bf n}(s)$.     As usual, the binormal 
vector ${\bf b}(s)$ is defined by 
$\det\left({\bf t}(s), {\bf n}(s), {\bf b}(s)\right)=1$.     
The Serret-Frenet formulas for the space curve ${\bf r}$ is as follows.      
\[
\left\{
\begin{array}{ccccc}
{\bf t}'(s) & = & { } & \kappa(s){\bf n}(s) & { } \\
{\bf n}'(s) & = & -\kappa(s){\bf t}(s)& { } & +\tau(s){\bf b}(s) \\ 
{\bf b}'(s) & = & { } & -\tau(s){\bf n}(s) & { } .     
\end{array}
\right. 
\] 
Define $\widetilde{\varphi}:\mathbb{R}^2\to \mathbb{R}^3$ and  
$\widetilde{\nu}:\mathbb{R}^2\to S^2$ 
by 
$\widetilde{\varphi}(s,u)={\bf r}(s)$ and $\widetilde{\nu}(s,u)={\bf b}(s)$ 
respectively.    
Then, the plane family 
$\mathcal{H}_{\left(\widetilde{\varphi}, \widetilde{\nu}\right)}
$ 
is the family of osculating planes of the space curve ${\bf r}$. 
Set $\widetilde{\gamma}(s,u)={\bf r}(s)\cdot {\bf b}(s)$.    
Then, all of the following six identities are clear.       
\[
\frac{\partial \widetilde{\gamma}}{\partial s}(s,u)  =  {\bf r}(s)\cdot 
\left(-\tau(s){\bf n}(s)\right), \quad 
 \frac{\partial \widetilde{\gamma}}{\partial u}(s,u)  =  0, 
\quad 
\frac{\partial \left(\Theta_{{\bf t}}\circ \widetilde{\nu}\right)}{\partial s}(s,u) =0, 
\quad 
\]
\[
\frac{\partial \left(\Theta_{{\bf t}}\circ \widetilde{\nu}\right)}{\partial u}(s,u) =0, 
\quad 
\frac{\partial \left(\Theta_{{\bf n}}\circ \widetilde{\nu}\right)}{\partial s}(s,u) 
=-\tau(s), 
\quad 
\frac{\partial \left(\Theta_{{\bf n}}\circ \widetilde{\nu}\right)}{\partial u}(s,u) =0.  
\] 
Therefore, we have the following.  
\begin{eqnarray*}
\frac{\partial \widetilde{\gamma}}{\partial s}(s,u) & = & 
\alpha_1(s,u)
\frac{\partial \left(\Theta_{{\bf t}}\circ \widetilde{\nu}\right)}{\partial s}(s,u)+ 
\left({\bf r}(s)\cdot {\bf n}(s)\right)
\frac{\partial \left(\Theta_{{\bf n}}\circ \widetilde{\nu}\right)}{\partial s}(s,u) 
,  \\
\frac{\partial \widetilde{\gamma}}{\partial u}(s,u)  
& = & 
\alpha_2(s,u)
\frac{\partial \left(\Theta_{{\bf t}}\circ \widetilde{\nu}\right)}{\partial u}(s,u)
+ 
\alpha_3(s,u)
\frac{\partial \left(\Theta_{{\bf n}}\circ \widetilde{\nu}\right)}{\partial u}(s,u),
\end{eqnarray*}
where $\alpha_1, \alpha_2, \alpha_3: \mathbb{R}^2\to \mathbb{R}$ 
are arbitrary functions.   
Thus, by Theorem \ref{theorem1}, 
the plane family $\mathcal{H}_{\left(\widetilde{\varphi}, \widetilde{\nu}\right)}$ 
creates an envelope if and only if $\left({\bf r}(s)\cdot {\bf n}(s)\right)
=\alpha_3(s,u)$  
and $\alpha_1(s,u)=\alpha_2(s,u)$.   
Therefore, again by {\color{black}Theorem \ref{theorem1}, }
we have the following concrete expression of the created envelopes.    
\begin{eqnarray*}
{ } & { } & \widetilde{f}(s,u) \\ 
{ } & = & 
\widetilde{\omega}(s,u)+\widetilde{\gamma}(s)\widetilde{\nu}(s) \\ 
{ } & = & 
 \left({\bf r}(s)\cdot {\bf n}(s)\right){\bf n}(s)+\alpha(s,u){\bf t}(s) 
+\left({\bf r}(s)\cdot {\bf b}(s)\right){\bf b}(s) \\ 
{ } & = & \left({\bf r}(s)\cdot {\bf n}(s)\right){\bf n}(s)+ 
\left({\bf r}(s)\cdot {\bf t}(s)\right){\bf t}(s)+ 
\left(\alpha(s,u)-\left({\bf r}(s)\cdot {\bf t}(s)\right)\right){\bf t}(s) 
+\left({\bf r}(s)\cdot {\bf b}(s)\right){\bf b}(s) \\ 
{ } & = & {\bf r}(s)+\beta(s,u){\bf t}(s),    
\end{eqnarray*} 
where $\alpha(s,u)=\alpha_1(s,u)=\alpha_2(s,u)$ 
and $\beta(s,u)=\alpha(s,u)-\left({\bf r}(s)\cdot {\bf t}(s)\right)$.   
All envelopes created by the osculating family 
$\mathcal{H}_{\left(\widetilde{\varphi}, \widetilde{\nu}\right)}$ can be exactly 
expressed as above.   
Hence, for example, both the tangent developable of ${\bf r}$ 
(in the case $\beta(s,u)=u$) and the space curve ${\bf r}$ (in the case 
$\beta(s,u)=0$) are envelopes of 
$\mathcal{H}_{\left(\widetilde{\varphi}, \widetilde{\nu}\right)}$.   
Not only these two, there are uncountably many envelopes created by 
$\mathcal{H}_{\left(\widetilde{\varphi}, \widetilde{\nu}\right)}$.   
All envelopes for the osculating plane family are created only by 
the given curve ${\bf r}$ and its unit tangent curve ${\bf t}$.   
\par 
Next, we consider envelopes created by  
$\mathcal{H}_{\left({\bf r}, {\bf b}\right)}$ and 
$\mathcal{H}_{\left(\widetilde{f}, {\bf n}\right)}$.   
Namely, we obtain all solutions ${\color{black}\widetilde{g}}(s,u)$ 
for the following system of PDEs with 
one constraint condition.   
\[
\left\{
\begin{array}{ccc}
\frac{\partial {\color{black}\widetilde{g}}}{\partial s}(s,u)\cdot {\bf b}(s) & = & 0, \\ 
\frac{\partial {\color{black}\widetilde{g}}}{\partial u}(s,u)\cdot {\bf b}(s) & = & 0, \\ 
\frac{\partial {\color{black}\widetilde{g}}}{\partial s}(s,u)\cdot {\bf n}(s) & = & 0, \\ 
\frac{\partial {\color{black}\widetilde{g}}}{\partial u}(s,u)\cdot {\bf n}(s) & = & 0, \\ 
\left({\color{black}\widetilde{g}}(s, u)-{\bf r}(s)\right)\cdot {\bf b} (s)& = & 0.  
\end{array}
\right.   
\]
Since $\kappa(s)>0$ for any $s\in \mathbb{R}$ and 
\begin{eqnarray*}
\frac{\partial \widetilde{f}}{\partial s}(s,u) & = & {\bf t}(s)
+\frac{\partial \beta}{\partial s}(s,u){\bf t}(s)
+\beta(s,u)\left(\kappa(s){\bf n}(s)\right), 
 \\ 
\frac{\partial \widetilde{f}}{\partial u}(s,u)  
& = & \frac{\partial \beta}{\partial u}(s,u){\bf t}(s), 
\end{eqnarray*} 
if $\widetilde{f}$ {\color{black}itself} is a solution of the above system of PDEs, 
then $\beta(s,u)$ must be constant 
$0$.    Conversely, it is clear that ${\bf r}$ {\color{black}itself} is a solution 
of the above system of PDEs with one constraint condition.   
Therefore, for the above system of PDEs with one constraint condition, 
there are no solutions except for the trivial solution ${\bf r}$.    
This implies that even for a space curve ${\bf r}: \mathbb{R}\to \mathbb{R}^3$, 
the inverse Legendre transformation  
\[
\mathcal{H}_{\left({\bf r}, \{{\bf b}, {\bf n}\}\right)}\mapsto {\bf r}
\]
works well.   
\par 
Finally, we consider envelopes created by 
$\mathcal{H}_{\left({\bf r}, {\bf b}\right)}$ and 
$\mathcal{H}_{\left(\widetilde{f}, {\bf t}\right)}$.   
Namely, we obtain all solutions ${\color{black}\widetilde{g}}(s,u)$ 
for the following system of PDEs with 
one constraint condition.   
\[
\left\{
\begin{array}{ccc}
\frac{\partial {\color{black}\widetilde{g}}}{\partial s}(s,u)\cdot {\bf b}(s) & = & 0, \\ 
\frac{\partial {\color{black}\widetilde{g}}}{\partial u}(s,u)\cdot {\bf b}(s) & = & 0, \\ 
\frac{\partial {\color{black}\widetilde{g}}}{\partial s}(s,u)\cdot {\bf t}(s) & = & 0, \\ 
\frac{\partial {\color{black}\widetilde{g}}}{\partial u}(s,u)\cdot {\bf t}(s) & = & 0, \\ 
\left({\color{black}\widetilde{g}}(s, u)-{\bf r}(s)\right)\cdot {\bf b} (s)& = & 0.  
\end{array}
\right.   
\]
By the above calculations, 
if $\widetilde{f}$ is a solution of the above system of PDEs, 
then both $1+\frac{\partial \beta}{\partial s}(s,u)=0$ and 
$\frac{\partial \beta}{\partial u}(s,u)=0$ must be satisfied.    
It follows $\beta(s,u)=-s+c$ $(c\in \mathbb{R})$.     
It is easily seen that for any $c\in \mathbb{R}$, 
the space curve $s\mapsto {\bf r}(s)+(-s+c){\bf t}(s)$ is a solution 
of the above system of PDEs with one constraint condition.   .   
Thus,  in this case, the system of PDEs with one constraint condition 
has uncountably many solutions.      
\end{enumerate}
\end{example}
}}
\begin{example}\label{example2}    
\begin{enumerate} 
{\color{black}
\item[(1)] (The shoe surface : Example 1 of \cite{gauss})  
In this example, along the general theory developed in this paper, 
we start from making several general formulas for the envelope created 
by the affine tangent plane family of the surface having the form 
$\widetilde{\varphi}: \mathbb{R}^2\to \mathbb{R}^3$,\, 
$\widetilde{\varphi}(x, y)=(x,\, y,\, \widetilde{\varphi}_1(x,y))$ such that 
the origin $(0,0)$ is a singular point of the function 
$\widetilde{\varphi}_1: \mathbb{R}^2\to \mathbb{R}$ 
and there are no other singular points of $\widetilde{\varphi}_1$.   
Then, by calculating the obtained general formulas in the case of 
the shoe surface 
$\widetilde{\varphi}(x, y)=\left(x,\, y,\, \frac{1}{3}x^3-\frac{1}{2}y^2\right)$,   
just by calculations,  
we confirm that the concrete representation form of the envelope 
created by the affine tangent plane family of the shoe surface 
$\widetilde{\varphi}$
is actually the shoe surface itself. 
\par     
Let $\widetilde{\varphi}: \mathbb{R}^2\to \mathbb{R}^3$ be the 
mapping having the form 
$\widetilde{\varphi}(x, y)=(x, y, \widetilde{\varphi}_1(x,y))$,  
where the function $\widetilde{\varphi}_1: \mathbb{R}^2\to \mathbb{R}$ has 
a unique singularity at the origin, namely 
$\frac{\partial \widetilde{\varphi}_1}{\partial x}(0, 0)=
\frac{\partial \widetilde{\varphi}_1}{\partial y}(0, 0)=0$  
and $\left(\frac{\partial \widetilde{\varphi}_1}{\partial x}(x, y), \,
\frac{\partial \widetilde{\varphi}_1}{\partial y}(x, y)\right)\ne (0,0)$ 
for any $(x, y)\in \mathbb{R}^2-\{(0,0)\}$.  
Then, the mapping $\widetilde{\nu}: \mathbb{R}^2\to S^2$ 
defined by 
\[
\widetilde{\nu}(x,y)= 
\frac
{
\frac{\partial \widetilde{\varphi}_1}{\partial x}(x,y)\times 
\frac{\partial \widetilde{\varphi}_1}{\partial y}(x,y)
}
{
\parallel
\frac{\partial \widetilde{\varphi}_1}{\partial x}(x,y)\times 
\frac{\partial \widetilde{\varphi}_1}{\partial y}(x,y)
\parallel
} 
= 
\frac{\left(-\frac{\partial \widetilde{\varphi}_1}{\partial x},\,  
-\frac{\partial \widetilde{\varphi}_1}{\partial y},\, 1\right)}
{\sqrt{\left(\frac{\partial \widetilde{\varphi}_1}{\partial x}\right)^2+
\left(\frac{\partial \widetilde{\varphi}_1}{\partial y}\right)^2+1}}   
\]
is a Gauss mapping of the tangent plane family 
of $\widetilde{\varphi}$.     
Here, the tangent plane family of $\widetilde{\varphi}$ is 
$\mathcal{H}_{\left(\widetilde{\varphi}, \widetilde{\nu}\right)}$.
Let $(x_0, y_0)$ be an arbitrary point of $\mathbb{R}^2-\{(0,0)\}$.    
Then, by the assumption on the function $\widetilde{\varphi}_1$, 
it follows that $\widetilde{\nu}(x_0, y_0)\ne (0,\,0,\,1)$.    
Set 
\begin{eqnarray*}
{\bf v}_0(x_0, y_0) & = & \widetilde{\nu}(x_0, y_0), \\ 
{\bf v}_1(x_0, y_0) & = & 
\frac{
(0,\,0,\,1)- \left((0,\,0,\,1)\cdot {\bf v}_0(x_0, y_0)\right){\bf v}_0(x_0, y_0)
}
{
\parallel 
(0,\,0,\,1)- \left((0,\,0,\,1)\cdot {\bf v}_0(x_0, y_0)\right){\bf v}_0(x_0, y_0)
\parallel
}, \\ 
{\bf v}_2(x_0, y_0) & = & {\bf v}_0(x_0, y_0) \times {\bf v}_1(x_0, y_0).      
\end{eqnarray*}
Then, 
$\langle{\bf v}_0(x_0, y_0), {\bf v}_1(x_0, y_0), {\bf v}_2(x_0, y_0)\rangle$ 
is an orthonormal basis of $\mathbb{R}^3$, 
and under the identification of two vector spaces 
$\mathbb{R}^3$ and $T_{\widetilde{\nu}(x_0, y_0)}\mathbb{R}^3$, 
$\langle{\bf v}_1(x_0, y_0), {\bf v}_2(x_0, y_0)\rangle$ is an orthonormal 
basis of the tangent vector space $T_{\widetilde{\nu}(x_0, y_0)}S^2$.    
Let $\varepsilon$ be a sufficiently small positive number and denote 
the set $\{\Theta_1{\bf v}_1(x_0, y_0)+\Theta_2{\bf v}_2(x_0, y_0)
\, |\, -\varepsilon < \Theta_1, \Theta_2 < \varepsilon\}$ by $V'$.   
Let $\exp: V'\to S^2$ be the restriction of the exponential mapping at 
$\widetilde{\nu}(x_0, y_0)$ to $V'$ and set $V=\exp(V')$.    
Let 
$\left(V, \left(\Theta_1, \Theta_2\right)\right)$ be the normal 
coordinate neighborhood at 
$\widetilde{\nu}(x_0, y_0)$ defined by $\exp^{-1}: V\to V'$. 
Set 
\[
\widetilde{\gamma}(x,y)= 
\widetilde{\varphi}(x,y)\cdot \widetilde{\nu}(x,y)= 
\frac{-x\frac{\partial \widetilde{\varphi}_1}{\partial x}-y
\frac{\partial \widetilde{\varphi}_1}{\partial y}+\widetilde{\varphi}_1(x,y)
}{\sqrt{\left(\frac{\partial \widetilde{\varphi}_1}{\partial x}\right)^2+
\left(\frac{\partial \widetilde{\varphi}_1}{\partial y}\right)^2+1}}.   
\]  
Since $\widetilde{\nu}:\mathbb{R}^2\to S^2$ is a Gauss mapping of 
$\widetilde{\varphi}: \mathbb{R}^2\to \mathbb{R}^3$, we have    
\begin{eqnarray*}
{ } & { } & \frac{\partial \widetilde{\gamma}}{\partial x}(x_0, y_0) \\ 
{ } & = & 
\widetilde{\varphi}(x_0, y_0)\cdot 
\frac{\partial \widetilde{\nu}}{\partial x}(x_0, y_0) \\ 
{ } & = & 
\left(\widetilde{\varphi}(x_0, y_0)\cdot {\bf v}_1(x_0, y_0)\right)
\frac{\partial \left(\Theta_1\circ \widetilde{\nu}\right)}{\partial x}(x_0, y_0) 
+ 
\left(\widetilde{\varphi}(x_0, y_0)\cdot {\bf v}_2(x_0, y_0)\right)
\frac{\partial \left(\Theta_2\circ \widetilde{\nu}\right)}{\partial x}(x_0, y_0)  
\end{eqnarray*}
and 
\begin{eqnarray*}
{ } & { } & \frac{\partial \widetilde{\gamma}}{\partial y}(x_0, y_0) \\ 
{ } & = & 
\widetilde{\varphi}(x_0, y_0)\cdot 
\frac{\partial \widetilde{\nu}}{\partial y}(x_0, y_0) \\ 
{ } & = & 
\left(\widetilde{\varphi}(x_0, y_0)\cdot {\bf v}_1(x_0, y_0)\right)
\frac{\partial \left(\Theta_1\circ \widetilde{\nu}\right)}{\partial y}(x_0, y_0) 
+ 
\left(\widetilde{\varphi}(x_0, y_0)\cdot {\bf v}_2(x_0, y_0)\right)
\frac{\partial \left(\Theta_2\circ \widetilde{\nu}\right)}{\partial y}(x_0, y_0).     
\end{eqnarray*}
Thus, as the equality of $2$-dimensional 
cotangent vectors of $T^*_{(x_0, y_0)}\mathbb{R}^2$, 
we have the following equality.  
\begin{eqnarray*}
d\widetilde{\gamma} & = & 
\frac{\partial \widetilde{\gamma}}{\partial x}(x_0, y_0)dx + 
\frac{\partial \widetilde{\gamma}}{\partial y}(x_0, y_0)dy  \\ 
{ } & = & 
\left(\widetilde{\varphi}(x_0, y_0)\cdot {\bf v}_1(x_0, y_0)\right)
d\left(\Theta_1\circ \widetilde{\nu}\right) 
+ 
\left(\widetilde{\varphi}(x_0, y_0)\cdot {\bf v}_2(x_0, y_0)\right)
d\left(\Theta_2\circ \widetilde{\nu}\right).    
\end{eqnarray*} 
Set $U=\mathbb{R}^2-\{(0,0)\}$ and assume that the singular set of 
$\widetilde{\nu}$ is of Lebesgue measure zero.   
Then, since $(x_0, y_0)$ is an arbitrary point of $U$, 
by Theorem \ref{theorem1} (1) and Theorem \ref{theorem2}, 
it follows that $\mathcal{H}_{(\widetilde{\varphi}|_U, \widetilde{\nu}|_U)}$ 
creates a unique envelope.   
Set 
\[
\widetilde{\omega}(x_0, y_0) 
= 
\left(\widetilde{\varphi}(x_0, y_0)\cdot {\bf v}_1(x_0, y_0)\right)
d\left(\Theta_1\circ \widetilde{\nu}\right) 
+ 
\left(\widetilde{\varphi}(x_0, y_0)\cdot {\bf v}_2(x_0, y_0)\right)
d\left(\Theta_2\circ \widetilde{\nu}\right).   
\]
Then, under the canonical identifications 
\[
T^*_{\widetilde{\nu}(x_0, y_0)}S^2\cong T_{\widetilde{\nu}(x_0, y_0)}S^2 
\subset T_{\widetilde{\nu}(x_0, y_0)}\mathbb{R}^3 \cong \mathbb{R}^3, 
\]
the $2$-dimensional cotangent vector 
\[
\widetilde{\omega}(x_0, y_0) 
= 
\left(\widetilde{\varphi}(x_0, y_0)\cdot {\bf v}_1(x_0, y_0)\right)
d\Theta_1 
+ 
\left(\widetilde{\varphi}(x_0, y_0)\cdot {\bf v}_2(x_0, y_0)\right)
d\Theta_2 
\] may be regarded 
as the following $3$-dimensional vector (denoted by the same symbol 
$\widetilde{\omega}(x_0, y_0)$).   
\[
\widetilde{\omega}(x_0, y_0) 
= 
\left(\widetilde{\varphi}(x_0, y_0)\cdot {\bf v}_1(x_0, y_0)\right)
{\bf v}_1(x_0, y_0) 
+ 
\left(\widetilde{\varphi}(x_0, y_0)\cdot {\bf v}_2(x_0, y_0)\right)
{\bf v}_2(x_0, y_0).   
\]
Therefore, by Theorem \ref{theorem1} (2), the envelope vector at $(x_0, y_0)$ 
must have the following form:   
\begin{eqnarray*}
\widetilde{f}(x_0, y_0)  
 & = & \widetilde{\omega}(x_0, y_0)
+\widetilde{\gamma}(x_0, y_0)\widetilde{\nu}(x_0, y_0) \\ 
{ } & = & 
\left(\widetilde{\varphi}(x_0, y_0)\cdot {\bf v}_1(x_0, y_0)\right)
{\bf v}_1(x_0, y_0) 
+ 
\left(\widetilde{\varphi}(x_0, y_0)\cdot {\bf v}_2(x_0, y_0)\right)
{\bf v}_2(x_0, y_0)  \\ 
{ } & { } & 
\qquad\qquad +
\left(\widetilde{\varphi}(x_0, y_0)\cdot {\bf v}_0(x_0, y_0)\right)
{\bf v}_0(x_0, y_0) \\ 
{ } & = & \widetilde{\varphi}(x_0, y_0).     
\end{eqnarray*}
By continuity, it follows that $\widetilde{f}=\widetilde{\varphi}$ is the unique 
envelope created by the given plane family 
$\mathcal{H}_{\left(\widetilde{\varphi}, \widetilde{\nu}\right)}$.      
\par 
\medskip     
Next, we apply the above formulas to the shoe surface.    
The shoe surface is the image of 
$\widetilde{\varphi}: \mathbb{R}^2\to \mathbb{R}^3$ 
defined by $\widetilde{\varphi}(x, y)
=\left(x, y, \frac{1}{3}x^3 -\frac{1}{2}y^2\right)$.   
Set $\widetilde{\varphi}_1(x,y)=\frac{1}{3}x^3 -\frac{1}{2}y^2$.   
Then, the origin $(0,0)$ is a unique singular point of $\widetilde{\varphi}_1$.   
For the given $\widetilde{\varphi}$,  we have 
$\widetilde{\nu}(x,y)=
\frac{\frac{\partial \widetilde{\varphi}}{\partial x}(x,y)\times 
\frac{\partial \widetilde{\varphi}}{\partial y}(x,y)}
{||\frac{\partial \widetilde{\varphi}}{\partial x}(x,y)\times 
\frac{\partial \widetilde{\varphi}}{\partial y}(x,y)||}
= 
\frac{\left(-x^2,\, y,\, 1\right)}{\sqrt{x^4+y^2+1}}$.   
It is easily confirmed that the set consisting of regular points of 
$\widetilde{\nu}$ is dense.   In fact, it is known that 
any singularity of $\widetilde{\nu}$ is a fold 
singularity (see \cite{gauss}).  
Set $U=\mathbb{R}^2-\{(0,0)\}$ and take an arbitrary point $(x_0, y_0)$ of 
$U$.   
For the shoe surface $\widetilde{\varphi}$, we set 
\begin{eqnarray*}
{\bf v}_0(x_0, y_0) & = & \widetilde{\nu}(x_0, y_0)= 
\frac{\left(-x_0^2,\, y_0,\, 1\right)}{\sqrt{x_0^4+y_0^2+1}}, \\ 
{\bf v}_1(x_0, y_0) & = & 
\frac{
(0,\,0,\,1)- \left((0,\,0,\,1)\cdot {\bf v}_0(x_0, y_0)\right){\bf v}_0(x_0, y_0)
}
{
\parallel 
(0,\, 0,\, 1)- \left((0,\, 0,\, 1)\cdot {\bf v}_0(x_0, y_0)\right){\bf v}_0(x_0, y_0)
\parallel
}= 
\frac{\left(x_0^2,\, -y_0,\, x_0^4+y_0^2\right)}
{\sqrt{\left(x_0^4+y_0^2\right)\left(x_0^4+y_0^2+1\right)}}, \\ 
{\bf v}_2(x_0, y_0) & = & {\bf v}_0(x_0, y_0) \times {\bf v}_1(x_0, y_0)
=\frac{\left(y_0, \, x_0^2,\, 0\right)}{\sqrt{x_0^4+y_0^2}}.      
\end{eqnarray*}
By calculation, we have   
\begin{eqnarray*}
\widetilde{\varphi}(x_0, y_0)\cdot {\bf v}_1(x_0, y_0) 
& = & 
\frac{x_0^3-y_0^2+\left(\frac{1}{3}x_0^3-\frac{1}{2}y_0^2\right)
\left(x_0^4+y_0^2\right)}
{\left(x_0^4+y_0^2\right)^{\frac{1}{2}}\left(x_0^4+y_0^2+1\right)^{\frac{1}{2}}}, \\ 
\widetilde{\varphi}(x_0, y_0)\cdot {\bf v}_2(x_0, y_0) 
& = & 
\frac{x_0y_0+x_0^2y_0}
{\left(x_0^4+y_0^2\right)^{\frac{1}{2}}}.  
\end{eqnarray*}
Let $(V, \left(\Theta_1, \Theta_2\right))$ be the normal 
cordinate neighborhood of $S^2$ defined above.     
By calculations using the following two identities 
\begin{eqnarray*}
\frac{\partial \widetilde{\nu}}{\partial x}(x_0, y_0) 
& = & 
{\bf v}_1(x_0, y_0)
\frac{\partial \left(\Theta_1\circ \widetilde{\nu}\right)}{\partial x}(x_0, y_0)
+ 
{\bf v}_2(x_0, y_0)
\frac{\partial \left(\Theta_2\circ \widetilde{\nu}\right)}{\partial x}(x_0, y_0), \\ 
\frac{\partial \widetilde{\nu}}{\partial y}(x_0, y_0) 
& = & 
{\bf v}_1(x_0, y_0)
\frac{\partial \left(\Theta_1\circ \widetilde{\nu}\right)}{\partial y}(x_0, y_0)
+ 
{\bf v}_2(x_0, y_0)
\frac{\partial \left(\Theta_2\circ \widetilde{\nu}\right)}{\partial y}(x_0, y_0),  
\end{eqnarray*}   
we have the following.  
\begin{eqnarray*}
\frac{\partial \left(\Theta_1\circ \widetilde{\nu}\right)}{\partial x}
(x_0, y_0) 
& = & 
\frac{-2x_0^3}
{\left(x_0^4+y_0^2\right)^{\frac{1}{2}}\left(x_0^4+y_0^2+1\right)}, \\ 
\frac{\partial \left(\Theta_2\circ \widetilde{\nu}\right)}{\partial x}
(x_0, y_0) 
& = & 
\frac{-2x_0y_0-2x_0y_0^3-2x_0^5y_0}
{\left(x_0^4+y_0^2\right)^{\frac{1}{2}}\left(x_0^4+y_0^2+1\right)^{\frac{3}{2}}}, \\ 
\frac{\partial \left(\Theta_1\circ \widetilde{\nu}\right)}{\partial y}
(x_0, y_0) 
& = & 
\frac{-y_0}
{\left(x_0^4+y_0^2\right)^{\frac{1}{2}}\left(x_0^4+y_0^2+1\right)}, \\ 
\frac{\partial \left(\Theta_2\circ \widetilde{\nu}\right)}{\partial y}
(x_0, y_0) 
& = & 
\frac{x_0^2}
{\left(x_0^4+y_0^2\right)^{\frac{1}{2}}\left(x_0^4+y_0^2+1\right)^{\frac{1}{2}}}. \\ 
\end{eqnarray*}
On the other hand, 
from the form  
$\widetilde{\gamma}(x, y)=\widetilde{\varphi}(x,y)\cdot \widetilde{\nu}(x,y) 
=\frac{-\frac{2}{3}x^3+\frac{1}{2}y^2}{\sqrt{x^4+y^2+1}}$, 
we have   
\begin{eqnarray*}
\frac{\partial \gamma}{\partial x}\left(x_0, y_0\right) 
& = & 
\frac{-2x_0^2-2x_0^2y_0^2-x_0^3y_0^2-\frac{2}{3}x_0^6}
{\left(x_0^4+y_0^2+1\right)^{\frac{3}{2}}},  \\ 
\frac{\partial \gamma}{\partial y}\left(x_0, y_0\right) 
& = & 
\frac{y_0+\frac{1}{2}y_0^3+\frac{2}{3}x_0^3y_0+x_0^4y_0}
{\left(x_0^4+y_0^2+1\right)^{\frac{3}{2}}}.   
\end{eqnarray*}      
Thus, we have the following desired identity at $(x_0, y_0)$.   
{\small 
\begin{eqnarray*}
{ } & { } & d\widetilde{\gamma} \\ 
{ } & = & \frac{\partial \widetilde{\gamma}}{\partial x}(x_0, y_0)dx 
+ \frac{\partial \widetilde{\gamma}}{\partial y}(x_0, y_0)dy \\ 
{ } & = & 
\frac{-2x_0^2-2x_0^2y_0^2-x_0^3y_0^2-\frac{2}{3}x_0^6}
{\left(x_0^4+y_0^2+1\right)^{\frac{3}{2}}} dx 
+ 
\frac{y_0+\frac{1}{2}y_0^3+\frac{2}{3}x_0^3y_0+x_0^4y_0}
{\left(x_0^4+y_0^2+1\right)^{\frac{3}{2}}} dy \\ 
{ } & = & 
\left(
\frac{x_0^3-y_0^2+\left(\frac{1}{3}x_0^3-\frac{1}{2}y_0^2\right)
\left(x_0^4+y_0^2\right)}
{\left(x_0^4+y_0^2\right)^{\frac{1}{2}}\left(x_0^4+y_0^2+1\right)^{\frac{1}{2}}}
\frac{-2x_0^3}
{\left(x_0^4+y_0^2\right)^{\frac{1}{2}}\left(x_0^4+y_0^2+1\right)}
+
\frac{\left(x_0y_0+x_0^2y_0\right)}
{\left(x_0^4+y_0^2\right)^{\frac{1}{2}}}
\frac{\left(-2x_0y_0-2x_0y_0^3-2x_0^5y_0\right)}
{\left(x_0^4+y_0^2\right)^{\frac{1}{2}}\left(x_0^4+y_0^2+1\right)^{\frac{3}{2}}}
\right) dx \\ 
{ } & { } & \quad + 
\left(
\frac{x_0^3-y_0^2+\left(\frac{1}{3}x_0^3-\frac{1}{2}y_0^2\right)
\left(x_0^4+y_0^2\right)}
{\left(x_0^4+y_0^2\right)^{\frac{1}{2}}\left(x_0^4+y_0^2+1\right)^{\frac{1}{2}}}
\frac{-y_0}
{\left(x_0^4+y_0^2\right)^{\frac{1}{2}}\left(x_0^4+y_0^2+1\right)}
+ 
\frac{\left(x_0y_0+x_0^2y_0\right)}
{\left(x_0^4+y_0^2\right)^{\frac{1}{2}}}
\frac{x_0^2}
{\left(x_0^4+y_0^2\right)^{\frac{1}{2}}\left(x_0^4+y_0^2+1\right)^{\frac{1}{2}}}
\right) dy  \\ 
{ } & = & 
\left(
\left(\widetilde{\varphi}(x_0, y_0)\cdot {\bf v}_1(x_0, y_0)\right)
\frac{\partial \left(\Theta_1\circ \widetilde{\nu}\right)}{\partial x}(x_0, y_0)
+ 
\left(\widetilde{\varphi}(x_0, y_0)\cdot {\bf v}_2(x_0, y_0)\right)
\frac{\partial \left(\Theta_2\circ \widetilde{\nu}\right)}{\partial x}(x_0, y_0)
\right) dx \\ 
{ } & { } & \quad + 
\left(
\left(\widetilde{\varphi}(x_0, y_0)\cdot {\bf v}_1(x_0, y_0)\right)
\frac{\partial \left(\Theta_1\circ \widetilde{\nu}\right)}{\partial y}(x_0, y_0)
+ 
\left(\widetilde{\varphi}(x_0, y_0)\cdot {\bf v}_2(x_0, y_0)\right)
\frac{\partial \left(\Theta_2\circ \widetilde{\nu}\right)}{\partial y}(x_0, y_0)
\right) dy \\ 
{ } & = & 
\left(\widetilde{\varphi}(x_0, y_0)\cdot {\bf v}_1(x_0, y_0)\right) 
d\left(\Theta_1\circ \widetilde{\nu}\right) 
+ 
\left(\widetilde{\varphi}(x_0, y_0)\cdot {\bf v}_2(x_0, y_0)\right) 
d\left(\Theta_2\circ \widetilde{\nu}\right).     
\end{eqnarray*} 
}
Hence, by Theorem \ref{theorem1} (1) and Theorem \ref{theorem2}, 
the plane family 
$\mathcal{H}_{\left(\widetilde{\varphi}|_{U}, \widetilde{\nu}|_{U}\right)}$ 
for the shoe surface  
$\widetilde{\varphi}(x, y)=\left(x, y, \frac{1}{3}x^3-\frac{1}{2}y^2\right)$ 
has a unique envelope $\widetilde{f}: 
{\color{black}U} \to \mathbb{R}^3$, 
where $U=\mathbb{R}^2-\{(0,0)\}$.    
Then, under the canonical identifications 
\[
T^*_{\widetilde{\nu}(x_0, y_0)}S^2\cong T_{\widetilde{\nu}(x_0, y_0)}S^2 
\subset T_{\widetilde{\nu}(x_0, y_0)}\mathbb{R}^3 \cong \mathbb{R}^3, 
\]
the $2$-dimensional cotangent vector 
\[
\widetilde{\omega}(x_0, y_0) 
= 
\left(\widetilde{\varphi}(x_0, y_0)\cdot {\bf v}_1(x_0, y_0)\right)
d\Theta_1 
+ 
\left(\widetilde{\varphi}(x_0, y_0)\cdot {\bf v}_2(x_0, y_0)\right)
d\Theta_2 
\] is identified with  
the following $3$-dimensional vector (denoted by the same symbol 
$\widetilde{\omega}(x_0, y_0)$).  
\begin{eqnarray*}
\widetilde{\omega}(x_0, y_0) 
 & = &  
\left(\widetilde{\varphi}(x_0, y_0)\cdot {\bf v}_1(x_0, y_0)\right)
{\bf v}_1(x_0, y_0) 
+ 
\left(\widetilde{\varphi}(x_0, y_0)\cdot {\bf v}_2(x_0, y_0)\right)
{\bf v}_2(x_0, y_0) \\ 
{ } & = & 
\frac{\left(x_0^3-y_0^2+\left(\frac{1}{3}x_0^3-\frac{1}{2}y_0^2\right)
\left(x_0^4+y_0^2\right)\right)}
{\left(x_0^4+y_0^2\right)^{\frac{1}{2}}\left(x_0^4+y_0^2+1\right)^{\frac{1}{2}}}
\frac{\left(x_0^2,\, -y_0,\, x_0^4+y_0^2\right)}
{{\left(x_0^4+y_0^2\right)^{\frac{1}{2}}\left(x_0^4+y_0^2+1\right)^{\frac{1}{2}}}} \\ 
{ } & { } & \qquad 
+ 
\frac{\left(x_0y_0+x_0^2y_0\right)}
{\left(x_0^4+y_0^2\right)^{\frac{1}{2}}}
\frac{\left(y_0,\, x_0^2,\, 0\right)}
{\left(x_0^4+y_0^2\right)^{\frac{1}{2}}}. 
\end{eqnarray*}
Therefore, by Theorem \ref{theorem1} (2), 
the unique envelope $\widetilde{f}$ 
must have the following desired parametric representation 
on $U=\mathbb{R}^2-\{(0,0)\}$.   
\begin{eqnarray*}
{ } & { } & \widetilde{f}(x_0, y_0) \\
{ } & = & 
\widetilde{\omega}(x_0,y_0)
+\widetilde{\gamma}(x_0,y_0)\widetilde{\nu}(x_0,y_0) \\ 
{ } & = & 
\frac{\left(x_0^3-y_0^2+\left(\frac{1}{3}x_0^3-\frac{1}{2}y_0^2\right)
\left(x_0^4+y_0^2\right)\right)}
{\left(x_0^4+y_0^2\right)^{\frac{1}{2}}\left(x_0^4+y_0^2+1\right)^{\frac{1}{2}}}
\frac{\left(x_0^2,\, -y_0,\, x_0^4+y_0^2\right)}
{{\left(x_0^4+y_0^2\right)^{\frac{1}{2}}\left(x_0^4+y_0^2+1\right)^{\frac{1}{2}}}} \\ 
{ } & { } & \qquad 
+ 
\frac{\left(x_0y_0+x_0^2y_0\right)}
{\left(x_0^4+y_0^2\right)^{\frac{1}{2}}}
\frac{\left(y_0,\, x_0^2,\, 0\right)}
{\left(x_0^4+y_0^2\right)^{\frac{1}{2}}} 
+ 
\frac{\left(-\frac{2}{3}x_0^3+\frac{1}{2}y_0^2\right)}
{\left(x_0^4+y_0^2+1\right)^{\frac{1}{2}}}
\frac{\left(-x_0^2,\, y_0,\, 1\right)}{\left(x_0^4+y_0^2+1\right)^{\frac{1}{2}}} \\ 
{ } & = & 
\left(x_0,\, y_0,\, \frac{1}{3}x_0^3-\frac{1}{2}y_0^2\right) \\ 
{ } & = & 
\widetilde{\varphi}(x_0, y_0).   
\end{eqnarray*}
By continuity, it follows that the given shoe surface 
$\widetilde{\varphi}$ itself is the unique envelope created by 
the tangent plane family 
$\mathcal{H}_{\left(\widetilde{\varphi}, \widetilde{\nu}\right)}$.   
\par 
\smallskip 
The set called the {\it parabolic line} of 
$\widetilde{\varphi}: \mathbb{R}^2\to \mathbb{R}^3$ 
consists of points $(x,y)\in \mathbb{R}^2$ at which 
$\widetilde{\nu}$ is singular.      For the shoe surface, 
the parabolic line is the $y$-axis $\{(0,y)\, |\, y\in \mathbb{R}\}$.      
Thus, as similar as the case of unit speed plane curves  
${\bf r}: \mathbb{R}\to \mathbb{R}^2$ with inflection points,  
the full discriminant of the tangent plane family 
$\mathcal{H}_{\left(\widetilde{\varphi}, \widetilde{\nu}\right)}$ for the 
shoe surface $\widetilde{\varphi}: \mathbb{R}^2\to \mathbb{R}^3$ 
is different from the unique desired envelope 
$\widetilde{\varphi}$ itself, 
since the full discriminant includes 
{\color{black}an affine tangent line 
$\left\{\left.\left(\lambda, y, -\frac{1}{2}y^2\right)\; 
\right|\; \lambda\in \mathbb{R}\right\}$} 
at any point $(0, y)$.     
Therefore, even in the case of surfaces in $\mathbb{R}^3$, 
by our method, one can distinguish the envelope 
in the sense of Definition \ref{definition1} and the full discriminant.    
This means that, in the case of surfaces in $\mathbb{R}^3$ as well, 
our method has an advantage.   }
\end{enumerate}
\begin{enumerate}
\item[(2)] (Example 4.1 of \cite{izumiya})  
Let $\widetilde{\nu}: \mathbb{R}^n\to S^n\subset \mathbb{R}^{n+1}$ be 
the mapping defined by 
$\widetilde{\nu}\left(p_1, \ldots, p_n\right)=
\frac{1}{\sqrt{\sum_{i=1}^n p_i^2+1}}\left(p_1, \ldots, p_n, -1\right)$.     
Then, $\widetilde{\nu}$ is non-singular and its inverse mapping 
$\widetilde{\nu}^{-1}: \widetilde{\nu}\left(\mathbb{R}^{n+1}\right)
\to \mathbb{R}^{n+1}$ is the central projection relative to the south pole 
$(0, \ldots, 0, -1)$ of $S^n$.    
Let $\widetilde{\varphi}: \mathbb{R}^n\to \mathbb{R}^{n+1}$ 
be an arbitrary mapping.   
Set 
$\widetilde{\gamma}(p)=\widetilde{\varphi}(p)\cdot \widetilde{\nu}(p)$ where 
$p=\left(p_1, \ldots, p_n\right)$ be a point of $\mathbb{R}^{n+1}$.   
Let $\left(X=\left(X_1, \ldots, X_n\right), Y\right)$ be a point of 
$\mathbb{R}^n\times\mathbb{R}$.    Since 
$J^1(\mathbb{R}^n, \mathbb{R})$ and 
$\mathbb{R}^n\times \mathbb{R}\times\mathbb{R}^n$ are identified, 
$(X, Y, p)$ may be 
regarded as the canonical coordinate system of 
$J^1\left(\mathbb{R}^n, \mathbb{R}\right)$.   
Since $\frac{X_i\circ\widetilde{\nu}(p)}{Y\circ \widetilde{\nu}(p)}=-p_i$ 
for any $i$ $(1\le i\le n)$ and any $p\in \mathbb{R}^{n+1}$, 
considering the first order differential equation 
\[
\left(\left(X, Y\right)-\widetilde{\varphi}(p)\right)\cdot \widetilde{\nu}(p)=0 
\]
is exactly the same as considering the following Clairaut equation 
\[
Y=\sum_{i=1}^n X_i p_i+\frac{\widetilde{\varphi}(p)\cdot \widetilde{\nu}(p)}
{Y\circ \widetilde{\nu}(p)}.
\]
Thus, for each $x\in \mathbb{R}^{n}$ the hyperplane 
$H_{\left(\widetilde{\varphi}(x), \widetilde{\nu}(x)\right)}$ 
is a complete solution of the above Clairaut equation.     
Since $\widetilde{\nu}$ is non-singular, by Theorem \ref{theorem1} 
and Theorem \ref{theorem2}, the above Clairaut equation 
has a unique singular solution $\widetilde{f}: \mathbb{R}^n\to \mathbb{R}^{n+1}$.    
By {\color{black}Theorem \ref{theorem1} again, }
the unique singular solution $\widetilde{f}$ 
has the following expression 
where $x$ is an arbitrary point of $\mathbb{R}^n$ and 
$\left(V, \left(\Theta_1, \ldots, \Theta_n\right)\right)$ is 
a sufficiently small normal coordinate neighborhood of $S^n$ 
at $\widetilde{\nu}(x)$.   
\[
\widetilde{f}(x)=
\sum_{i=1}\frac{\partial \left(\widetilde{\gamma}\circ \widetilde{\nu}^{-1}\right)}
{\partial \Theta_{i}}\left(\widetilde{\nu}(x)\right)
\frac{\partial}{\partial \Theta_{i}}
+
\widetilde{\gamma}(x)\widetilde{\nu}(x).   
\]
By this expression, for instance, it is easily seen that 
when $\widetilde{\gamma}(x)\equiv c(\ne 0)$ 
for any $x\in \mathbb{R}^{n+1}$, then the unique singular solution 
$Y: U_c\to \mathbb{R}$ 
must be an explicit solution 
with the following expression where $U_c=\{X\; |\; \Vert X\Vert<|c|\}$.    
\[
Y(X)=
\left\{
\begin{array}{rr}
-\sqrt{|c|^2-\sum_{i=1}^n X_i^2} & (\mbox{ if } c>0)  \\ 
\sqrt{|c|^2-\sum_{i=1}^n X_i^2} & (\mbox{ if } c<0).
\end{array}
\right.
\]
\end{enumerate}
\end{example}
{\color{black}
\section{Appendix: Alternative proof of Theorem \ref{theorem1} except for 
the assertion (3) \\ in the case $n=1$}\label{section5} 
Let $N$ be a $1$-dimensional manfold and let 
$\widetilde{\varphi}: N\to \mathbb{R}^2$, 
$\widetilde{\nu}: N\to S^1$ be mappings.   
Define the function $\widetilde{\Theta}: N\to \mathbb{R}$ by 
$\widetilde{\nu}(t)=\left(\cos\widetilde{\Theta}(t), 
\sin\widetilde{\Theta}(t)\right)$.     
Define also $\widetilde{\tau}(t)
:=\left(\sin\widetilde{\Theta}(t), -\cos\widetilde{\Theta}(t)\right)$.   
Then, the following trivially holds.  
\begin{fact}
For any $h: N\to \mathbb{R}^2$, 
\[
h(t)=\left(h(t)\cdot \widetilde{\tau}(t)\right)\widetilde{\tau}(t) + 
\left(h(t)\cdot \widetilde{\nu}(t)\right)\widetilde{\nu}(t).   
\]
\end{fact}
\par 
We first show that the creative condition can be naturally obtained from 
an envelope by introducing a gauge theoretic approach.    
Suppose that 
$\widetilde{f}: N\to \mathbb{R}^2$ is an envelope created by 
the line family 
$\mathcal{H}_{\left(\widetilde{\varphi}, \widetilde{\nu}\right)}$.    
Then, we have the following.    
\[
\widetilde{\gamma}'(t)  =  
\left(\widetilde{f}(t)\cdot \widetilde{\nu}(t)\right)'  
 =  
\widetilde{f}'(t)\cdot \widetilde{\nu}(t) + \widetilde{f}(t)\cdot \widetilde{\nu}'(t) 
 =  
0 -\left(\widetilde{f}(t)\cdot \widetilde{\tau}(t)\right)\widetilde{\Theta}'(t) .
\]
Let $h: N\to N$ be a bijective mapping.   
Then, notice that 
\[
\mathcal{H}_{\left(\widetilde{\varphi}, \widetilde{\nu}\right)}
= 
\mathcal{H}_{\left(\widetilde{\varphi}\circ h, \widetilde{\nu}\circ h\right)}
\] 
and 
\[
\left(\widetilde{\gamma}\circ h\right)'(t)= 
-\left(\widetilde{f}(h(t))\cdot \widetilde{\tau}(h(t))\right)
\widetilde{\Theta}'(h(t))h'(t).  
\] 
From these simple observations, 
we see that it is important to extract a significant quantity 
which does not depend on the particular choice of $h$.   
Then, we naturally reach the following setting.   
\[
\widetilde{\omega}(t):=
-\left(\widetilde{f}(t)\cdot \widetilde{\tau}(t)\right)d\widetilde{\Theta}.   
\]   
and we trivially have $d\widetilde{\gamma}=\widetilde{\omega}$.    
Take an arbitrary point $t_0$ of $N$ and fix it.   
Let $(V, \Theta)$ be a normal coordinate neighborhood of $S^1$ at 
$\widetilde{\nu}(t_0)$ such that $\Theta\left(\widetilde{\nu}(t_0)\right)=0$ and  
$\widetilde{\Theta}(t)=
\left(\Theta\circ \widetilde{\nu}\right)(t)$ for any $t\in \widetilde{\nu}^{-1}(V)$.    
In other words, $\left(\Theta\circ \widetilde{\nu}\right)(t)$ 
$\left(t\in \widetilde{\nu}^{-1}(V)\right)$
is just the radian (or its negative) between two unit vectors 
$\widetilde{\nu}(t_0)$ and $\widetilde{\nu}(t)$.    
By using the function $\Theta: V\to \mathbb{R}$, 
the $1$-form  
$\widetilde{\omega}(t)$ may be written as follows.   
\[
\widetilde{\omega}(t)=
-\left(\widetilde{f}(t)\cdot \widetilde{\tau}(t)\right)
\widetilde{\nu}^*d\Theta, 
\]  
where $\widetilde{\nu}^*d{\Theta}$ stands for the pullback of the 
$1$-form $d\Theta$ by $\widetilde{\nu}$.    
Hence, we naturally reach the following $1$-form which is denoted by 
the same symbol $\widetilde{\omega}$.
\[
\widetilde{\omega}(t)=
-\left(\widetilde{f}(t)\cdot \widetilde{\tau}(t)\right)d\Theta.   
\]   
It is easily seen that for any $t\in \widetilde{\nu}^{-1}(V)$, 
under the canonical identifications 
\[
T^*_{\widetilde{\nu}(t)}S^1\cong T_{\widetilde{\nu}(t)}S^1 
\subset T_{\widetilde{\nu}(t)}\mathbb{R}^2\cong \mathbb{R}^2, 
\]
the $1$-dimensional cotangent vector 
\[
\widetilde{\omega}(t)=
-\left(\widetilde{f}(t)\cdot \widetilde{\tau}(t)\right)d\Theta 
\in T^*_{\widetilde{\nu}(t)}S^1
\]
is identified with the $2$-dimensional vector 
\[
\widetilde{\omega}(t)=
\left(\widetilde{f}(t)\cdot \widetilde{\tau}(t)\right)\widetilde{\tau}(t) 
\in \mathbb{R}^2.   
\]
Since $t_0$ is an arbitrary point of $N$, 
we naturally see that the creative condition is satisfied for 
$\mathcal{H}_{\left(\widetilde{\varphi, \widetilde{\nu}}\right)}$ and 
the following horizontal-vertical decomposition formula holds for any 
$t\in N$.      
\begin{fact}
\[
\widetilde{f}(t)=
\left(\widetilde{f}(t)\cdot \widetilde{\tau}(t)\right)\widetilde{\tau}(t) + 
\left(\widetilde{f}(t)\cdot \widetilde{\nu}(t)\right)\widetilde{\nu}(t)  
=\widetilde{\omega}(t)+\widetilde{\gamma}(t)\widetilde{\nu}(t).     
\]
\end{fact}
\smallskip 
\par 
Conversely, suppose that 
$\mathcal{H}_{\left(\widetilde{\varphi}, \widetilde{\nu}\right)}$ is creative.   
Then, there exists a function $\alpha: N\to \mathbb{R}$ such that 
$d\widetilde{\gamma}=\alpha d\widetilde{\Theta}$.     
Set $\widetilde{\omega}=\alpha d\widetilde{\Theta}$.   
Let $t_0\in N$ be an arbitrary point.    
Then, under the canonical identifications 
\[
T^*_{\widetilde{\nu}(t)}S^1\cong T_{\widetilde{\nu}(t)}S^1 
\subset T_{\widetilde{\nu}(t)}\mathbb{R}^2\cong \mathbb{R}^2, 
\]
the $1$-dimensional cotangent vector 
\[
\widetilde{\omega}(t)=
\alpha(t)d\Theta 
\in T^*_{\widetilde{\nu}(t)}S^1
\]
is identified with the $2$-dimensional vector 
\[
\widetilde{\omega}(t)=
-\alpha(t)\widetilde{\tau}(t) 
\in \mathbb{R}^2, 
\]
where $(V, \Theta)$ is a normal coordinate system 
of $S^1$ at $\widetilde{\nu}(t_0)$ such that $\Theta(\widetilde{\nu}(t_0))=0$ 
and $t\in \widetilde{\nu}^{-1}(V)$.    
Set 
\[
\widetilde{f}(t)=\widetilde{\omega}(t)
+\widetilde{\gamma}(t)\widetilde{\nu}(t) 
= 
-\alpha(t)\widetilde{\tau}(t)+\widetilde{\gamma}(t)\widetilde{\nu}(t).   
\]
Then, 
$\widetilde{f}$ clearly satisfies the condition (b) of Definition 
\ref{definition1} for any $t\in \widetilde{\nu}^{-1}(V)$.   
Moreover we have the following.   
\begin{lemma}\label{lemma5.1}
For any $t\in \widetilde{\nu}^{-1}(V)$, 
$\widetilde{f}'(t)\cdot \widetilde{\nu}(t)=0$ holds.   
\end{lemma}  
\noindent 
\underline{\it Proof of Lemma \ref{lemma5.1}}\quad 
We have 
\[
\widetilde{\gamma}'(t) = 
\left(\widetilde{f}(t)\cdot \widetilde{\nu}(t)\right)'
= 
\widetilde{f}'(t)\cdot \widetilde{\nu}(t)
-\left(\widetilde{f}{\color{black}(t)}
\cdot \widetilde{\tau}(t)\right)\widetilde{\Theta}'(t) 
= 
\widetilde{f}'(t)\cdot \widetilde{\nu}(t)
+\alpha(t)\widetilde{\Theta}'(t).    
\]
Thus, we have the following.   
\begin{eqnarray*}
\widetilde{\omega}(t)=d\widetilde{\gamma}=\widetilde{\gamma}'(t)dt 
& = & 
\left(\widetilde{f}'(t)\cdot \widetilde{\nu}(t)\right)dt 
+\alpha(t)\widetilde{\Theta}'(t)dt \\ 
{ } & = & 
\left(\widetilde{f}'(t)\cdot \widetilde{\nu}(t)\right)dt 
+\alpha(t)d\widetilde{\Theta} \\ 
{ } & = & 
\left(\widetilde{f}'(t)\cdot \widetilde{\nu}(t)\right)dt 
+\widetilde{\omega}(t).   
\end{eqnarray*} 
It follows $\left(\widetilde{f}'(t)\cdot \widetilde{\nu}(t)\right)dt =0$.    
Since $t$ is a coordinate function on an open set 
$\widetilde{\nu}^{-1}(V)$ of $N$, 
for any fixed $t\in \widetilde{\nu}^{-1}(V)$, 
the $1$-dimensional cotangent vector 
$dt$ at $t$ is not zero.  
Therefore, the number 
$\left(\widetilde{f}'(t)\cdot \widetilde{\nu}(t)\right)
$ is always zero for any $t\in \widetilde{\nu}^{-1}(V)$.    
Since $t_0$ is an arbitrary point of $N$, 
Theorem \ref{theorem1} (1) holds.    
By the above decomposition of $\widetilde{f}(t)$, 
Theorem \ref{theorem1} (2) holds as well.    
\qed    
}

\textcolor{black}
{
\section*{Acknowledgements}
{\color{black}The author is most grateful to two anonymous reviewers.     
Their comments/suggestions are hard to replace.   }               
The author would like to thank Richard Montgomery for appropriate comments.     
His suggestions improved this paper.         
}
%
%

\end{document}